\documentclass[11pt,a4paper]{article}
\usepackage[utf8]{inputenc}
\usepackage{appendix}
\usepackage{subfig}
\usepackage{graphicx}
\usepackage{ifthen}  

\usepackage[most]{tcolorbox}
\usepackage{tikz,tikz-cd}
\usetikzlibrary{decorations.pathmorphing, decorations.pathreplacing, patterns, angles, quotes}

\usepackage{xspace}
\usepackage{comment}

\usepackage{amsmath,amssymb,amsfonts}
\usepackage{bm}
\usepackage{amsthm}

\usepackage{multirow}
\usepackage{derivative}
\usepackage{xcolor}
\usepackage{pifont}

\usepackage{hyperref}

\usepackage{authblk}
\usepackage[a4paper, margin=0.8in]{geometry} 
\usepackage[numbers]{natbib}
\newcommand{\greencheck}{{\color{green}\ding{51}}} 
\newcommand{\redcross}{{\color{red}\ding{55}}} 

\newtheorem{remark}{Remark}
\newtheorem{proposition}{Proposition}

\renewcommand\d{\ensuremath{\mathrm{d}}}

\newcommand{\bbR}{\mathbb{R}}
\newcommand{\bbV}{\mathbb{V}}
\newcommand{\bbS}{\mathbb{S}}
\newcommand{\bbM}{\mathbb{M}}

\newcommand{\bbP}{\mathbb{P}}

\DeclareMathOperator{\tr}{tr}

\newcommand*{\norm}[1]{\ensuremath{\left\|#1\right\|}}

\newcommand{\inpr}[3][]{\ensuremath{( #2, \, #3 )_{#1}}}

\DeclareMathOperator*{\grad}{grad}
\renewcommand{\div}{\operatorname{div}}

\DeclareMathOperator*{\Hess}{Hess}
\DeclareMathOperator*{\curl}{curl}
\DeclareMathOperator*{\rot}{rot}

\DeclareMathOperator*{\Air}{Air}

\DeclareMathOperator*{\sgrad}{def}
\DeclareMathOperator*{\sym}{sym}

\newcommand{\firedrake}{\textsc{Firedrake}\xspace}

\newcommand{\rewone}[1]{\textcolor{black}{#1}}
\newcommand{\rewtwo}[1]{\textcolor{black}{#1}}
\newcommand{\rewother}[1]{\textcolor{black}{#1}}

\newcommand{\keywords}[1]{\textbf{Keywords:}\hspace*{5pt}#1}

\title{A linearly-implicit energy-momentum preserving scheme for geometrically nonlinear mechanics based on non-canonical Hamiltonian formulations}
\date{}
\author[1]{Andrea Brugnoli}
\author[2]{Denis Matignon}
\author[1]{Joseph Morlier}

\affil[1]{ICA, Universit\'e de Toulouse, ISAE–SUPAERO, INSA, CNRS, MINES ALBI, UPS, Toulouse, France}
\affil[2]{F\'ed\'eration ENAC ISAE-SUPAERO ONERA, Universit\'e de Toulouse, Toulouse, France}

\begin{document}

\maketitle


\abstract{
This work presents a novel formulation and numerical strategy for the simulation of geometrically nonlinear structures. First, a non-canonical Hamiltonian (Poisson) formulation is introduced by including the dynamics of the stress tensor. This framework is developed for von-Kármán nonlinearities in beams and plates, as well as geometrically nonlinear elasticity with Saint-Venant material behavior. In the case of plates, both negligible and non-negligible membrane inertia are considered. For the former case the two-dimensional elasticity complex is leveraged to express the dynamics in terms of the Airy stress function. The finite element discretization employs a mixed approach, combining a conforming approximation for displacement and velocity fields with a discontinuous stress tensor representation. A staggered, linear implicit time integration scheme is proposed, establishing connections with existing explicit-implicit energy-preserving methods. The stress degrees of freedom are statically condensed, reducing the computational complexity to solving a system with a positive definite matrix. \rewtwo{The integration strategy preserves energy and angular momentum exactly.}  The methodology is validated through numerical experiments on the Duffing oscillator, a von-Kármán beam, and a column undergoing finite deformations. Comparisons with fully implicit energy-preserving method and the leapfrog scheme demonstrate that the proposed approach achieves superior accuracy while maintaining energy stability. Additionally, it enables larger time steps compared to explicit schemes and exhibits computational efficiency comparable to the leapfrog method. 
}

\vspace{10pt}
\keywords{geometrically-nonlinear mechanics, geometric numerical integration, exact energy conservation, Hamiltonian dynamics, mixed finite elements}

\section{Introduction}
In mechanics, nonlinearities may arise from material behavior, contact and friction phenomena, or large deformations and displacements, in which case they are referred to as geometrical nonlinearities \cite{lacarbonara2013nonlinear}. These nonlinearities occur in many real-world engineering applications \cite{touzé2021review} such as aeronautics \cite{patil2004importance,kerschen2013nonlinear}, wind energy systems \cite{manolas2015importance}, musical acoustics \cite{chaigne2005gongs,jossic2018internal} or microelectromechanical devices \cite{lazarus2012fem,gobat2021rom}. The accurate time-domain simulation of geometrically nonlinear systems is an essential tool for their analysis and design, widely employed across various disciplines, including computer animation \cite{bertails2006super}, sound synthesis \cite{bilbao2009numerical}, model reduction \cite{vizzaccaro2021rom,jain2022compute} and control \cite{artola2021modal}. \\

When dealing with flexible structures, spatial discretization is required to obtain a system of ordinary differential equations (ODEs). In the context of musical acoustics finite difference methods are typically used \cite{bilbao2001phdthesis}, because they yield diagonal mass matrices after spatial discretization and this feature can be exploited to develop fast time integration schemes \cite{ducceschi2022string}. However, implementing finite difference methods for complex geometrical domains is challenging and their application to general three-dimensional elasticity is limited to the linear case \cite{lipnikov2014mimetic}. For nonlinear elasticity problems, the finite element method remains the most widely used approach \cite{deborst2012nonlinear}. \\

The time integration of the resulting ODEs is typically performed using finite difference methods, with the Newmark method being the most well-known \cite{newmark1959}. Its success is attributed to its variational properties, meaning that it originates from discrete variational principles and thus inherits many desirable properties, such as symplecticity and momentum preservation \cite{kane2000variational}. A specific instance of the Newmark method, the implicit midpoint method, also preserves energy when applied to linear systems; however, this property does not extend to the nonlinear case. To achieve exact energy preservation in geometrically nonlinear problems, the stress tensor should not be evaluated at the midpoint but instead computed as an average between the current and next time instants, as shown in \cite{simo1992conserving}. This approach is a special case of the discrete gradient method \cite{mclachlan1999discrete,franke2023} and leads to an implicit scheme that requires a root-finding method, such as Newton's method. \\

In recent years, exact energy-preserving schemes that require solving only a linear system have been developed. The first of these methods, known as invariant energy quadratisation (IEQ) approaches, was initially applied to phase-field models \cite{yang2016phase}. The method was further developed and simplified in \cite{shen2018sav}, where a spatially distributed field was replaced by a single scalar variable. For this reason, it is referred to as the scalar auxiliary variable (SAV) method. This framework was later extended to Hamiltonian systems \cite{bilbao2023explicit}, particularly in the context of mechanical models in acoustics. The authors consider a finite difference discretization in space and diagonal mass matrices. Thanks to this assumption, the scheme can be made fully explicit using the Sherman-Morrison inversion theorem \cite{sherman1950adjustment}. The methodology was further enhanced to higher-order accuracy and extended to systems with multiple first integrals in \cite{andrews2024high}, utilizing a finite element variational discretization in time. However, this strategy is implicit and therefore computationally demanding. \\

This work builds upon and extends the preliminary ideas introduced in \cite{brugnoli2024exact}, serving a dual purpose. First, geometrically nonlinear models are formulated in a non-canonical Hamiltonian form (also known as the Poisson formulation) by considering the dynamics of the stress tensor. Second, a discretization strategy employing mixed finite elements and a linearly implicit time integration scheme is presented. The Poisson formulation is detailed for von Kármán-type nonlinearities in beams and plates (preliminary results can be found in \cite{brugnoli2021vonkarman,brugnoli2022enoc}), as well as for geometrically nonlinear elasticity with Saint-Venant material behavior. For plates, both the cases of negligible and non-negligible membrane inertia are considered. In the former case, the two-dimensional elasticity complex—first developed by Kröner in the context of linear-elastic dislocation theory \cite{kroner1959}—is leveraged to express the dynamics in terms of the Airy stress function \cite{touzé2021review}. The Poisson formulation for geometrically nonlinear elasticity has already been detailed in \cite{thomas2024velocity} within the framework of port-Hamiltonian systems. \rewother{An analogous Poisson formulation has been used in \cite{kinon2023discrete,kinon2025energy} for the specific case of systems of geometrically nonlinear oscillators and planar geometrically exact beams. In the latter case, their findings suggest that the employed formulation is advantageous because:
    \begin{enumerate}
        \item It does not suffer from shear locking as the shear stress is discretized using a discontinuous finite element;
        \item Even if the formulation is dynamical, it can be used as an iterative method to solve quasi static problem when the inertia is set to zero;
        \item It requires less Newton iterations than other well-established schemes used in the literature.
    \end{enumerate}}

The finite element discretization employs a conforming approximation of the displacement and velocity fields, along with a discontinuous approximation of the stress tensor. This discontinuous space must be carefully chosen to accurately capture geometrically nonlinear effects. The time integration scheme is a staggered version of the method proposed in \cite{bilbao2023explicit}, and the connection between these two approaches is established. \rewtwo{The scheme is shown to preserve exactly energy and angular momentum.} The degrees of freedom associated with the stress variable can be statically condensed at the discrete time  level, reducing the problem to solving a system involving a positive definite matrix to compute the velocity field at the next time step. \\

The methodology is tested on the Duffing oscillator, the vibrations of a von Kármán beam, and the bending of a column in geometrically nonlinear elasticity. It is compared against the exact energy-preserving method by Simo \cite{simo1992conserving}, that corresponds to a discrete gradient method, and the explicit central difference Newmark method (also known as the St\"ormer-Verlet or leapfrog scheme). The proposed method exhibits higher accuracy when measuring the error relative to the exact solution of the Duffing oscillator. Furthermore, it is energy stable, allowing for a larger time step than the explicit Newmark method. When applied to flexible structures, the scheme demonstrates computational efficiency comparable to the leapfrog method. However, mass lumping strategies —~which would significantly enhance the efficiency of the leapfrog method~— are not considered in this work. Nonetheless, the present approach can be further optimized through appropriate numerical linear algebra techniques when mass lumping is incorporated. Based on the numerical experiments, Table \ref{tab:comp_methods} highlights the advantages and limitations of each method.  \\

The paper is structured as follows. Section~\ref{sec:models} details the non-canonical Hamiltonian formulation, first  for a finite-dimensional nonlinear oscillator, then in a more abstract framework that is specialized to the aforementioned models. Then,  the space and time discretizations are addressed in Sec.~\ref{sec:discr}. Numerical results are presented in Sec.~\ref{sec:num_results}. Finally, perspectives are given and conclusions are drawn in Sec.~\ref{sec:concl}.

\begin{table}[htb]
\centering
\begin{tabular}{lccc}
\hline
Method & Accuracy & Stability & Efficiency \\
\hline        
Discrete Gradient & \greencheck & \greencheck & \redcross \\
Linear Implicit & \greencheck & \greencheck & \greencheck \\
Leapfrog & \greencheck & \redcross & \greencheck \\
\hline
\end{tabular}
    \caption{Comparison of time integration methods for geometrically nonlinear mechanics}
    \label{tab:comp_methods}
\end{table}

\section{A non-canonical Hamiltonian structure for geometrical nonlinear mechanics}
\label{sec:models}
Geometrically nonlinear problems in continuum mechanics can be framed in a general non-canonical Hamiltonian structure. This statement will be motivated through several examples in mechanics. For the sake of clarity, a one degree-of-freedom geometrical nonlinear oscillator is first considered. Then, an abstract formulation is introduced. Finally, examples from continuum mechanics justifying the abstract formulation are considered.

\subsection{An introductory finite-dimensional example}

  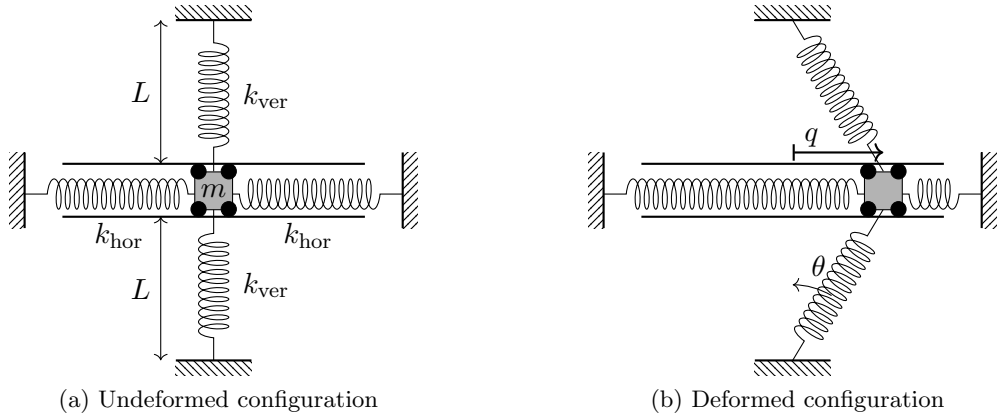
\begin{figure*}[tbh]
\centering
\subfloat[Undeformed configuration]{%
\begin{tikzpicture}
\fill [pattern = north west lines] (-.5,0) rectangle ++(1,.2);
\draw[thick] (-.5,0) -- ++(1,0);
\draw
[
    decoration={
        coil,
        aspect=0.3, 
        segment length=1.2mm, 
        amplitude=2mm, 
        pre length=3mm,
        post length=3mm},
    decorate
] (0,0) -- ++(0,-2) 
    node[midway,right=0.25cm,black]{$k_{\rm ver}$}; 
\draw[<->] (-.7,0) -- ++ (0,-1.9) node[midway, left] {$L$};

\fill [pattern = north west lines] (-.5,-4.5) rectangle ++(1,-0.2);
\draw[thick] (-.5,-4.5) -- ++(1,0);
\draw
[
    decoration={
        coil,
        aspect=0.3, 
        segment length=1.2mm, 
        amplitude=2mm, 
        pre length=3mm,
        post length=3mm},
    decorate
] (0,-4.5) -- ++(0,+2) 
    node[midway,right=0.25cm,black]{$k_{\rm ver}$}; 

\draw[<->] (-.7,-2.6) -- ++ (0,-1.9) node[midway, left] {$L$};

\fill [pattern = north east lines] (-2.5,-1.75) rectangle ++(-0.2,-1);
\draw[thick] (-2.5,-1.75) -- ++(0,-1);

\draw
[
    decoration={
        coil,
        aspect=0.3, 
        segment length=1.2mm, 
        amplitude=2mm, 
        pre length=3mm,
        post length=3mm},
    decorate
] (-2.5,-2.3) -- ++(2.5,0) 
    node[midway,below=0.25cm,black]{$k_{\rm hor}$}; 

\fill [pattern = north east lines] (+2.5,-1.75) rectangle ++(+0.2,-1);
\draw[thick] (+2.5,-1.75) -- ++(0,-1);

\draw
[
    decoration={
        coil,
        aspect=0.3, 
        segment length=1.2mm, 
        amplitude=2mm, 
        pre length=3mm,
        post length=3mm},
    decorate
] (+2.5,-2.3) -- ++(-2.5,0) 
    node[midway,below=0.25cm,black]{$k_{\rm hor}$}; 

\node[draw,
    fill=gray!60,
    minimum width=.5cm,
    minimum height=.5cm,
    anchor=north,
    label=center:$m$] at (0,-2) {};

\draw[fill=black] (-0.2,-2.5) circle (0.1);
\draw[fill=black] (+0.2,-2.5) circle (0.1);
\draw[fill=black] (-0.2,-2) circle (0.1);
\draw[fill=black] (+0.2,-2) circle (0.1);

\draw[thick] (-2,-2.6) -- ++ (4,0);
\draw[thick] (-2,-1.9) -- ++ (4,0);

\end{tikzpicture}    
\label{fig:undeformed_1dof}
}
\hspace{2cm}
\subfloat[Deformed configuration]{%
\begin{tikzpicture}
\def\xdisp{1.2}
\fill [pattern = north west lines] (-.5,0) rectangle ++(1,.2);
\draw[thick] (-.5,0) -- ++(1,0);
\draw
[
    decoration={
        coil,
        aspect=0.3, 
        segment length=1.2mm, 
        amplitude=2mm, 
        pre length=3mm,
        post length=3mm},
    decorate
] (0,0) -- ++(0+\xdisp,-2) 
    node[midway,right=0.25cm,black]{}; 

\fill [pattern = north west lines] (-.5,-4.5) rectangle ++(1,-0.2);
\draw[thick] (-.5,-4.5) -- ++(1,0);
\draw
[
    decoration={
        coil,
        aspect=0.3, 
        segment length=1.2mm, 
        amplitude=2mm, 
        pre length=3mm,
        post length=3mm},
    decorate
] (0,-4.5) -- ++(0+\xdisp,+2) 
    node[midway,right=0.25cm,black]{};

\fill [pattern = north east lines] (-2.5,-1.75) rectangle ++(-0.2,-1);
\draw[thick] (-2.5,-1.75) -- ++(0,-1);

\draw
[
    decoration={
        coil,
        aspect=0.3, 
        segment length=1.2mm, 
        amplitude=2mm, 
        pre length=3mm,
        post length=3mm},
    decorate
] (-2.5,-2.3) -- ++(2.5+\xdisp,0) 
    node[midway,below=0.25cm,black]{}; 

\fill [pattern = north east lines] (+2.5,-1.75) rectangle ++(+0.2,-1);
\draw[thick] (+2.5,-1.75) -- ++(0,-1);

\draw
[
    decoration={
        coil,
        aspect=0.3, 
        segment length=1.2mm, 
        amplitude=2mm, 
        pre length=3mm,
        post length=3mm},
    decorate
] (+2.5,-2.3) -- ++(-2.5+\xdisp,0) 
    node[midway,below=0.25cm,black]{}; 

\node[draw,
    fill=gray!60,
    minimum width=.5cm,
    minimum height=.5cm,
    anchor=north] at (0+\xdisp,-2) {};

\draw[|->,thick] (0,-1.75) -- ++ (\xdisp,0) node[midway, left=10pt, above=-2pt] {$q$};

\draw[fill=black] (-0.2+\xdisp,-2.5) circle (0.1);
\draw[fill=black] (+0.2+\xdisp,-2.5) circle (0.1);
\draw[fill=black] (-0.2+\xdisp,-2) circle (0.1);
\draw[fill=black] (+0.2+\xdisp,-2) circle (0.1);

\draw[thick] (-2,-2.6) -- ++ (4,0);
\draw[thick] (-2,-1.9) -- ++ (4,0);

\coordinate (O) at (0,-4.5);
\coordinate (A) at (0,-2.5);
\coordinate (B) at (\xdisp,-2.5);
\pic [draw, ->, "$\theta$", angle radius=1cm, angle eccentricity=1.3] {angle=B--O--A};
\end{tikzpicture}
\label{fig:deformed_1dof}
}
\caption{Geometrically nonlinear oscillations of a sliding mass}
\label{fig:1dof_oscillator}
\end{figure*}

Consider a mass sliding frictionless in a horizontal plane. The mass is attached to two horizontal and two vertical springs whose material behavior is linear and given by the stiffness $k_{\rm hor}$ and $k_{\rm ver}$ respectively. In the undeformed configuration the two springs have length $L$ (cf. Fig. \ref{fig:1dof_oscillator}). Denoting by $q$ the horizontal displacement of the mass, the equations of motion are

\begin{equation}\label{eq:nonlinear_oscillator}
m \ddot{q} = - 2 k_{\rm hor} q - 2 k_{\rm ver} \delta \sin(\theta),
\end{equation}
In this expression $\displaystyle \delta = \sqrt{L^2 + q^2} - L$ is the elongation of the vertical springs and $\displaystyle \sin{\theta} = \frac{q}{\sqrt{L^2 + q^2}}$. 
Given the Hooke law $\sigma_{\rm hor} = k_{\rm hor} q, \; \sigma_{\rm ver} = k_{\rm ver} \delta$ and introducing the velocity $v$, the total energy of the system is given by the quadratic form
\begin{equation*}
\begin{aligned}
    H &= \frac{1}{2} m \dot{q}^2 + \frac{1}{2} 2k_{\rm hor}^{-1}\sigma_{\rm hor}^2 + \frac{1}{2} 2k_{\rm ver}^{-1}\sigma_{\rm ver}^2 \\
      &= 
    \frac{1}{2} 
    \begin{pmatrix}
        v \\ 
        \sigma_{\rm hor} \\
        \sigma_{\rm ver} \\
    \end{pmatrix}^\top
    \begin{bmatrix}
            m & 0 & 0 \\
            0 & 2k^{-1}_{\rm hor} & 0 \\
            0 & 0 & 2k^{-1}_{\rm ver}
        \end{bmatrix}\begin{pmatrix}
        v \\
        \sigma_{\rm hor} \\
        \sigma_{\rm ver} \\
    \end{pmatrix}.    
\end{aligned}
\end{equation*}
The time derivative of the elongation of the vertical springs is given by
\begin{equation*}
   \dot{\delta} = \frac{q}{\sqrt{L^2 + q^2}} \dot{q} = \sin(\theta) \dot{q}.
\end{equation*}
The dynamics \eqref{eq:nonlinear_oscillator} can therefore be rewritten as a first order system using the horizontal displacement of the mass $q$, its velocity $v$ and the axial stresses $\sigma_{\rm hor}, \; \sigma_{\rm ver}$ in the following way
\begin{equation*}\label{eq:Ham_nonlinear_oscillator}
\begin{aligned}
    \dot{q} &= \mathbf{e}_1^\top \mathbf{x}, \\
    \mathbf{H}\dot{\mathbf{x}} &= \mathbf{J}(q) \mathbf{x},
\end{aligned}
\end{equation*}
where $\mathbf{e}_1 = [1\quad 0\quad 0]^\top$ is the first element of the Euclidean canonical basis, $\mathbf{x} = [v \; \sigma_{\rm hor} \; \sigma_{\rm ver}]^\top$ and
\begin{equation}
\mathbf{H} = \mathrm{Diag}
        \begin{bmatrix}
            m \\
            2k_{\rm hor}^{-1}\\
            2k_{\rm ver}^{-1}
        \end{bmatrix}, \quad 
\mathbf{J} = \begin{bmatrix}
            0 & -2 & - 2 \sin\theta\\ 
            2 & 0 & 0 \\
            2 \sin\theta & 0 & 0
        \end{bmatrix}.
\end{equation}
The matrix $\mathbf{H}= \mathbf{H}^\top > 0$ is symmetric positive definite, whereas the matrix $\mathbf{J} = -\mathbf{J}^\top$ is skew-symmetric. The energy is given by the quadratic form
$$
H = \frac{1}{2}\mathbf{x}^\top \mathbf{H} \mathbf{x}.
$$
At the price of adding additional variables, the dynamics can be rewritten in this special form. This has important consequences at the numerical level, as it will be shown in section~\ref{sec:time_int}. The system can be approximated by considering a Taylor expansion
\begin{equation*}
    (L^2 + q^2)^{-1/2} = \frac{1}{L} - \frac{q^2}{2L^3} + \frac{3 q^4}{8L^5} + \mathcal{O}(q^6).
\end{equation*}
If only the quadratic term is retained, then the Duffing oscillator \cite{strogatz2018nonlinear} is obtained
\begin{equation}\label{eq:Duffing_oscillator}
m\ddot{q} = - 2k_{\rm hor} q - k_{\rm ver}\frac{q^3}{L^2}.
\end{equation}
The Hamiltonian form for the Duffing oscillator arises from the linearization of the matrix $\mathbf{J}$ in system \eqref{eq:Ham_nonlinear_oscillator} 
\begin{equation}\label{eq:Ham_Duffing_oscillator}
    \begin{aligned}
        \dot{q} &= v, \\
        \mathrm{Diag}
        \begin{bmatrix}
            m \\
            2k_{\rm hor}^{-1}\\
            2k_{\rm ver}^{-1}
        \end{bmatrix}
        \begin{pmatrix}
            \dot{v} \\
            \dot{\sigma}_{\rm hor} \\
            \dot{\sigma}_{\rm ver}
        \end{pmatrix} &= 
        \begin{bmatrix}
            0 & -2 & - \frac{2q}{L} \\
            2 & 0 & 0 \\
            \frac{2q}{L} & 0 & 0
        \end{bmatrix}
        \begin{pmatrix}
            v \\
            \sigma_{\rm hor} \\
            \sigma_{\rm ver}
        \end{pmatrix}.
    \end{aligned}
\end{equation}

\subsection{A general abstract framework}
The previous discussion is instrumental to illustrate that geometrically nonlinear models in mechanics may be written using a total Lagrangian formulation in the following abstract non-canonical Hamiltonian form
\begin{equation}\label{eq:abstract_structure_full}
\begin{aligned}
\partial_t {\bm{q}} &= \bm{v}, \\
    \begin{bmatrix}
        \rho & 0 \\
        0 & \bm{C} 
    \end{bmatrix}
    \pdv{}{t}
    \begin{pmatrix}
    \bm{v} \\
    \bm{S}
    \end{pmatrix} &= 
    \begin{bmatrix}
        0 & -\mathcal{L}^*(\mathcal{D} \bm{q}) \\
        \mathcal{L}(\mathcal{D} \bm{q}) & 0
    \end{bmatrix}
    \begin{pmatrix}
        {\bm{v}} \\
        \bm{S}
    \end{pmatrix}.
\end{aligned}
\end{equation}
Parameters $\rho$ and $\bm{C}$ are related to the density and compliance of the material (in general the latter is a fourth order symmetric tensor), $\bm{q}$ represents a generalized coordinate (or displacement), $\bm{v}$ is the velocity field, $\bm{S}$ a stress-like variable. Furthermore $\mathcal{L}$ is an (unbounded) differential operator, that depends on $\mathcal{D} \bm{q}$, where $\mathcal{D}$ is a differential operator associated to a deformation measure. When $\mathcal{D} \bm{q}$ is regarded as a parameter, the notation $\mathcal{L}^*(\mathcal{D} \bm{q})$ stands for the (formal) adjoint of $\mathcal{L}(\mathcal{D}\bm{q})$, characterized by the relation
\begin{equation}\label{eq:D_adjoint}
    \inpr[\Omega]{\bm{S}}{\mathcal{L}(\mathcal{D}\bm{q})\bm{v}} = \inpr[\Omega]{\mathcal{L}^*(\mathcal{D}\bm{q}) \bm{S}}{\bm{v}}, \qquad \forall \bm{S},\; \forall \bm{v},
\end{equation}
where $\inpr[\Omega]{\bm{f}}{\bm{g}} = \int_{\Omega} \bm{f} \cdot \bm{g} \; \d\Omega$ denotes the inner product of two (generally vector-valued) functions over the domain $\Omega \subset \bbR^d$. For this equation to be true the fields need to satisfy appropriate homogeneous boundary conditions that will be prescribed precisely in the following examples. Because of \eqref{eq:D_adjoint}, the operator
\begin{equation*}
\mathcal{J}(\mathcal{D}\bm{q}) := 
\begin{bmatrix}
0 & -\mathcal{L}^*(\mathcal{D}\bm{q}) \\
\mathcal{L}(\mathcal{D}\bm{q}) & 0
\end{bmatrix},
\end{equation*}
is formally skew-adjoint, meaning that
\begin{equation*}
    \inpr[\Omega]{\bm{\alpha}}{\mathcal{J}(\mathcal{D}\bm{q})\bm{\beta}} = -\inpr[\Omega]{\mathcal{J}^*(\mathcal{D}\bm{q}) \bm{\alpha}}{\bm{\beta}}, \qquad \forall \bm{\alpha}, \forall \bm{\beta}.
\end{equation*}
 Because of this property the energy 
\begin{equation}
    H = \frac{1}{2} \inpr[\Omega]{\bm{v}}{\rho \bm{v}} + \frac{1}{2} \inpr[\Omega]{\bm{S}}{\bm{C} \bm{S}},
\end{equation}
  is conserved, i.e. $\dot{H}= 0$. System $\eqref{eq:abstract_structure_full}$ can be written compactly as follows
  \begin{equation}\label{eq:abstract_structure}
\begin{aligned}
\partial_t {\bm{q}} &= \bm{v}, \\
\mathcal{H}\partial_{t}
    \bm{x} &= \mathcal{J}(\bm{q})\bm{x},
\end{aligned}
\end{equation}
with $\bm{x}= (\bm{v} \; \bm{S})^\top$ and $\mathcal{H} = \mathrm{Diag}[\rho \quad \bm{C}]$. This system can be deduced from the Euler-Lagrange equations arising from the following kinetic and potential energies
\begin{equation*}
    T = \frac{1}{2} \int_\Omega \rho \norm{\pdv{\bm{q}}{t}}^2 \, \d{\Omega}, \quad V = \frac{1}{2}\int_\Omega \bm{E} : \bm{K} \bm{E} \, \d\Omega.
\end{equation*}
The parameter $\bm{K} = \bm{C}^{-1}$ is related to the stiffness tensor and relates stress and strain tensors $\bm{S} = \bm{K}\bm{E}$. The strain $\bm{E} = \bm{E}(\mathcal{D}\bm{q})$ is a nonlinear function of $\mathcal{D}\bm{q}$. Formulation \eqref{eq:abstract_structure} is completely local and requires the specification of the boundary conditions for the overall motion. These will be specified in the following examples.

\subsection{Von K\'arm\'an beams}

\begin{figure}[tbh]
    \centering
    \includegraphics[width=.4\textwidth]{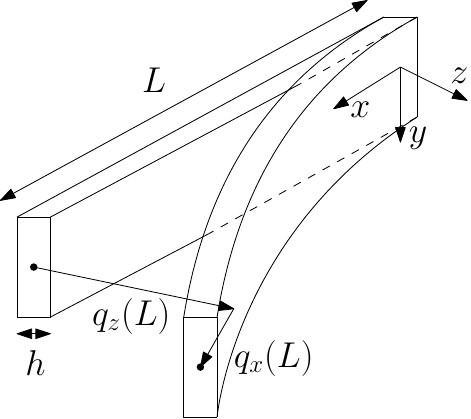}
    \caption{Notation for the kinematic fields of the beam.}
    \label{fig:vK_beam_sketch}
\end{figure}

The von K\'arm\'an model for thin beams is built upon two main geometrical assumptions:
\begin{itemize}
    \item the out of plane deflection is comparable to the thickness;
    \item the squares of axial stretching terms are negligible compared to the square of rotations.
\end{itemize} 
This two assumptions imply that the quadratic term responsible for the bending membrane coupling are retained in the expression of the axial deformation
\begin{equation*}
\varepsilon = \pdv{q_x}{x} +  \frac{1}{2} \left(\pdv{q_z}{x} \right)^2 - z \pdv[order=2]{q_z}{x},
\end{equation*}
where $q_x, \; q_z$ are the horizontal and vertical displacements cf. Fig. \ref{fig:vK_beam_sketch}.
The axial strain and the (linearized) curvature are given by
$$
\varepsilon_a: =\pdv{q_x}{x} +  \frac{1}{2} \left(\pdv{q_z}{x} \right)^2, \qquad \kappa := \pdv[order=2]{q_z}{x}.
$$
Consider the kinetic and potential energy
\begin{equation*}
\begin{aligned}
    T &= \frac{1}{2} \int_0^L \rho A \left\{ \left(\pdv{q_x}{t}\right)^2 + \left(\pdv{q_z}{t}\right)^2 \right\} \d{x}, \\
    V &= \frac{1}{2} \int_0^L \{ EA \varepsilon_a^2 + EI \kappa^2 \} \d{x}.
\end{aligned}
\end{equation*}
The Euler-Lagrange equations obtained by the principle of least-action are given by 
\begin{equation*}
\begin{aligned}
\rho A\, \partial_{tt}{q_x} &= \partial_x  N, \\
\rho A\, \partial_{tt}{q_z} &= -\partial^2_{xx} M + \partial_x(N\,\partial_x q_z),
\end{aligned} 
\end{equation*}
where the axial and bending stress resultant have been introduced
\begin{equation*}
    N: = EA \varepsilon_a, \qquad M:= EI \kappa.
\end{equation*}
The time derivative of the stress variables gives
\begin{align*}
    C_a \partial_t N &= \partial_x v_x + (\partial_x q_z) \partial_x v_z, \\
    C_b \partial_t M &= \partial_{xx} v_z,
\end{align*}
where the bending compliance $C_b := (EI)^{-1}$  has been introduced. By including the dynamics of these variables the 
Hamiltonian formulation \eqref{eq:abstract_structure} is obtained \cite{brugnoli2021vonkarman}, where the state and operators take the following specific form 
\begin{equation*}
\begin{aligned}
\bm{x} &:= \begin{pmatrix}
{v}_x & {v}_z  & {N} & {M}
\end{pmatrix}^\top, \\
\mathcal{H} &:= \mathrm{Diag}
\begin{bmatrix}
\rho A & \rho A & C_a & C_b
\end{bmatrix}^\top, \\
\mathcal{J}(q_z) &:= 
\begin{bmatrix}
    0 & 0 & \partial_x & 0\\
    0 & 0 & {\partial_x(\circ \, \partial_x q_z)} & -\partial_{xx}^2 \\
    \partial_x &  (\partial_x q_z)  \partial_x  \circ  & 0 & 0 \\
    0 & \partial_{xx}^2 & 0 & 0 \\ 
\end{bmatrix}.
\end{aligned}
\end{equation*}
Given a splitting of the boundary $\partial\Omega=\Gamma_q \cup \Gamma_\sigma$, the boundary conditions are given by 
\begin{equation*}
    \begin{aligned}
        q_x|_{\Gamma_q} = 0, \\
        q_z|_{\Gamma_q} = 0, \\
        \partial_x q_z|_{\Gamma_q} = 0,
    \end{aligned} \qquad
    \begin{aligned}
        N|_{\Gamma_\sigma} =0, \\
        (N\partial_x q_z - \partial_x M)|_{\Gamma_\sigma} =0, \\
        M|_{\Gamma_\sigma} =0, \\
    \end{aligned}
\end{equation*}
\begin{remark}
    The aforementioned boundary conditions describe a clamped-free beam. Simply supported boundary conditions may also be considered. Imagine then a splitting of the form $\partial\Omega=\Gamma_q \cup \Gamma_\sigma \cup \Gamma_{\rm ss}$, then the simply supported boundary conditions on $\Gamma_{\rm ss}$ are 
    $$
    q_x|_{\Gamma_{\rm ss}} = 0, \qquad 
    q_z|_{\Gamma_{\rm ss}} = 0, \qquad M|_{\Gamma_{\rm ss}} =0.
    $$
\end{remark}

For this example, the operator $\mathcal{L}(\partial_x q_z)$ and its adjoint are given by
\begin{equation*}
\begin{aligned}
    \mathcal{L}(\partial_x q_z)& = 
    \begin{bmatrix}
    \partial_x &  (\partial_x q_z) \partial_x \circ \\
    0 & \partial_{xx}^2
    \end{bmatrix}, \\ 
    \mathcal{L}^*(\partial_x q_z) &= -
    \begin{bmatrix}
    \partial_x &  0\\
    \partial_x( \circ \partial_x q_z) & -\partial_{xx}^2
    \end{bmatrix}.
\end{aligned}
\end{equation*}

\subsection{Von K\'arm\'an plate model}\label{sec:vK}

This model is the two-dimensional extension of the previous one, and it is built on the same geometrical assumptions.    This means that the strain tensor $\bm{\varepsilon} \in \bbR^{2\times 2}_{\rm sym}$, takes the following simplified expression 
\begin{equation}
    \bm{\varepsilon} = \mathrm{def}(\bm{q}_m)  + \frac{1}{2} \nabla q_z \otimes \nabla q_z - z\Hess q_z.
\end{equation}
In this equation $\mathrm{def} := \frac{1}{2}(\nabla + \nabla^\top)$, which stands for {\em deformation},  is the symmetric gradient and corresponds to the infinitesimal strain tensor in linear elasticity, and $\otimes$ is the dyadic product of the two vectors, i.e. $\mathbf{a} \otimes \mathbf{b} = \mathbf{a}\mathbf{b}^\top$ where $\mathbf{a}, \mathbf{b} \in \bbR^2$ are column vectors and $\Hess$ denotes the Hessian operator. The displacement vector $\bm{q}$ has been split into the membrane displacement $\bm{q}_m = (q_x \quad q_y)^\top$ and the out-of-plane component $q_z$. Again the strain tensor can be split into the membrane strain and the (linearized) curvature
\begin{equation*}
    \bm{\varepsilon}_m = \mathrm{def}(\bm{q}_m)  + \frac{1}{2} \nabla q_z \otimes \nabla q_z, \quad \bm{\kappa} = \Hess q_z.
\end{equation*}
\subsubsection{Full model}
The following model has been presented in a concise manner in \cite{brugnoli2022enoc}. The wording {\em full} refers to the fact that both the membrane and bending behavior are considered in this model. Given a two-dimensional domain $\Omega \subset \bbR^2$, the kinetic and potential energies are 
$$
\begin{aligned}
T &= \frac{1}{2} \int_{\Omega} \rho h   \left\{\norm{\pdv{\bm{q}_m}{t}}^2 + \left(\pdv{q_z}{t}\right)^2\right\} \d\Omega, \\
V &= \frac{1}{2} \int_{\Omega} \{ \bm{\varepsilon}_m : \bm{K}_m \bm{\varepsilon}_m + \bm{\kappa} : \bm{K}_b \bm{\kappa}\}\,  \d\Omega,
\end{aligned}
$$
where $h$ is the plate thickness and $\bm{K}_m, \; \bm{K}_b$ are fourth order tensors representing the membrane and bending stiffness. The contraction of two matrices is denoted as $\bm{A} : \bm{B} = \sum_{i, j=1}^n A_{ij} B_{ij}$. In the case of an isotropic homogeneous material in-plane stress condition, their action on a symmetric second-order tensor $\bm{S} \in \bbR^{2\times 2}_{\rm sym}$  takes the following form
\begin{equation*}
\bm{K}_m(\bm{S}) = Eh\bm{\Phi}(\bm{S}), \qquad \bm{K}_b = \frac{Eh^3}{12}\bm{\Phi}(\bm{S}).
\end{equation*}
where $\bm{\Phi}: \bbR^{2\times 2}_{\rm sym} \rightarrow \bbR^{2\times 2}_{\rm sym}$ is the following linear map between symmetric tensors
$$
\bm{\Phi}(\bm{S}) = \frac{1}{1 - \nu^2}\{(1-\nu)\bm{S} + \nu \bm{I}_2 \tr(\bm{S}) \}.
$$
The membrane and bending stresses are given by $\bm{N}= \bm{K}_m \bm\varepsilon_m, \; \bm{M} = \bm{K}_b \bm{\kappa}$ respectively. The Euler-Lagrange equations for this model are \rewother{(cf. \cite{bilbao2015conservative})}
	\begin{equation*}
		\begin{aligned}
			\rho h\, \partial_{tt}{\bm{q}_m} &= \div \bm{N}, \\
			\rho h\, \partial_{tt}{q_z} &= -\div\div\bm{M} + \div(\bm{N}\,\nabla q_z),
		\end{aligned} 
	\end{equation*}
\rewtwo{The operator $\div$ applied to a tensor is its row-wise divergence. Its coordinates expression reads 
$$[\div \bm{A}]_i = \sum_{j=1}^3 \partial_{x_j} \bm{A}_{ij},$$
where $\bm{A}: \Omega \rightarrow \mathbb{R}^{3\times 3}$ is a tensor field.} The time derivative of the strain tensors gives
\begin{equation*}
    \begin{aligned}
        \bm{C}_m\, \partial_t \bm{N} &= \mathrm{def} \bm{v}_m + \mathrm{sym}(\nabla q_z \otimes \nabla v_z), \\
        \bm{C}_b\, \partial_t \bm{M} &= \Hess v_z,
    \end{aligned}
\end{equation*}
where $\bm{C}_m := \bm{K}_m^{-1}, \; \bm{C}_b := \bm{K}_b^{-1}$ and $\mathrm{sym}$ is the symmetrization operator. Consequently the dynamics assumes the same form as in \eqref{eq:abstract_structure}, with state and operators taking the following specific form 
\begin{equation}\label{eq:vK_full}
\begin{aligned}
\bm{x} &:= \begin{pmatrix}
{v}_x & {v}_z  & \bm{N} & \bm{M}
\end{pmatrix}^\top, \\
\mathcal{H} &:= \mathrm{Diag}
\begin{bmatrix}
\rho h & \rho h & \bm{C}_m & \bm{C}_b
\end{bmatrix}^\top, \\
\mathcal{J}(q_z) &:= 
{\scriptstyle
\begin{bmatrix} 
0 & 0 & \div & 0 \\
0 & 0 &  \mathcal{C}(q_z) & -\div\div \\
\mathrm{def} & -\mathcal{C}^*(q_z) & 0 & {0} \\
0 & \Hess & 0 & 0 \\
\end{bmatrix},
}
\end{aligned}
\end{equation}
where the operator $\mathcal{C}(q_z)(\cdot): L^2(\Omega, \bbR^{2\times 2}_{\text{sym}}) \rightarrow L^2(\Omega)$ acting on symmetric tensors is defined by
\begin{equation*}
	\mathcal{C}(q_z)(\bm{N}) = \div(\bm{N}\,\nabla q_z).
\end{equation*}
We are going to show in Proposition \ref{prop:adj_C} that its adjoint is given by 
$$\mathcal{C}(q_z)^*(\cdot) = - \mathrm{sym}\left[\nabla (\cdot) \otimes \nabla(q_z)\right].$$
The kinematic boundary conditions are applied on $\Gamma_q$ whereas the dynamic boundary conditions are applied on $\Gamma_\sigma$ ($\partial\Omega=\Gamma_q \cup \Gamma_\sigma$)
\begin{equation*}
    \begin{aligned}
    \bm{q}_m|_{\Gamma_q} = 0, \\
    q_z|_{\Gamma_q} = 0, \\
    \partial_{\bm{n}}q_z |_{\Gamma_q} = 0,
    \end{aligned} \qquad 
    \begin{aligned}
    \bm{N}\bm{n}|_{\Gamma_\sigma} = 0, \\
    \gamma_{\perp\perp, 1} \bm{M} +  (\bm{N}\bm{n})^\top\nabla q_z|_{\Gamma_\sigma} = 0, \\
    \bm{n}^\top \bm{M} \bm{n}|_{\Gamma_\sigma} = 0,
    \end{aligned}
\end{equation*}
where $\bm{n}$ is the outward normal vector and $\gamma_{\perp\perp, 1} \bm{M} = - \bm{n}^\top \div \bm{M} - \partial_{\bm{s}}(\bm{n}^\top \bm{M} \bm{s})\vert_{\partial\Omega}$ is the effective shear force at the boundary ($\bm{s}$ is the tangent vector at the boundary).
\begin{proposition}\label{prop:adj_C}
In this case, the differential operator and its adjoint read
$$
\begin{aligned}
\mathcal{L}(\nabla q_z) = \begin{bmatrix}
    \mathrm{def} & \mathrm{sym}\left(\nabla q_z \otimes \nabla \circ\right) \\
    0 & \Hess
\end{bmatrix}, \\
\mathcal{L}^*(\nabla q_z) = -\begin{bmatrix}
    \div & 0 \\
    \div(\circ \nabla q_z) & -\div\div
\end{bmatrix}.    
\end{aligned}
$$    
\end{proposition}

\begin{proof}
    The fact that $\mathrm{def}$ is the adjoint of $-\div$ acting on symmetric tensors is known, see e.g.  \cite{brugnoli2019mindlin}. Moreover, the proof that $\div\div$ and $\Hess$ are adjoint operators is given in \cite{brugnoli2019kirchhoff}.\\
Consider a smooth scalar field  $u \in C^\infty_0(\Omega)$ and a smooth symmetric tensor field $\bm{N} \in C^\infty_0(\Omega, \bbR^{2\times 2}_{\text{sym}})$ with compact support. The formal adjoint of $\mathcal{C}(q_z)(\cdot)$ satisfies the relation
\begin{equation*}
\inpr[L^2(\Omega)]{u}{\mathcal{C}(q_z)(\bm{N})} = \inpr[L^2(\Omega, \bbR^{2\times 2}_{\text{sym}})]{\mathcal{C}(q_z)^*(u)}{\bm{N}}.
\end{equation*}
The adjoint is deduced by the following computation
\begin{equation*}
\begin{aligned}
\inpr[\Omega]{u}{\mathcal{C}(q_z)(\bm{N})} &= \inpr[\Omega]{u}{\div(\bm{N} \nabla q_z)},\\
&= \inpr[\Omega]{-\nabla u}{\bm{N} \nabla q_z},\\
&= \inpr[\Omega]{-\nabla u \otimes \nabla q_z}{\bm{N}},\\
&= \inpr[\Omega]{-\mathrm{sym}(\nabla u \otimes \nabla q_z)}{\bm{N}}.
\end{aligned} 
\end{equation*}
The first equality follows from integration by parts, the second by the dyadic product properties and the third by the symmetry of $\bm{N}$.
This means
\begin{equation*}
\mathcal{C}(q_z)^*(\cdot) = - \mathrm{sym} \left[\nabla (\cdot) \otimes \nabla(q_z)\right],
\end{equation*}
leading to the final result.
\end{proof}

\subsubsection{Von K\'arm\'an plate in Airy form}

Consider now the case of negligible membrane inertia. This means that the following constraint is imposed on the membrane stress tensor
$$
\div \bm{N} = 0.
$$
To simplify the problems the mathematical structure of the elasticity complex is exploited. For the reader convenience, the notation for the differential operators and the main results for the two-dimensional elasticity complex are summarized in Appendix \ref{app:diff_2D} and \ref{app:complex} respectively.  For a simply connected two-dimensional domain the elasticity complex \eqref{eq:air_complex} is exact and the Airy stress potential can be used to deduce the membrane stress as
$$
\bm{N} = \Air{\varphi}.
$$
This expression can be used to simplify the dynamics. In particular the bending-membrane coupling term can be simplified as follows
\begin{equation*}
\begin{aligned}
\div(\bm{N}\,\nabla q_z) = \div(\Air\varphi\,\nabla q_z), \\
 = \div(\Air\varphi) \cdot \nabla q_z + \Air\varphi : \Hess q_z, \\
 =\Air\varphi : \Hess q_z,
\end{aligned}
\end{equation*}
since $\div \Air \equiv 0$. As is customary in the literature, the last term is denoted using the following bilinear operator
\begin{equation*}
\begin{aligned}
    \mathcal{B}(f, g) := \Air f : \Hess g \\
    =(\partial_{22}f) (\partial_{11} g) + (\partial_{11}f) (\partial_{22}g) - 2 (\partial_{12}f) (\partial_{12}g).
\end{aligned}
\end{equation*}
The following properties of this bilinear form have been proven in \cite{bilbao2008conservative}
\begin{itemize}
    \item Symmetry: $\mathcal{B}(f, g) = \mathcal{B}(g, f)$;
    \item Self-adjointness: $\inpr[\Omega]{\mathcal{B}(f, g)}{h} = \inpr[\Omega]{g}{\mathcal{B}(f, h)}$ (function $f$ is here regarded as a parameter).      
\end{itemize}
The expression of the membrane strain still contains the contribution of the in-plane displacement
$$
\bm{C}_m \Air \varphi = \sgrad\bm{q}_m + \frac{1}{2} \nabla q_z \otimes \nabla q_z.
$$
The idea is to exploit the relation 
$$\rot\rot \sgrad=0$$
from complex \eqref{eq:rotrot_complex} (see Appendix \ref{app:diff_2D}) to eliminate the membrane bending coupling.
$$
\rot\rot \bm{C}_m \Air \varphi = \frac{1}{2} \rot\rot(\nabla q_z \otimes \nabla q_z).
$$
Notice that $\Air^* = \rot\rot$, so the operator $\rot\rot \bm{C}_m \Air =  \mathrm{Air}^* \bm{C}_m \Air$ is self-adjoint. A little algebra provides 
$$
\rot\rot(\nabla q_z \otimes \nabla q_z)= - \mathcal{B}(q_z, q_z).  
$$
So the Airy potential is related to the out-of-plane displacement via
\begin{equation}\label{eq:airy_displacement}
(\mathrm{Air}^*\bm{C}_m \Air)\,\varphi = - \frac{1}{2} \mathcal{B}(q_z, q_z).
\end{equation}
If the material is isotropic, the above relation simplifies into
$$
\Delta^2 q_z =  - \frac{Eh}{2} \mathcal{B}(q_z, q_z),
$$
where $\Delta^2$ is the bi-Laplacian. For the sake of generality, a generic compliance tensor is considered. The derivative in time of Equation \eqref{eq:airy_displacement} provides the dynamics of the Airy potential $\varphi$
$$
(\mathrm{Air}^* \bm{C}_m \Air)\,\pdv{\varphi}{t} = - \mathcal{B}(q_z, v_z),
$$
where the symmetry and bilinearity of $\mathcal{B}$ have been exploited. The von-K\'arm\'an plate can then be put in Hamiltonian form considering all the aforementioned simplification as in \eqref{eq:abstract_structure} with
\begin{equation}\label{eq:vK_airy}
\begin{aligned}
\bm{x} &:= \begin{pmatrix}
{v}_z \quad \varphi \quad \bm{M}
\end{pmatrix}, \\
\mathcal{H} &:=
\mathrm{Diag}
\begin{bmatrix}
\rho h \quad \mathrm{Air}^* \bm{C}_m \Air \quad \bm{C}_b
\end{bmatrix},\\
\mathcal{J} &:=
\begin{bmatrix}
    0 &  \mathcal{B}(q_z, \circ) & -\div\div \\
        - \mathcal{B}(q_z, \circ) & 0 & 0 \\
        \Hess & 0 & 0 \\ 
\end{bmatrix}.
\end{aligned}
\end{equation}
Since the operator $\mathcal{B}$ is self adjoint, the dynamics is ruled by a skew-adjoint operators as in the previous case. A self-adjoint differential operator takes the place of an algebraic energy matrix. The differential operator and its adjoint read
$$
\begin{aligned}
\mathcal{L}(\Air q_z) &= \begin{bmatrix}
-\mathcal{B}(q_z, \circ) \\
\Hess
\end{bmatrix}, \\
\mathcal{L}^*(\Air q_z) &= -\begin{bmatrix}
\mathcal{B}(q_z, \circ) \quad -\div\div
\end{bmatrix}.
\end{aligned}
$$

\subsection{Geometrically nonlinear elasticity}

In geometrically nonlinear elasticity different tensors may be used to describe deformations. In this article we focus on the Green-Lagrange tensor 
$$
\bm{E} := \frac{1}{2}(\bm{F}^\top\bm{F} - \bm{I}), \qquad \bm{F} :=  \bm{I} + \nabla \bm{q}.
$$
where $[\nabla \bm{q}]_{ij} = \partial_j q_i$ is the gradient of a vector defined row-wise. The kinetic and potential energies are given by 
$$
\begin{aligned}
T &= \frac{1}{2} \int_\Omega \rho \norm{\pdv{\bm{q}}{t}}^2 \d\Omega, \\
V &= \frac{1}{2} \int_\Omega \bm{E}:\bm{K}\bm{E} \; \d\Omega.
\end{aligned}
$$
For the potential energy a Saint-Venant \rewtwo{Kirchhoff} material model has been used. It is well known that this material behaviour exhibits numerical instabilities under compressive loading \cite{sifakis2012}. Nevertheless it remains of interest in many applications as it can be computed under severe real time limitations. The Euler-Lagrange equations are then given by
\begin{equation*}
\rho\, \partial_{tt}\bm{q} = \div(\bm{F}\bm{S}),
\end{equation*}
where $\bm{S}=\bm{K}\bm{E}$ is the \rewtwo{second Piola-Kirchhoff} stress tensor.  By introducing the dynamical equation for the \rewtwo{second Piola-Kirchhoff} stress tensor, the Hamiltonian structure of the equations can be highlighted \cite{thomas2024velocity}:
\begin{equation*}
\begin{aligned}
\partial_t \bm{q} &= \bm{v}, \\
\begin{bmatrix}
    \rho & 0 \\
    0 & \bm{C}
\end{bmatrix}
\pdv{}{t}
\begin{pmatrix}
    \bm{v} \\
    \bm{S}
\end{pmatrix} &= 
\begin{bmatrix}
    0 & \div (\bm{F} \; \circ) \\
    \sym( \bm{F}^\top \nabla \; \circ) & 0 \\
\end{bmatrix}
\begin{pmatrix}
    \bm{v} \\
    \bm{S}
\end{pmatrix},    
\end{aligned}
\end{equation*}
where $\bm{C}:=\bm{K}^{-1}$ is the compliance tensor. The boundary conditions are 
\begin{equation*}
    \bm{q}|_{\Gamma_q} = 0, \qquad 
    \bm{F}\bm{S}\bm{n}|_{\Gamma_\sigma} = 0.
\end{equation*}
The differential operator and its adjoint read
$$
\begin{aligned}
    \mathcal{L}(\nabla \bm{q}) &= \sym( \bm{F}^\top \nabla \; \circ), \\
\mathcal{L}^*(\nabla \bm{q}) &= -  \div (\bm{F} \; \circ).
\end{aligned}
$$

\section{Linearly implicit energy-preserving integration}\label{sec:discr}
We  consider here a mixed finite element discretization strategy together with a time integration method to preserve the energy of the system. First we detail the space discretization methodology. The general procedure is first explained using the abstract formulation, and then specialized for the two models considered in the numerical examples, i.e. the von-K\'arm\'an beam model and geometrically nonlinear elasticity. Next the time integration is discussed.

\subsection{Semi-discretization in space with mixed finite element}\label{sec:mfem}

\begin{figure*}[ht]
\centering
\begin{tikzpicture}[scale=0.5]
\draw [->] (0,0) -- (0,5) ;

\draw [->] (0,0) -- (4,0) -- (8,0)  -- (12,0) -- (16,0)  -- (18,0) ;

\draw[<->] (0,-0.5) -- (4,-0.5) node[midway, below] {$h$};

\node[circle, fill=black, inner sep=2pt] at (0,0) {};
\node[circle, fill=black, inner sep=2pt] at (4,0) {};
\node[circle, fill=black, inner sep=2pt] at (8,0) {};
\node[circle, fill=black, inner sep=2pt] at (12,0) {};
\node[circle, fill=black, inner sep=2pt] at (16,0) {};

\draw[thick] (4,0) -- (8,4) node[above]{$\varphi_{3}$} ;
\draw[thick] (8,4) -- (12,0) ;

\draw (0,4) node[left]{$1$} [dotted] -- (16,4) ;

\draw [->] (0,-9) -- (0,-3) ;

\draw [->] (0,-6) -- (4,-6) -- (8,-6)  -- (12,-6) -- (16,-6) -- (18,-6) ;

\node[circle, fill=black, inner sep=2pt] at (0,-6) {};
\node[circle, fill=black, inner sep=2pt] at (4,-6) {};
\node[circle, fill=black, inner sep=2pt] at (8,-6) {};
\node[circle, fill=black, inner sep=2pt] at (12,-6) {};
\node[circle, fill=black, inner sep=2pt] at (16,-6) {};

\draw[fill=gray!20, dashed] (4,-6) -- (4,-4) -- (8,-4) -- (8,-6) -- cycle ;

\draw (0,-4) node[left]{$\frac{1}{h}$} [dotted] -- (16,-4) ;

\draw[fill=gray!20, dashed] (8,-6) -- (12,-6) -- (12,-8) -- (8,-8) -- cycle ;

\draw (0,-8) node[left]{$-\frac{1}{h}$} [dotted] -- (16,-8) ;

\draw[thick] (4,-4) -- (8,-4) node[above]{$\partial_x \varphi_{3}$} ;
\draw[thick] (8,-8) -- (12,-8) ;

\draw[->, bend right] (20, 2) to[out=90, in=90] (20,-6);

\node at (24,-2) {$\pdv{}{x}$};

\end{tikzpicture}
\caption{Derivative of a Lagrange space $\mathrm{CG}_1$, leading to a piecewise constant function.}
\label{fig:derivative_Lagrange}
\end{figure*}
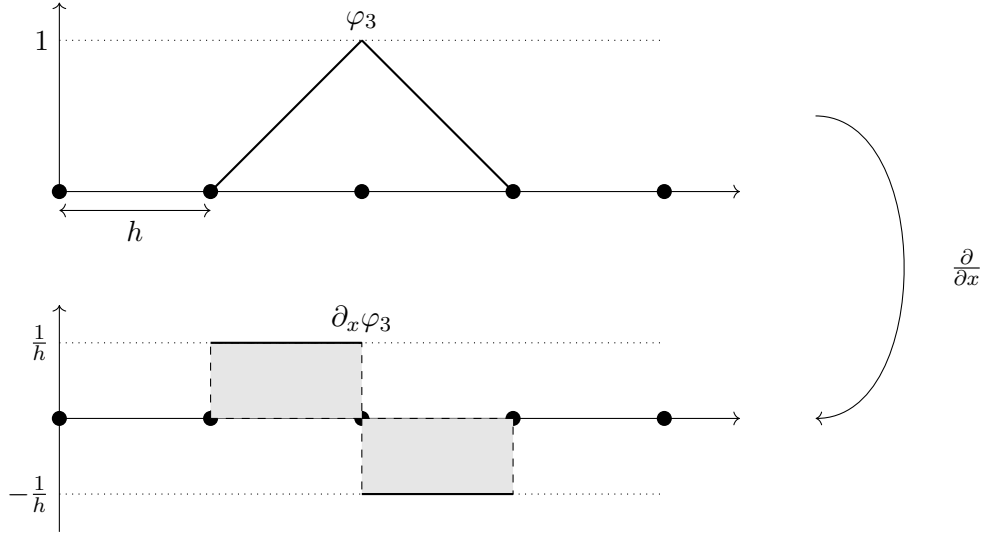

\rewone{
To illustrate the idea behind mixed finite element formulations, consider the wave equation in 1D, describing the longitudinal wave propagation in a bar under Neumann boundary conditions
$$
\begin{aligned}
\rho A\, \partial_{tt} q_x - \partial_x (EA\, \partial_x q_x) &= 0, \quad \text{in } \Omega=[0, L], \\
EA\, \partial_x q_x|_{x=0} = EA\, \partial_x q_x|_{x=L} &= 0,
\end{aligned}
$$
where $\rho$ is the density, $A$ the cross section area, $E$ the Young modulus and $L$ the length of the bar. The unknown $q_x$ represents the longitudinal displacement in the bar. A classical finite element discretization using Lagrange polynomials leads to  
\begin{equation*}
\mathbf{M} \ddot{\mathbf{q}}_x + \mathbf{K} \mathbf{q}_x = 0.
\end{equation*}
The $(i,j)$ component of the mass and stiffness matrix are obtained as follows
$$
\begin{aligned}
    [\mathbf{M}]_{ij} &= \int_0^L \rho A \varphi_i \varphi_j \d{x}, \\
    [\mathbf{K}]_{ij} &= \int_0^L (\partial_x \varphi_i) EA (\partial_x \varphi_j) \d{x},
\end{aligned}
$$
where $\varphi_i(x), \; \varphi_j(x)$ is the Lagrange basis associated with the node $i$ and $j$ respectively. 
As illustrated in the case of the Duffing oscillator, the proposed formulation uses as variables the velocity and the axial stress resultant, defined by
$$
v_x:=\partial_t q_x, \qquad N := EA \partial_x q_x.
$$
As done in the Duffing example we also introduce the dynamics of the stress by moving to the left the axial stiffness and taking the derivative
$$
(EA)^{-1} \partial_t N := \partial_t \partial_x q_x = \partial_x v_x
$$
The latter equality is obtained from the fact that for sufficiently regular functions higher order derivatives commute. The wave equation can then rewritten as a first order system in space and time (system of conservation laws)
$$
\begin{aligned}
\rho A\, \partial_t v_x &= \partial_x \sigma, \\    
C_a\, \partial_t N &= \partial_x v_x,
\end{aligned} 
$$
where $C_a := (EA)^{-1}$ is the axial compliance. To obtain a finite element formulation we consider a weak formulation (multiplication by test function $w, \; \tau$ and integration over the domain) and integrate by parts the first line only
\begin{align}\label{eq:wave_mixed}
    \int_0^L w\, \rho A \partial_t v_x \; \d{x} &= - \int_0^L (\partial_x w)\, N \; \d{x}, \qquad \forall w, \\
    \int_0^L \tau\, C_a \partial_t N \; \d{x} &= + \int_0^L \tau\, (\partial_x v_x) \; \d{x}, \qquad \forall \tau. \label{eq:N_dynamics}
\end{align} 
The boundary terms arising from the integration by parts vanish because of the Neumann conditions. To discretize the system a Galerkin formulation is used. However, velocity and stress are discretized with different bases
\begin{align*}
    v_{x, h} &= \sum_{i=1}^{N_v} \varphi_i(x) v_{x, i}, \qquad  v_{x, h} \in V_h := \mathrm{span}\{\varphi_i\}_i^{N_v}, \\
    N_h &= \sum_{i=1}^{N_\sigma} \psi_i(x) n_{i}, \qquad  N_h \in \Sigma_h := \mathrm{span}\{\psi_j\}_j^{N_\sigma}.
\end{align*}
Test functions belong to the same space of the associated variable, i.e. $w \in V_h, \; \tau \in \Sigma_h$. This is an important feature of mixed finite elements: they require the simultaneous approximation of different variables. Since $v$ undergoes differentiation, it can be approximate by Lagrange finite element, i.e hat functions. For the approximation of axial stress resultant $N$ it is important to select a space that satisfy the following inclusion
$$
\partial_x V_h \subset \Sigma_h.
$$
Indeed if the inclusion hold, then Eq. \eqref{eq:N_dynamics} holds pointwise
$$
C_a \partial_t N_h = \partial_x v_h.
$$
Suppose that Lagrange elements of order 1, denoted by $\mathrm{CG}_1$, are used for $V_h$. Then their derivative give rise to a piecewise constant space, as illustrated in Fig.~\ref{fig:derivative_Lagrange}. 
}

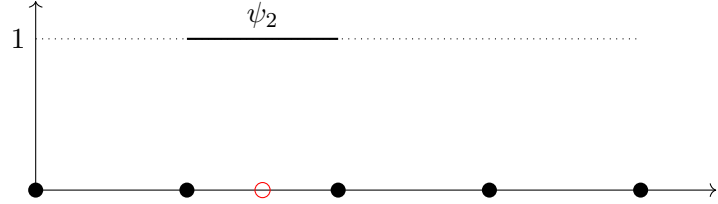
\begin{figure*}[hbt]
\centering
\begin{tikzpicture}[scale=0.5]
\draw [->] (0,0) -- (0,5) ;
\draw [->] (0,0) -- (4,0) -- (8,0)  -- (12,0) -- (16,0)  -- (18,0) ;
\node[circle, fill=black, inner sep=2pt] at (0,0) {};
\node[circle, fill=black, inner sep=2pt] at (4,0) {};
\node[circle, fill=black, inner sep=2pt] at (8,0) {};
\node[circle, fill=black, inner sep=2pt] at (12,0) {};
\node[circle, fill=black, inner sep=2pt] at (16,0) {};

\draw (0,4) node[left]{$1$} [dotted] -- (16,4) ;
\draw[thick] (4,4) -- (8,4) node[midway,above]{$\psi_{2}$};
\draw[red] (6,0) circle (.2cm);
\end{tikzpicture}
    \caption{The second basis function $\psi_2$ of the space $\mathrm{DG}_0$}
    \label{fig:basis_DG0}
\end{figure*}

\rewone{
A discontinuous and constant finite element space, can be used to represent the stress $N$. This space is called discontinuous Galerkin (or Lagrange) space and it is denoted via $\mathrm{DG}_0$ (see Fig. \ref{fig:basis_DG0}).
These continuous Lagrange space of degree 1 and the piecewise constant discontinuous space satisfy $\partial_x \mathrm{CG}_1 \subset \mathrm{DG}_{0}$. This is true  no matter the  polynomial degree $k$ used for the two basis
$$
\partial_x \mathrm{CG}_k \subset \mathrm{DG}_{k-1}.
$$
So the space for the stress is completely local. In the following we assume that 
$$
V_h = \mathrm{CG}_k, \qquad \Sigma_h = \mathrm{DG}_{k-1}.
$$
The algebraic version of system \eqref{eq:wave_mixed} is given by
\begin{equation}\label{eq:discrete_wave_mixed}
    \begin{bmatrix}
        \mathbf{M}_{\rho A} & 0 \\
        0 & \mathbf{M}_{C_a}
    \end{bmatrix}
    \odv{}{t}
    \begin{pmatrix}
        \mathbf{v}_x \\
        \mathbf{n}
    \end{pmatrix} = 
    \begin{bmatrix}
        0 & -\mathbf{D}_{\partial_x}^\top \\
        \mathbf{D}_{\partial_x} & 0 
    \end{bmatrix}
    \begin{pmatrix}
        \mathbf{v}_x \\
        \mathbf{n}
    \end{pmatrix}.
\end{equation}
The matrices are computed as follows
\begin{equation*}
\begin{aligned}
    [\mathbf{M}_{\rho A}]_{ij} &= \int_0^L \rho A \varphi_i \cdot \varphi_j \; \d{x}, \\
    [\mathbf{M}_{C_a}]_{ij} &= \int_0^L C_a \psi_i \cdot \psi_j \;\d{x}, \\
    [\mathbf{D}_{\partial_x}]_{ij} &= \int_0^L \varphi_i \pdv{\psi_j}{x} \; \d{x}.
\end{aligned}
\end{equation*}
All the matrices appearing in this system are sparse. Moreover because of the local nature of the discontinuous Galerkin space, the mass matrix
$\mathbf{M}_{C_a}$, that is associated with the mechanical compliance, is block diagonal. For instance, if piecewise constant function are chosen for the stress, i.e. $N_h \in \mathrm{DG}_0,$ then 
$$\mathbf{M}_{C_a} = (C_a h)\, \mathbf{I},$$ 
where $h$ is the element size and $\mathbf{I}$ the identity in $\mathbb{R}^{N_\sigma \times N_\sigma}$. Furthermore for $k>0$ if Legendre nodes are used to construct the space $\mathrm{DG}_k$ the matrix $\mathbf{M}_{C_a}$ remains diagonal. In the case of constant physical coefficients, it is possible to prove that the classical finite element stiffness $\mathbf{K}$ is related by the matrices appearing in system \eqref{eq:discrete_wave_mixed} by the following decomposition 
$$
\mathbf{K} = \mathbf{D}_{\partial_x}^\top \mathbf{M}_{C_a}^{-1} \mathbf{D}_{\partial_x}.
$$
The inverse $\mathbf{M}_{C_a}^{-1}$ can be easily computed as the matrix is block diagonal. The fact that the stiffness matrix is singular when free boundary conditions are applied comes from the fact that $\mathbf{D}_{\partial_x}$ has kernel given by constant vectors (rigid translation), indeed 
$$
\mathbf{D}_{\partial_x}\bm{1}=0, \qquad \text{where } \bm{1}:=[1 \; \dots \; 1]^\top.
$$
This does not pose any problem for the dynamic resolution. Indeed if the implicit midpoint rule is used, the matrix to be solved is
$$
\begin{aligned}
    \mathbf{A} &= \mathbf{M} - \frac{\Delta t}{2} \mathbf{J},\\
                &=  \begin{bmatrix}
\mathbf{M}_{\rho A} & 0 \\
        0 & \mathbf{M}_{C_a}
\end{bmatrix}  - \frac{\Delta t}{2} \begin{bmatrix}
        0 & -\mathbf{D}_{\partial_x}^\top \\
        \mathbf{D}_{\partial_x} & 0 
    \end{bmatrix}.
\end{aligned}
$$
Matrix $\mathbf{A}$ is invertible. The proof will be provided in Proposition \ref{pr:invertibility} for a more general case.
}

\subsection{Discretization of the abstract formulation}
\rewtwo{As explained in the previous section, the discretization is built upon a weak formulation. The integrations by parts is performed on the the equation describing the linear momentum balance.} Variables are then discretized via a suitable finite element space. The discrete abstract problem in weak form is: find $\bm{q}_h, \bm{v}_h \in V_h, \; \bm{S} \in \Sigma_h$ such that 
\begin{equation}\label{eq:weak_abstract}
\begin{aligned}
    \partial_t \bm{q}_h &= \bm{v}_h, \\
    \inpr[\Omega]{\bm{\psi}}{\rho \partial_t \bm{v}_h} &= - \inpr[\Omega]{\mathcal{L}(\mathcal{D}\bm{q}_h)\bm{\psi}}{\bm{S}_h}, \\
    \inpr[\Omega]{\bm{\Psi}}{\bm{C} \partial_t \bm{S}_h} &= +\inpr[\Omega]{\bm{\Psi}}{\mathcal{L}(\mathcal{D}\bm{q}_h)\bm{v}_h}, 
\end{aligned}
\end{equation}
forall $\bm{\psi} \in V_h$ and $\bm{\Psi} \in \Sigma_h$. If the problem were linear (meaning that the operator $\mathcal{L}$ would not depend on the displacement), then one could consider spaces $V_h, \; \Sigma_h$ that would respect the inclusion $\mathcal{L} V_h \subset \Sigma_h$. In this way the second line of the system \eqref{eq:weak_abstract} would be satisfied pointwise. This is clearly not true in a nonlinear context. So a choice will be made in order to enrich the polynomial space for the stress, so that the nonlinear terms can be faithfully represented. Once a finite element basis has been selected, the weak formulation is converted into the following algebraic system
\begin{equation*}
\begin{aligned}
\dot{\mathbf{q}} &= \mathbf{v}, \\
\begin{bmatrix}
    \mathbf{M}_\rho & 0 \\
    0 & \mathbf{M}_{C}
\end{bmatrix}
\odv{}{t}
\begin{pmatrix}
    \mathbf{v} \\
    \mathbf{s}\\
\end{pmatrix} &= 
\begin{bmatrix}
    0 & -\mathbf{L}^\top(\mathbf{q}) \\
    \mathbf{L}(\mathbf{q}) & 0 \\
\end{bmatrix}
\begin{pmatrix}
    {\mathbf{v}} \\
    {\mathbf{s}} \\
\end{pmatrix}. 
\end{aligned}
\end{equation*}
\rewtwo{The matrices $\mathbf{M}_\rho, \; \mathbf{M}_{C}$ 
are mass matrices arising from the finite element formulation. Their expression is given by 
\begin{equation*}
\begin{aligned}
    [\mathbf{M}_{\rho}]_{ij} &= \int_0^L \rho \bm{\varphi}_i \cdot \bm{\varphi}_j \; \d{x}, \\
    [\mathbf{M}_{C}]_{ij} &= \int_0^L  \bm{\Xi}_i : \bm{C} \bm{\Xi}_j \;\d{x}, 
\end{aligned}
\end{equation*}
where $\bm{\varphi}_i, \;  \bm{\Xi}_i$ represent the finite element basis chosen for $\bm{v}$ and $\bm{S}$. The coefficients $\mathbf{v}, \; \mathbf{s}$ are the degrees of freedom the velocity and stress variable.
}
The algebraic system can be written more compactly in the following form
\begin{equation}\label{eq:alg_system_abstract}
    \begin{aligned}
\dot{\mathbf{q}} &= \mathbf{v}, \\
\mathbf{H}
\dot{\mathbf{x}} &= 
\mathbf{J}(\mathbf{q})\mathbf{x}.
\end{aligned}
\end{equation}
The energy of the system is given by $H = \frac{1}{2}\mathbf{x}^\top \mathbf{H}\mathbf{x}$, in complete analogy with the nonlinear oscillator example.

\subsubsection{Von-K\'arm\'an beam}
The physical domain for this example is an interval $\Omega=[0, L]$. Denote by $E$ a generic element in a mesh $\mathcal{I}_h$. For this example the space $V_h$ contains two  finite element spaces, the longitudinal and vertical deflections/velocities. These are discretized using linear and cubic (Hermite) polynomials respectively \cite{deborst2012nonlinear}
\begin{equation*}
\begin{aligned}
    V_h &= \mathrm{CG}_1 \times \mathrm{Her}, \\
    \mathrm{CG}_1 &:= \{v_h \in C^0([0, L]), \; v_h|_E \in \bbP_1, \; \forall E \in \mathcal{I}_h \}, \\
    \mathrm{Her} &:= \{v_h \in C^1([0, L]), \; v_h|_E \in \bbP_3, \; \forall E \in \mathcal{I}_h \}, 
\end{aligned}
\end{equation*}
where $ \mathrm{CG}_1$ is the linear Lagrange space and $\mathrm{Her}$ is the space of Hermite polynomials. The space $\Sigma_h$ contains the axial and bending stress resultant. These are discretized using quartic and linear discontinuous shape functions respectively
\begin{equation}
    \begin{aligned}
        \Sigma_h &= \mathrm{DG}_4 \times \mathrm{DG}_1, \\
        \mathrm{DG}_{k} &= \{v_h|_{E} \in \bbP_{k}, \; \forall E \in \mathcal{I}_h \}, 
    \end{aligned}
\end{equation}
where $ \mathrm{DG}_{k}$ is the Discontinuous Galerkin of order $k$. The quartic choice is due to the fact that the axial stress is proportional to the square power of the derivative of the vertical displacement 
$$N = EA\left(\partial_x q_x + \frac{1}{2} (\partial_x q_z)^2\right).$$
Since the vertical displacement is discretized via Hermite cubic polynomial, its squared derivative is a continuous quartic polynomial. Selecting a discontinuous quartic polynomial guarantees a correct representation of the axial stress, while choosing a smaller discretization space would lead to locking phenomena The weak formulation for this problem then becomes: find $v_{x, h} \in \mathrm{CG}_1, \; (q_{z, h},\,v_{z, h}) \in \mathrm{Her}, \; N_h \in \mathrm{DG}_4, \; M_h \in \mathrm{DG}_1$ such that
\begin{equation}\label{eq:weak_form_vK_beam}
\begin{aligned}
\partial_t q_{z, h} &= v_{z, h}, \\
\inpr[\Omega]{\psi_x}{\rho A \, \partial_t{v}_{x, h}} = &-\inpr[\Omega]{\partial_x \psi_{x}}{N_h}. \\
\inpr[\Omega]{\psi_z}{\rho A\,\partial_t {v}_{z, h}} = &-\inpr[\Omega]{\partial_x \psi_z\,\partial_x q_{z, h}}{N_h}  \\
&-\inpr[\Omega]{\partial_{xx} \psi_z}{M_h}, \\
\inpr[\Omega]{\psi_N}{C_a \, \partial_t N_h} = &+ \inpr[\Omega]{\psi_N}{\partial_x q_z^h\,\partial_x v_z^h} \\
&+ \inpr[\Omega]{\psi_N}{\partial_{x} v_{x, h}}, \\
\inpr[\Omega]{\psi_M}{C_b \, \partial_t M_h} = &+\inpr[\Omega]{\psi_M}{\partial_{xx} v_z^h}, 
\end{aligned} 
\end{equation}
for all $\psi_x \in \mathrm{CG}_1, \; \psi_z \in \mathrm{Her}, \; \psi_N \in \mathrm{DG}_0, \; \psi_M \in \mathrm{DG}_1$.
The algebraic system arising from the finite element discretization has the same form as in \eqref{eq:alg_system_abstract} with 
$\mathbf{x} := (\mathbf{v}^\top_x \; \mathbf{v}^\top_z \; \mathbf{n}^\top \; \mathbf{m}^\top)^\top$ and
\begin{equation*}
\begin{aligned}
\mathbf{H}&:= 
\mathrm{Diag}
\begin{bmatrix}
\mathbf{M}_{\rho A} \\
\mathbf{M}_{\rho A} \\
\mathbf{M}_{C_a} \\
\mathbf{M}_{C_b}  \\
\end{bmatrix}, \\
\mathbf{J}&:=
\begin{bmatrix}
    0 & 0 & -\mathbf{D}_{\partial_x}^\top & 0 \\
    0 & 0 & -\mathbf{L}^\top(\mathbf{q}_z) & -\mathbf{D}_{\partial_{xx}}^\top \\
    \mathbf{D}_{\partial_x} & \mathbf{L}(\mathbf{q}_z) & 0 & 0  \\
    0 & \mathbf{D}_{\partial_{xx}} & 0 & 0 \\
\end{bmatrix}.
\end{aligned}
\end{equation*}

\subsubsection{Geometrically nonlinear elasticity}
Denote with $T$ a generic cell of the computational mesh $\mathcal{T}_h$. For geometrically nonlinear elasticity, displacement and velocity are discretized using Lagrange polynomials of order~$1$. The stress is instead discretized using a tensor-valued symmetric discontinuous space of polynomials of degree~$0$
\begin{equation*}
    \begin{aligned}
    V_h &= \mathrm{CG}_1(\bbR^d), \\ 
    \Sigma_h &= \mathrm{DG}_{0}(\bbR^{d \times d}_{\rm sym}) = \{\bm{S}_h|_{T} \in [\bbP_{0}]^{d \times d}_{\mathrm{sym}}, \; \forall T \in \mathcal{T}_h \}.
    \end{aligned}
\end{equation*}
where $d= \{2,3\}$ is the geometric dimension of the problem. The resulting discrete formulation reads
\begin{equation*}
\begin{aligned}
    \partial_t \bm{q}_h &= \bm{v}_h, \\
    \inpr[\Omega]{\bm{\psi}}{\rho\,\partial_t \bm{v}_h} &= - \inpr[\Omega]{\bm{F}_h^\top\nabla\bm{\psi}}{\bm{S}_h},  \\
    \inpr[\Omega]{\bm{\Psi}}{\bm{C}\,\partial_t \bm{S}_h} &= +\inpr[\Omega]{\bm{\Psi}}{\bm{F}_h^\top\nabla\bm{v}_h},  \\
\end{aligned}
\end{equation*}
forall $\bm{\psi} \in V_h$ and $\bm{\Psi} \in \Sigma_h$. \rewother{The symmetrization operator has been omitted as the inner product with a symmetric tensor naturally enforces  symmetry.} The ordinary differential equation arising from the discretization has the same form as \eqref{eq:alg_system_abstract}.

\subsection{Time integration method}\label{sec:time_int}

For the time integration the combination of two well-known symplectic methods is considered: the St\"ormer-Verlet (or leapfrog) scheme and the implicit midpoint method. These schemes are particular instances of the Newmark method in structural mechanics. Consider the following system 
$$
\begin{aligned}
\dot{\mathbf{q}} &= \mathbf{v}, \\
\mathbf{M}\dot{\mathbf{v}} &= \mathbf{f}(\mathbf{q}).
\end{aligned}
$$
If $\mathbf{f}(\mathbf{q}) = -\nabla_{\mathbf{q}} V$, then the system represents a canonical Hamiltonian formulation. The St\"ormer-Verlet method (also called leapfrog method, central difference or explicit Newmark in mechanics) can be then written in the following form
$$
\begin{aligned}
    \frac{\mathbf{q}_{n+\frac{1}{2}} - \mathbf{q}_{n-\frac{1}{2}}}{\Delta t} &=  \mathbf{v}_{n}, \\
    \mathbf{M}\frac{(\mathbf{v}_{n+1}-\mathbf{v}_{n})}{\Delta t} &=  \mathbf{f}(\mathbf{q}_{n+\frac{1}{2}}).
\end{aligned}
$$
To start the iterations, the first value for the velocity is obtained via a second order Taylor approximation
\begin{equation}\label{eq:init_displacement}
    \mathbf{q}_{\frac{1}{2}} = \mathbf{q}_0 + \frac{\Delta t}{2} \mathbf{v}_{0} + \frac{1}{8}\Delta t^2\mathbf{a}_0, \qquad \mathbf{a}_0 := \mathbf{M}^{-1}\mathbf{f}(\mathbf{q}_0).
\end{equation}
The implicit midpoint method applied to a generic ODE of the form $\dot{\mathbf{x}} = \mathbf{g}(\mathbf{x})$ gives 
$$
\frac{\mathbf{x}_{n+1} - \mathbf{x}_n}{\Delta t} = \mathbf{g}\left(\mathbf{x}_{n+\frac{1}{2}}\right),
$$
with $ \mathbf{x}_{n+\frac{1}{2}} : = (\mathbf{x}_{n+1} + \mathbf{x}_n)/2$. The idea is to apply these two methods together to the special form given by system \eqref{eq:alg_system_abstract}, leading to the following relations 
\begin{equation}\label{eq:time_int_scheme}
\begin{aligned}
\frac{\mathbf{q}_{n+\frac{1}{2}} -\mathbf{q}_{n-\frac{1}{2}}}{\Delta t} &= \mathbf{v}_{n}, \\
\mathbf{H}\frac{(\mathbf{x}_{n+1} - \mathbf{x}_{n})}{\Delta t}&= 
\mathbf{J}(\mathbf{q}_{n+\frac{1}{2}})\frac{(\mathbf{x}_{n+1} + \mathbf{x}_{n})}{2}
\end{aligned}
\end{equation}
The initial position $\mathbf{q}_{n+\frac{1}{2}}$ is again computed via a second order Taylor expansion as in \eqref{eq:init_displacement}.

\rewtwo{
\begin{remark}[Imposition of the Dirichlet boundary conditions]
For the linearly implicit integration, the boundary conditions are imposed on the velocity field only. The enforcement of the boundary conditions at the displacement follows automatically as the initial condition satisfy the boundary condition and the next displacement is computed via update rule 
$$
\mathbf{q}_{n+1/2} = \mathbf{q}_{n-1/2} + \Delta t \mathbf{v}_{n}.
$$
Since $\mathbf{q}_{n-1/2}|_{\Gamma_D} = 0, \; \mathbf{v}_{n}|_{\Gamma_D} = 0$, then it follows that $\mathbf{q}_{n+1/2}|_{\Gamma_D} = 0$.
\end{remark}
}

\subsection{Properties of the scheme}

The proposed time discretization scheme conserves the energy and angular momentum exactly.

\begin{proposition}
    The discrete energy 
    $$H_n := \frac{1}{2} \mathbf{x}_{n}^\top\mathbf{H} \mathbf{x}_{n}$$
    is conserved by the scheme.
\end{proposition}
\begin{proof}
 This is shown by considering the scalar multiplication of the second equation by $\mathbf{x}_{n+\frac{1}{2}}$
$$
\mathbf{x}_{n+\frac{1}{2}}^\top\mathbf{H} \frac{(\mathbf{x}_{n+1} - \mathbf{x}_{n})}{\Delta t} = \mathbf{x}_{n+\frac{1}{2}}^\top\mathbf{J}(\mathbf{q}_{n+\frac{1}{2}})\mathbf{x}_{n+\frac{1}{2}}= 0,
$$
by the skew-symmetry of $\mathbf{J}(\mathbf{q}_{n+1})$. This implies
\begin{equation}\label{eq:energy_conser}
    \mathbf{x}_{n+1}^\top\mathbf{H}\mathbf{x}_{n+1} = \mathbf{x}_{n}^\top\mathbf{H}\mathbf{x}_{n}, \qquad \text{i.e.} \quad H_{n+1} = H_{n}.
\end{equation}
\end{proof}
\begin{proposition}\label{pr:invertibility}
The discrete dynamics \eqref{eq:time_int_scheme} can be rewritten in a recursive form as follows
\begin{equation}
\begin{pmatrix}
   \mathbf{q}_{n+\frac{1}{2}} \\
   \mathbf{x}_{n+1}
\end{pmatrix} =
\begin{bmatrix} 
    \mathbf{I} & \Delta t\mathbf{B} \\
    0 & \displaystyle \mathrm{Cay}\left(\frac{\Delta t}{2}\mathbf{H}^{-1}\mathbf{J}\right)
\end{bmatrix}
\begin{pmatrix}
   \mathbf{q}_{n-\frac{1}{2}} \\
   \mathbf{x}_{n}
\end{pmatrix},
\end{equation}
where $\mathbf{B}\,\mathbf{x}_n:=\mathbf{v}_n$ and the Cayley transform have been introduced
$$
\mathrm{Cay}(\mathbf{M}):= (\mathbf{I} - \mathbf{M})^{-1} (\mathbf{I} + \mathbf{M}).
$$    
\end{proposition}
\begin{proof}
Consider system \eqref{eq:time_int_scheme} and the recursion rule 
$$
\left[\mathbf{H} - \frac{\Delta t}{2}\mathbf{J} \right]\mathbf{x}_{n+1} = \left[\mathbf{H} - \frac{\Delta t}{2}\mathbf{J} \right]\mathbf{x}_{n}.
$$
where the explicit dependence on $\mathbf{q}$ for the matrix $\mathbf{J}$ has been removed for simplicity. The matrix $\left[\mathbf{H} - \frac{\Delta t}{2}\mathbf{J} \right]$ is invertible. To prove this, consider the Cholesky factorization $\mathbf{H}=\mathbf{Q}^\top\mathbf{Q}$. Then it holds
\begin{align*}
    \mathbf{Q}^\top\mathbf{Q} - \frac{\Delta t}{2}\mathbf{J} = \mathbf{Q}^\top \left[\mathbf{I} - \frac{\Delta t}{2}\mathbf{Q}^{-\top}\mathbf{J}\mathbf{Q}^{-1} \right]  \mathbf{Q}.
\end{align*}
The matrix $\displaystyle \mathbf{Q}^{-\top}\mathbf{J}\mathbf{Q}^{-1}$ is skew-symmetric. Therefore $\displaystyle \mathbf{I} - \frac{\Delta t}{2}\mathbf{Q}^{-\top}\mathbf{J}\mathbf{Q}^{-1}$ is invertible and consequently $\displaystyle \mathbf{H} - \frac{\Delta t}{2}\mathbf{J}$ since it is the product of three invertible matrices. The state transition matrix $\mathbf{x}_{n+1} = \mathbf{A}(\mathbf{q}_{n+\frac{1}{2}})\mathbf{x}_{n}$ is then obtained as 
$$
\begin{aligned}
\mathbf{A}(\mathbf{q}_{n+\frac{1}{2}})&:= \left[\mathbf{H} - \frac{\Delta t}{2}\mathbf{J} \right]^{-1}\left[\mathbf{H} - \frac{\Delta t}{2}\mathbf{J}\right], \\ 
&= \mathrm{Cay}\left(\frac{\Delta t}{2}\mathbf{H}^{-1}\mathbf{J}\right),
\end{aligned}
$$
where the equality follows from factorizing $\mathbf{H}$ from both terms and using the inverse of the product.
\end{proof}
The energy conservation \eqref{eq:energy_conser} is then rewritten as
$$
||\mathbf{x}||_{\mathbf{H}} = ||\mathrm{Cay}\left(\mathbf{H}^{-1}\mathbf{J}\mathbf{x}\right)||_{\mathbf{H}}, 
$$
where $||\mathbf{x}||_{\mathbf{H}}:=\mathbf{x}^\top\mathbf{H}\mathbf{x}$ is the norm induced by the positive definite symmetric matrix $\mathbf{H}$. This means that the Cayley transform of $\mathbf{H}^{-1}\mathbf{J}$ is a unitary matrix in the $\mathbf{H}$ norm and therefore has eigenvalues lying on the unit circle. Given its block upper triangular structure, the overall state transition matrix has a spectrum given by the union of the spectrum of its diagonal blocks $\mathbf{I}$ and $\mathrm{Cay}(\frac{\Delta t}{2}\mathbf{H}^{-1}\mathbf{J}(\mathbf{q}_{n+\frac{1}{2}}))$. Therefore all its eigenvalues lie on the unit circle. 

\rewtwo{
\begin{proposition}
When no Dirichlet boundary conditions apply, the method preserves angular momentum.     
\end{proposition}
\begin{proof}
    The fact that the linear implicit scheme preserves the angular momentum comes from the employment of the St\"ormer-Verlet integrator for the dynamics of the displacement and velocity. At the continuous time level, the preservation of the angular momentum can be proved using the argument presented in \cite[Section 2.3.2]{simo1992conserving}. Consider the weak form of the proposed scheme
\begin{equation*}
\begin{aligned}
    \partial_t \bm{q}_h &= \bm{v}_h, \\
    \inpr[\Omega]{\bm{\psi}}{\rho\,\partial_t \bm{v}_h} &= - \inpr[\Omega]{\bm{F}_h^\top\nabla\bm{\psi}}{\bm{S}_h},  \\
\end{aligned}
\end{equation*}
The dynamics of the stress is omitted as it does not play any role in the argument. Let $[\bm{\xi}]_{\times}$ denotes the skew-symmetric matrix arising from an $\mathbb{R}^3$-vector
\begin{equation*}
[\bm{\xi}]_\times :=
\begin{bmatrix}
0 & - \xi_3 & \xi_2 \\
\xi_3 & 0 & - \xi_1 \\
-\xi_2 & \xi_1 & 0
\end{bmatrix}.
\end{equation*}
For the linear momentum balance consider the test function to be $\bm{\psi} = [\bm{e}_i]_{\times} \bm{r}_h$ where $\bm{e}_i, \; i=\{1,2,3\}$ is the canonical basis of $\mathbb{R}^3$ and $\bm{r}_h = \bm{X} + \bm{q}_h$ denotes the position vector ($\bm{X}$ is the position of the undeformed configuration). This choice of test function is allowed only if no Dirichlet boundary condition is imposed. A direct computation gives 
$$\nabla ([\bm{e}_i]_{\times} \bm{r}_h) = [\bm{e}_i]_{\times} \bm{F}_h,$$ 
where $\bm{F}_h:= \bm{I} + \nabla \bm{q}_h$ is the deformation. For the $i$-th component of the angular momentum $J_i$ is obtained
\begin{equation*}
\begin{aligned}
\odv{J_i}{t} &= \inpr[\Omega]{[\bm{e}_i]_{\times} \bm{r}_h}{\rho\,\partial_t \bm{v}_h}, \\
     &= - \inpr[\Omega]{\bm{F}_h^\top\nabla([\bm{e}_i]_{\times} \bm{r}_h)}{\bm{S}_h}, \\
    &= - \inpr[\Omega]{\bm{F}_h^\top [\bm{e}_i]_{\times} \bm{F}_h}{\bm{S}_h} = 0,
\end{aligned}
\end{equation*}
since the product of a skew-symmetric  $\bm{F}_h^\top [\bm{e}_i]_{\times} \bm{F}_h$ and symmetric tensor $\bm{S}_h$ vanishes. For the time discrete conservation, we follow the same argument as in \cite[Theorem 3.5]{hairer2003verlet}, using the weak form of the equation. The leapfrog method is rewritten as the composition of two symplectic Euler methods SE1, SE2
\begin{subequations}
\begin{align}
\bm{q}_{n+1/2} &= \bm{q}_{n} + \frac{\Delta t}{2} \bm{v}_h^n, \tag{SE1 [a]} \\
\inpr[\Omega]{\bm{\psi}}{\rho\,(\bm{v}_{n+1/2}-\bm{v}_{n})} &= - \frac{\Delta t}{2}\inpr[\Omega]{\bm{F}_{n+1/2}^\top\nabla\bm{\psi}}{\bm{S}}, \tag{SE1 [b]} \\
\inpr[\Omega]{\bm{\psi}}{\rho\,(\bm{v}_{n+1}-\bm{v}_{n+1/2})} &= - \frac{\Delta t}{2}\inpr[\Omega]{\bm{F}_{n+1/2}^\top\nabla\bm{\psi}}{\bm{S}}, \tag{SE2 [a]} \\
\bm{q}_{n+1} &= \bm{q}_{n+1/2} + \frac{\Delta t}{2} \bm{v}_h^n. \tag{SE2 [b]}
\end{align}
\end{subequations}
The actual point at  which the stress is evaluated does not play any role in the proof, so we omit it. We consider only the $\mathrm{SE}1$ part as the proof for the second part is analogous. Adding $\bm{X}$ to both side of the displacement update leads to 
\begin{equation}\label{eq:update_position}
\bm{r}_{n+1/2} = \bm{r}_n + \frac{\Delta t}{2} \bm{v}_n.  
\end{equation}
Choosing $\bm{\psi} = [\bm{e}_i]_\times \bm{r}_{n+1/2}$ in the linear momentum balance provides
\begin{multline*}
\inpr[\Omega]{[\bm{e}_i]_\times \bm{r}_{n+1/2}}{\rho\,(\bm{v}_{n+1/2}-\bm{v}_{n})} \\
= - \frac{\Delta t}{2}\inpr[\Omega]{\bm{F}_{n+1/2}^\top\nabla([\bm{e}_i]_\times \bm{r}_{n+1/2})}{\bm{S}}\,, \\
= - \frac{\Delta t}{2}\inpr[\Omega]{\bm{F}_{n+1/2}^\top[\bm{e}_i]_\times \bm{F}_{n+1/2}}{\bm{S}} = 0.
\end{multline*}
The last equality is again due to the fact that the inner product of symmetric and skew-symmetric tensor vanishes. Using the update of the position \eqref{eq:update_position}, it is obtained
$$
\begin{aligned}
J_i^{n+1/2} &= \inpr[\Omega]{[\bm{e}_i]_\times \bm{r}_{n+1/2}}{\rho\,\bm{v}_{n+1/2}}, \\
 &= \inpr[\Omega]{[\bm{e}_i]_\times \bm{r}_{n+1/2}}{\rho\,\bm{v}_{n}}, \\
 &= \inpr[\Omega]{[\bm{e}_i]_\times \bm{r}_{n}}{\rho\,\bm{v}_{n}} + \frac{\Delta t}{2}\inpr[\Omega]{[\bm{e}_i]_\times \bm{v}_{n}}{\rho\,\bm{v}_{n}}, \\
 &= J_i^{n}
\end{aligned}
$$
since $\inpr[\Omega]{[\bm{e}_i]_\times \bm{v}_{n}}{\rho\,\bm{v}_{n}}=0$. The same argument holds for $\mathrm{SE}2$ hence $$J_i^{n+1} = J_i^n.$$
\end{proof}
}

\subsection{Static condensation}
Since the finite element space for the stress variable is discontinuous, the associated mass matrix is block-diagonal. To speed up the solver, this variable can be statically condensed. In the following the midpoint values of the velocity and stress are denoted by
$$
\mathbf{v}_{n+\frac{1}{2}}:=\frac{\mathbf{v}_{n+1}+\mathbf{v}_{n}}{2}, \qquad \mathbf{s}_{n+\frac{1}{2}}:=\frac{\mathbf{s}_{n+1}+\mathbf{s}_{n}}{2}.
$$
From system \eqref{eq:time_int_scheme}, the velocity and stress discrete system reads
\begin{align}
    \mathbf{M}_\rho \frac{(\mathbf{v}_{n+1}-\mathbf{v}_{n})}{\Delta t} &= -\mathbf{L}^\top(\mathbf{q}_{n+\frac{1}{2}})\mathbf{s}_{n+\frac{1}{2}}, \label{eq:time_int_velocity}\\
    \mathbf{M}_C \frac{(\mathbf{s}_{n+1}-\mathbf{s}_{n})}{\Delta t} &= +\mathbf{L}(\mathbf{q}_{n+\frac{1}{2}})\mathbf{v}_{n+\frac{1}{2}}. \label{eq:time_int_stress}
\end{align}
The expression of the midpoint value for the stress $\mathbf{s}_{n+\frac{1}{2}}$ is given by
$$
\mathbf{s}_{n+\frac{1}{2}} = \mathbf{s}_{n} + \frac{\Delta t}{2} \mathbf{M}_C^{-1} \mathbf{L}(\mathbf{q}_{n+\frac{1}{2}})\mathbf{v}_{n+\frac{1}{2}},
$$
Replacing this value into \eqref{eq:time_int_velocity} leads to
$$
\begin{aligned}
\mathbf{M}_\rho \frac{(\mathbf{v}_{n+1}-\mathbf{v}_{n})}{\Delta t} = &- \frac{\Delta t}{2}\mathbf{K}(\mathbf{q}_{n+\frac{1}{2}})\mathbf{v}_{n+\frac{1}{2}} \\
&- \mathbf{L}^\top(\mathbf{q}_{n+\frac{1}{2}})\mathbf{s}_{n},  
\end{aligned}
$$
where $\mathbf{K}(\mathbf{q}_{n+\frac{1}{2}}) := \mathbf{L}^\top(\mathbf{q}_{n+\frac{1}{2}}) \mathbf{M}_C^{-1} \mathbf{L}(\mathbf{q}_{n+\frac{1}{2}})$.
The time integration recursion is then expressed as: given $(\mathbf{q}_{n-\frac{1}{2}}, \; \mathbf{v}_{n}, \; \mathbf{s}_{n})$, the values for the next time step $(\mathbf{q}_{n+\frac{1}{2}}, \; \mathbf{v}_{n+1},\; \mathbf{s}_{n+1})$ are given by 
    \begin{equation}
    \begin{aligned}
\mathbf{q}_{n+\frac{1}{2}} =& \mathbf{q}_{n-\frac{1}{2}} + \Delta t\, \mathbf{v}_{n}, \\
    \left[\mathbf{M}_\rho + \frac{\Delta t^2}{4} \mathbf{K}\right]
    \mathbf{v}_{n+1}
=&  \left[\mathbf{M}_{\rho} - \frac{\Delta t^2}{4} \mathbf{K}\right]\mathbf{v}_{n} \\
&- \Delta t \mathbf{L}^\top \mathbf{s}_{n}, \\
\mathbf{s}_{n+1} =& \mathbf{s}_{n} + \Delta t\mathbf{M}_C^{-1} \mathbf{L}\frac{(\mathbf{v}_{n+1} + \mathbf{v}_{n})}{2},
\end{aligned} 
\end{equation}
where the explicit dependence of matrices $\mathbf{L}, \; \mathbf{K}$ from $\mathbf{q}_{n+\frac{1}{2}}$ has been suppressed for simplicity.

\subsection{Connection with the scalar auxiliary variable approach}

The scalar auxiliary variable approach has been first introduced in \cite{shen2018sav} in the context of gradient flows and then extended to the case of Hamiltonian dynamics in \cite{bilbao2023explicit}. Therein the authors apply the approach to Hamiltonian system with diagonal mass matrix, thus obtaining a scheme that is both computationally efficient and energy stable. To illustrate the connection with the scalar auxiliary method and the present framework consider a separable Hamiltonian system in canonical form
\begin{equation*}
    \begin{bmatrix}
        \mathbf{I} & 0 \\
        0 & \mathbf{M}
    \end{bmatrix}
    \odv{}{t}
    \begin{pmatrix}
        \mathbf{q} \\
        \mathbf{v}
    \end{pmatrix} = \begin{bmatrix}
        0 & \mathbf{I} \\
        -\mathbf{I} & 0
    \end{bmatrix}
    \begin{pmatrix}
        \nabla_\mathbf{q} V \\
        \mathbf{v}
    \end{pmatrix}.
\end{equation*}
The dynamics is written in terms of the velocity instead of the linear momentum to avoid taking the inverse of the mass matrix in a finite element formulation. Assume that the potential energy is positive $V\ge 0$ (as in the case of continuum mechanics). The non negativity condition allows writing the potential energy as
$$
V = \frac{1}{2} \xi^2.
$$
The gradient of the potential can then be written in terms of $\xi$ as 
$$
\nabla_\mathbf{q} V = \xi \nabla_{\mathbf{q}} \xi.
$$
The chain rule provides
$$
\dot{\xi} = (\nabla_{\mathbf{q}} \xi)^\top \dot{\mathbf{q}} = (\nabla_{\mathbf{q}} \xi)^\top \mathbf{v}.
$$
The dynamics can then be rewritten in terms of $\mathbf{q}, \mathbf{v}, \xi$ in non-canonical Hamiltonian form as follows
\begin{equation*}
\begin{aligned}
\dot{\mathbf{q}} &= \mathbf{v}, \\
    \begin{bmatrix}
        \mathbf{M} & 0 \\
        0 & 1
    \end{bmatrix}
    \odv{}{t}
    \begin{pmatrix}
        \mathbf{v} \\
        \xi
    \end{pmatrix} &= \begin{bmatrix}
        0 & - \mathbf{g}(\mathbf{q})\\
        \mathbf{g}(\mathbf{q})^\top & 0
    \end{bmatrix}
    \begin{pmatrix}
        \mathbf{v} \\
        \xi
    \end{pmatrix},
\end{aligned}
\end{equation*}
where $\mathbf{g}(\mathbf{q}) = \nabla_{\mathbf{q}} \xi$. Notice that the system can be rewritten as follows
\begin{equation*}
\begin{aligned}
\dot{\mathbf{q}} &= \mathbf{v}, \\
\mathbf{H} \dot{\mathbf{x}} &= \mathbf{J}(\mathbf{q}) \mathbf{x}, \qquad \mathbf{x} :=  [\mathbf{v}^\top \; \xi]^\top,
\end{aligned}
\end{equation*}
and the total energy is given by $H = \frac{1}{2}\mathbf{x}^\top \mathbf{H} \mathbf{x}$. This is in complete analogy with the non-canonical formulation presented in the present work. In \cite[Sec. 3.3]{bilbao2023explicit} the time integration is presented as follows
$$
\begin{aligned}
    \frac{\mathbf{q}_{n+1} - \mathbf{q}_{n}}{\Delta t} &= \mathbf{v}_{n+\frac{1}{2}} , \\
   \mathbf{M}\frac{\mathbf{v}_{n+\frac{1}{2}} -  \mathbf{v}_{n-\frac{1}{2}}}{\Delta t} &= -\mathbf{g}(\mathbf{q}_n)\frac{(\xi_{n+\frac{1}{2}} + \xi_{n-\frac{1}{2}})}{2}, \\
   \xi_{n+\frac{1}{2}}-\xi_{n-\frac{1}{2}} &= +\mathbf{g}(\mathbf{q}_{n})^\top \frac{(\mathbf{q}_{n+1}+\mathbf{q}_{n-1})}{2}.
\end{aligned}
$$
By the first line it holds 
$$\frac{\mathbf{q}_{n+1}+\mathbf{q}_{n-1}}{2\Delta t} = \frac{\mathbf{v}_{n+\frac{1}{2}}+\mathbf{v}_{n-\frac{1}{2}}}{2}.$$
Therefore the scheme is equivalently rewritten as follows
$$
\begin{aligned}
    \frac{\mathbf{q}_{n+1} - \mathbf{q}_{n}}{\Delta t} &= \frac{\mathbf{v}_{n+1} + \mathbf{v}_n}{2}, \\
    \mathbf{H}\frac{(\mathbf{x}_{n+\frac{1}{2}}-\mathbf{x}_{n-\frac{1}{2}})}{\Delta t} &= \mathbf{J}(\mathbf{q}_{n}) \frac{\mathbf{x}_{n+\frac{1}{2}}+\mathbf{x}_{n-\frac{1}{2}}}{2}.
\end{aligned}
$$

This recursion is simply a staggered version of \eqref{eq:time_int_scheme}. 

\section{Numerical results}\label{sec:num_results}
In this section, we test the proposed methodology for three different cases: 
\begin{itemize}
    \item the Duffing oscillator;
    \item vibrations of a von-K\'arm\'an beam;
    \item bending of a column in geometrically nonlinear elasticity.
\end{itemize}
The linearly implicit scheme is compared against the St\"ormer-Verlet method and the exact energy conserving scheme presented in \cite{simo1992conserving}.
The former scheme is symplectic, energy preserving only in an approximate sense and explicit, thus it requires the fulfillment of a CFL-like condition. The latter method is implicit and exactly energy preserving and unconditionally stable. It requires however an iterative procedure to solve the resulting nonlinear system. 
\rewtwo{The idea of momentum and energy preserving algorithm dates back to \cite[Section 3.2.2.]{simo1992conserving}. The idea applies  to generic constitutive laws. To illustrate the idea, consider a mechanical system of the form
$$
\begin{aligned}
\dot{\mathbf{q}} &= \mathbf{v}, \\
\mathbf{M}\dot{\mathbf{v}} &= -\mathbf{L}^\top(\mathbf{q}) \bm{\sigma}
\end{aligned}
$$
where $\bm{\sigma}$ is the stress. This system arises from a classical finite element discretization of geometrically nonlinear mechanical models considered in this paper. 
If only geometrical nonlinearity are considered, the deformation energy is a quadratic form in the strain
$$
V_{\rm def} = \frac{1}{2} \bm{\varepsilon}^\top \mathbf{W} \bm{\varepsilon}
$$
The stress variable is given by
$$
\bm{\sigma} = \pdv{V_{\rm def}}{\bm\varepsilon} = \mathbf{W}\bm{\varepsilon}
$$
The strain is a nonlinear function of the displacement, i.e. $\bm{\varepsilon} = \mathbf{G}(\mathbf{q})$ and so the stress
$$
\bm{\sigma} = \mathbf{W} \mathbf{G}(\mathbf{q}).
$$
Now to obtain exact energy conservation in the case of Saint Venant-Kirchhoff material it is sufficient to use the midpoint rule with a slight modification. Instead of evaluating the stress at the midpoint, the average of the previous and next value is taken. Following the notation used in the paper for the weak formulation (cf. \cite[Section 3.3.1.]{simo1992conserving})
\begin{equation*}
    \widehat{\bm{\sigma}}_{n+1/2} := \mathbf{W}\frac{\mathbf{G}(\mathbf{q}_{n+1}) + \mathbf{G}(\mathbf{q}_{n})}{2}.
\end{equation*}
This leads to the following discretization
$$
\begin{aligned}
\frac{\mathbf{q}_{n+1} - \mathbf{q}_n}{\Delta t} &= \frac{\mathbf{v}_{n+1/2} + \mathbf{v}_{n}}{2}, \\
\mathbf{M}\frac{\mathbf{v}_{n+1} - \mathbf{v}_n}{\Delta t} &= -\mathbf{L}^\top\left(\frac{\mathbf{q}_{n+1} + \mathbf{q}_n}{2}\right) \widehat{\bm{\sigma}}_{n+1/2}.
\end{aligned}
$$
As an illustrative example consider the Duffing oscillator
$$
\ddot{q} = -\alpha q - \beta q^3.
$$
In this case the deformation (or potential) energy is given by
$$
V = \frac{1}{2}\alpha q^2 + \frac{1}{4} \beta q^4.
$$
As shown in the introductory example in Section~\ref{sec:models}, the strain field is two dimensional for this example
$$
\varepsilon_1 = q, \qquad \varepsilon_2 = q^2.
$$
The expression of the energy in terms of this stress is given by
$$
V_{\rm def} = \frac{1}{2} \begin{pmatrix}
    \varepsilon_1 \\
    \varepsilon_2
\end{pmatrix}^\top
\begin{bmatrix}
    \alpha & 0 \\
    0 & \beta/2
\end{bmatrix}
\begin{pmatrix}
    \varepsilon_1 \\
    \varepsilon_2
\end{pmatrix}.
$$
So the stresses are given by
$$
\sigma_1 = \alpha q, \qquad \sigma_2 = \frac{\beta}{2} q^2.
$$
The average value of $\sigma_1$ coincide with the application of the midpoint scheme. This is not surprising as this is the linear part of the equation and the implicit midpoint exactly preserves the energy for linear systems. For $\sigma_2$ the average is given by
$$
\widehat{\sigma_2}^{n+1/2} =\beta \frac{q^2_{n+1} + q^2_{n}}{4}.
$$
So the discrete scheme for the Duffing oscillator becomes
$$
\begin{aligned}
\frac{{q}_{n+1} - {q}_n}{\Delta t} &= \frac{{v}_{n+1} + {v}_{n}}{2}, \\
\frac{{v}_{n+1} - {v}_n}{\Delta t} &= - \alpha q_{n+1/2} - 2 q_{n+1/2} \left(\beta \frac{q^2_{n+1} + q^2_{n}}{4}\right).
\end{aligned}
$$
where $\displaystyle q_{n+1/2}:=\frac{q_{n+1} + q_{n}}{2}$. It can be verified that this final scheme coincides with the employment of the mean value discrete gradient applied to the Hamiltonian form of the Duffing oscillator (see also \cite{franke2023} for a discussion in a more general setting). 
}

For the test of convergence the error are computed using an $L^2$ norm in time. The $L^2$ space norm is replaced by the Euclidian inner product for simplicity
$$
\begin{aligned}
    \mathrm{Error} \; \mathbf{q} &= \sqrt{\sum_{n=0}^{N_t} \Delta t ||\mathbf{q}_n - \mathbf{q}_{\mathrm{ref}}(t=n\Delta t)||^2}, \\
    \mathrm{Error} \; \mathbf{v} &= \sqrt{\sum_{n=0}^{N_t} \Delta t ||\mathbf{v}_n - \mathbf{v}_{\mathrm{ref}}(t=n\Delta t)||^2}.
\end{aligned}
$$
The reference solution is either an exact solution (in the case of the Duffing oscillator) or a solution computed using a leapfrog method with time step given by:
\begin{itemize}
    \item ${\Delta t}_{\rm ref} = \Delta t_{\rm base}/2^6$ for the von-K\`arm\`an beam;
    \item ${\Delta t}_{\rm ref} = \Delta t_{\rm base}/2^7$ for geometrically nonlinear elasticity.
\end{itemize}
In both cases $\Delta t_{\rm base}$ is the time step for the coarsest simulation. 

The finite element library \firedrake \cite{rathgeber2017firedrake} is used for the numerical investigation.
The \firedrake component Slate \cite{gibson2020slate} is used to implement the static condensation and the local solvers for stress reconstruction. 

\subsection{The Duffing oscillator}

\begin{figure*}[htbp]%
\centering
\subfloat[][Displacement $q$]{%
	\label{fig:err_q_duff}%
\includegraphics[width=0.45\textwidth]{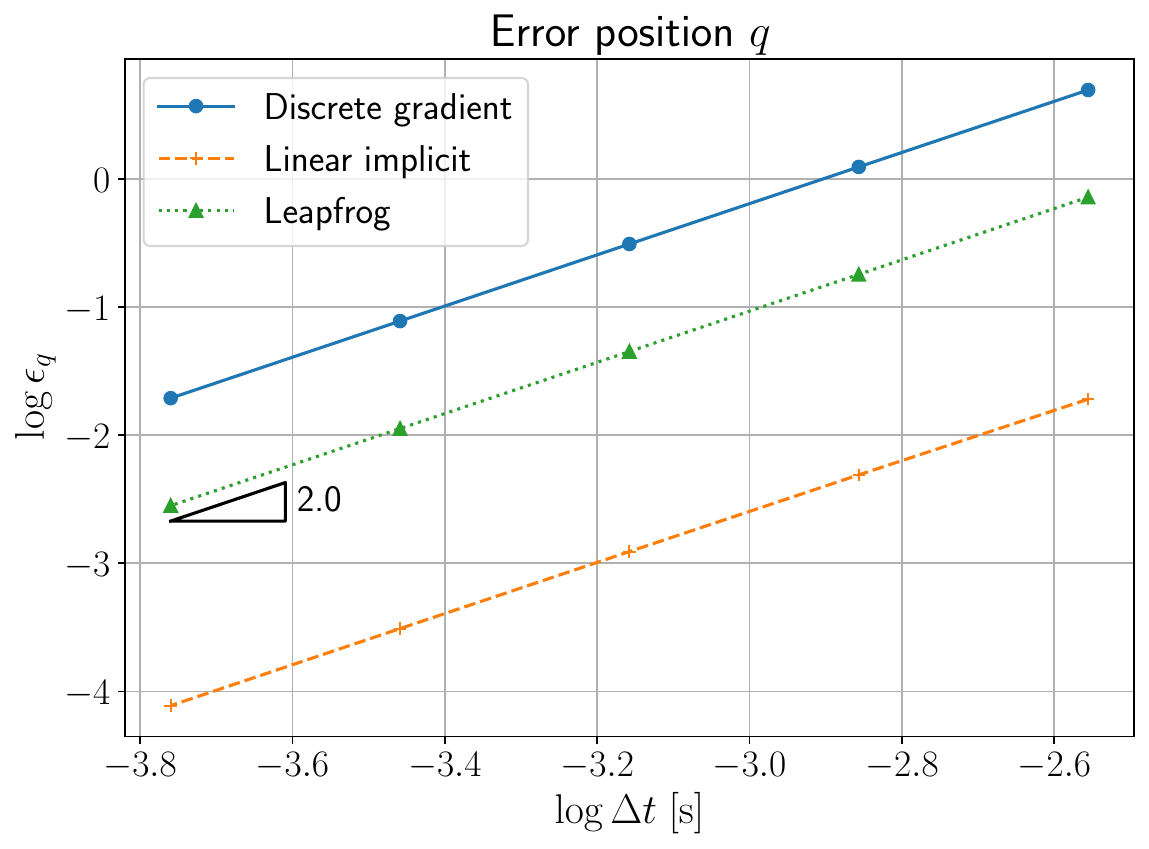}}
\hspace{8pt}
\subfloat[][Velocity $v$]{%
	\label{fig:err_v_duff}%
\includegraphics[width=0.45\textwidth]{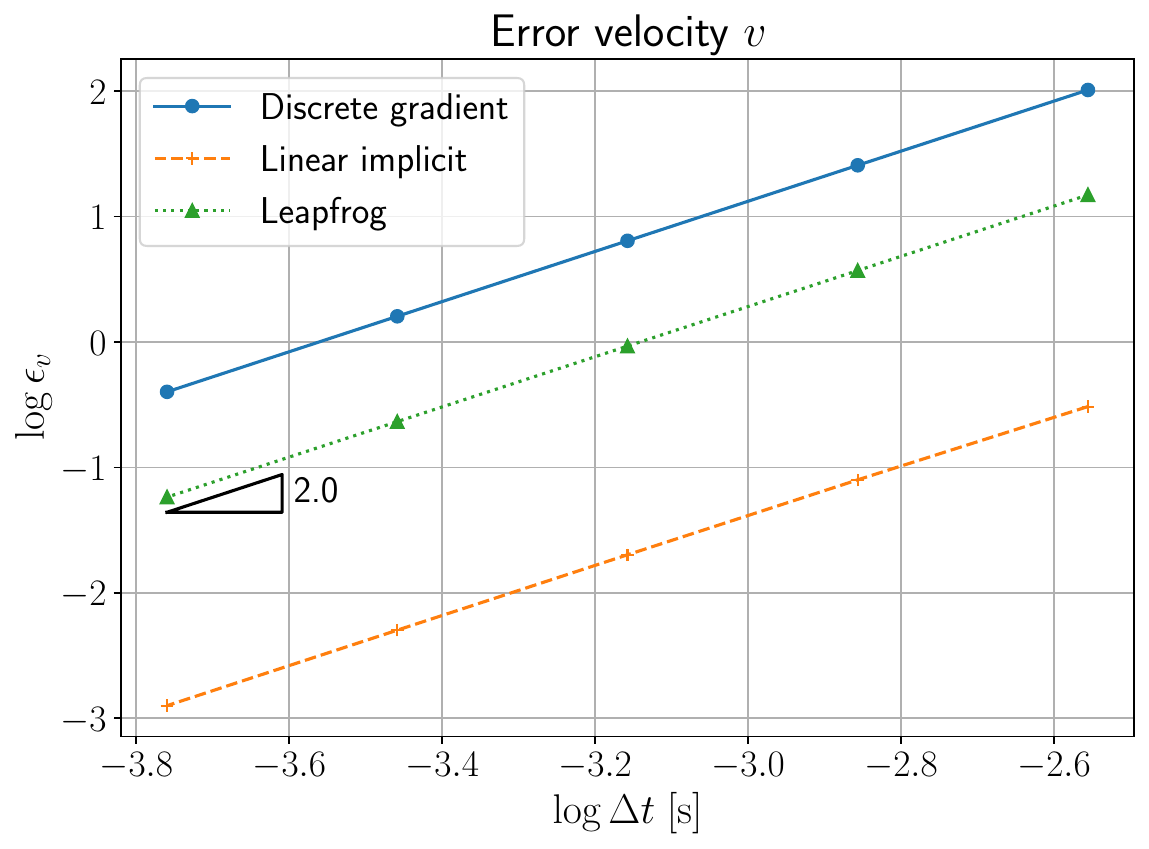}}
\caption{Convergence rate for the Duffing oscillator}%
\label{fig:conv_duff}%
\end{figure*}

\begin{table}[tb]
\centering
\begin{tabular}{|c|c|}
\hline
Parameter & Value\\
\hline
$\alpha$ & 10 \\
$\beta$ & 5 \\
$q_0$ & 10 \\
$T$ & $\displaystyle \frac{2\pi}{\sqrt{\alpha + \beta q_0^2}}$ \\
$T_{\rm end}$ & 100 $T$ \\
$\Delta t_{\rm base}$ & $0.278 \; \mathrm{[ms]}$\\
\hline
\end{tabular}
\caption{Parameters for the Duffing oscillator}
\label{tab:par_Duffing}
\end{table}
The presented scheme is compared with the exact solution of the unforced and undamped Duffing oscillator, the dynamics of which is ruled by the ODE
\begin{equation*}
\ddot{q} = -\alpha q - \beta q^3.
\end{equation*}
Notice the this ODE is equivalent to system \eqref{eq:Duffing_oscillator} with 
$$\alpha=2 k_{\rm hor}/m, \qquad \beta= k_{\rm ver}/(m L^2).$$ 
For the initial conditions we consider
$$q(t=0)=q_0, \qquad v(t=0)=0.$$
The analytic solution is given by the following expression
$$
\begin{aligned}
q(t) &= q_0 \mathrm{cn}\left(\omega_0 t; k^2\right), \\
v(t) &= -\omega_0 q_0\, \mathrm{sn}\left(\omega_0 t; k^2\right) \mathrm{dn}\left(\omega_0 t; k^2\right),  
\end{aligned}
$$
where 
$$\omega_0:=\sqrt{\alpha + \beta q_0^2}, \qquad k^2:= \frac{\beta q_0^2}{2(\alpha + \beta q_0^2)},$$
and $\mathrm{cn}(z; m), \; \mathrm{sn}(z; m)\; \mathrm{dn}(z; m)$ are the Jacobi elliptic functions of argument $z$ and parameter $m$. The parameters for the simulation are reported in table \ref{tab:par_Duffing}, where $T_{\rm end}$ is final simulation time and $\Delta t$ the time step for the simulation. The convergence plot against the exact solution is presented in Fig. \ref{fig:conv_duff}. It can be noticed that the all methods exhibit a second order convergence. However the linear implicit method is much more precise than both the discrete gradient and leapfrog method. \\

\begin{figure*}[htbp]%
\centering
\subfloat[][Displacement $q$]{%
	\label{fig:trend_q_duff}%
\includegraphics[width=0.45\textwidth]{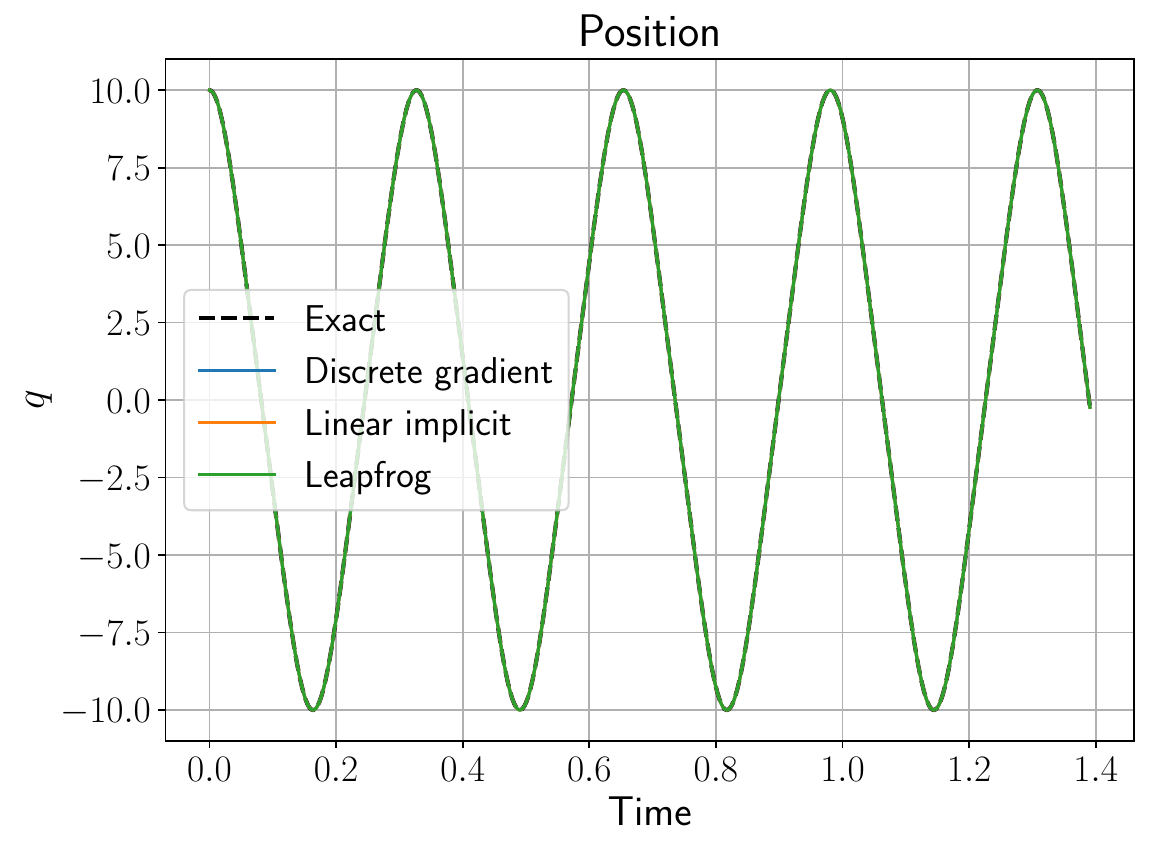}}
\hspace{8pt}
\subfloat[][Velocity $v$]{%
	\label{fig:trend_v_duff}%
\includegraphics[width=0.45\textwidth]{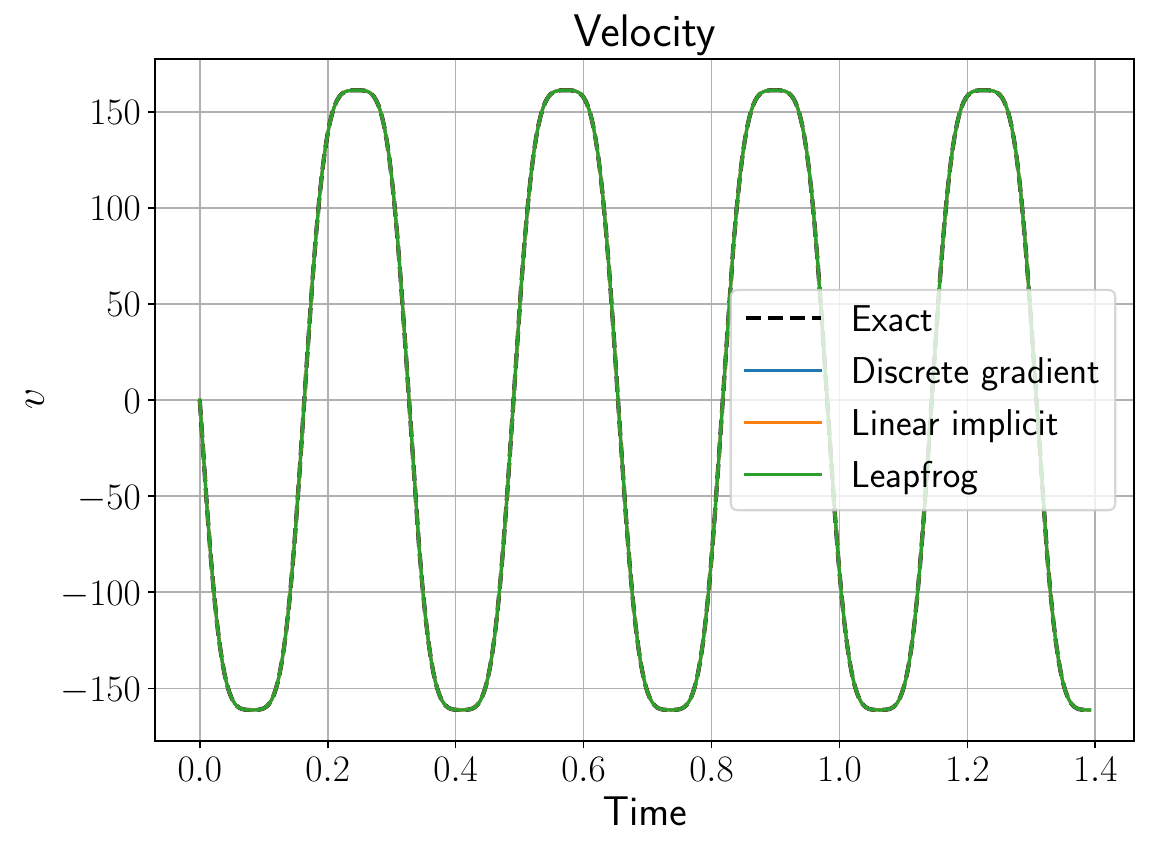}}
\caption{Time trend of the position and velocity for the Duffing oscillator}%
\label{fig:trend_variables_Duffing}%
\end{figure*}

\begin{figure*}[htbp]
    \centering
    \begin{minipage}{0.45\textwidth}
    \includegraphics[width=0.98\textwidth]{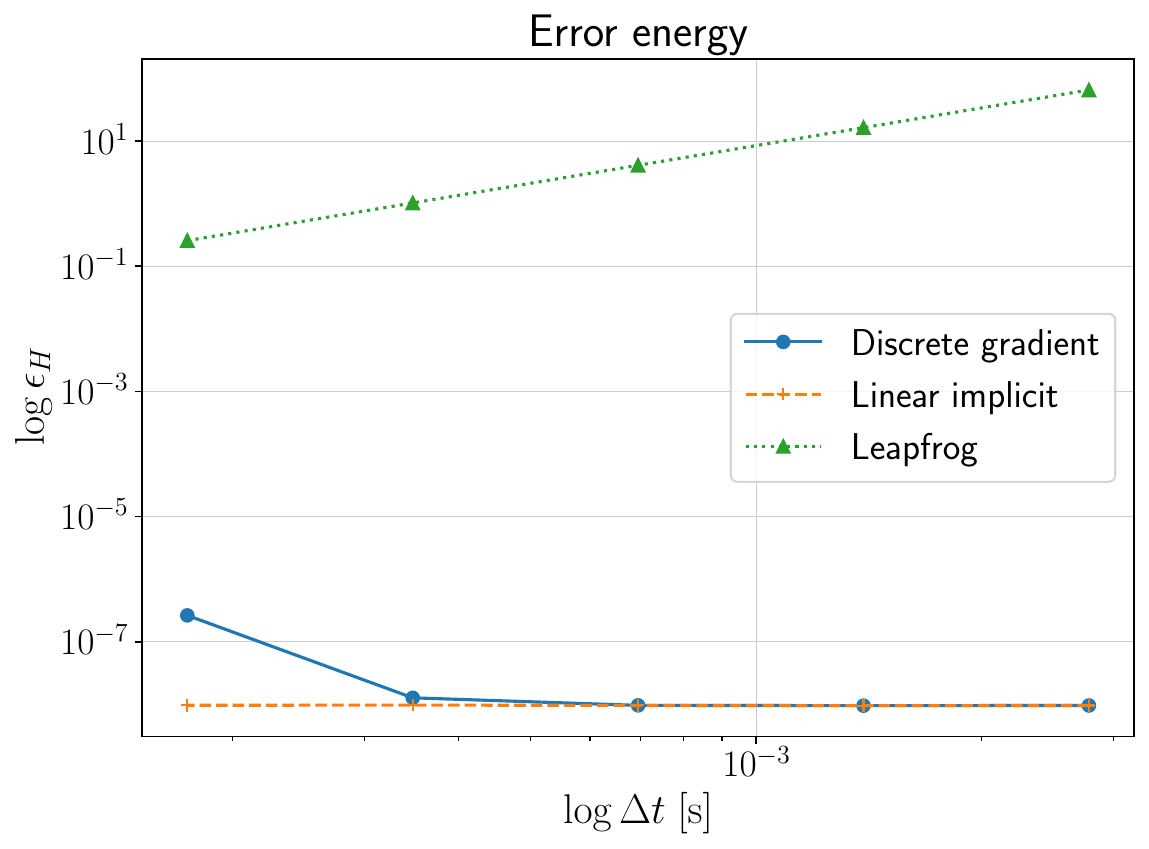}
    \caption{Energy error for the Duffing oscillator}
    \label{fig:err_energy_Duffing}
    \end{minipage}\hspace{8pt}
    \begin{minipage}{0.45\textwidth}
    \includegraphics[width=0.98\textwidth]{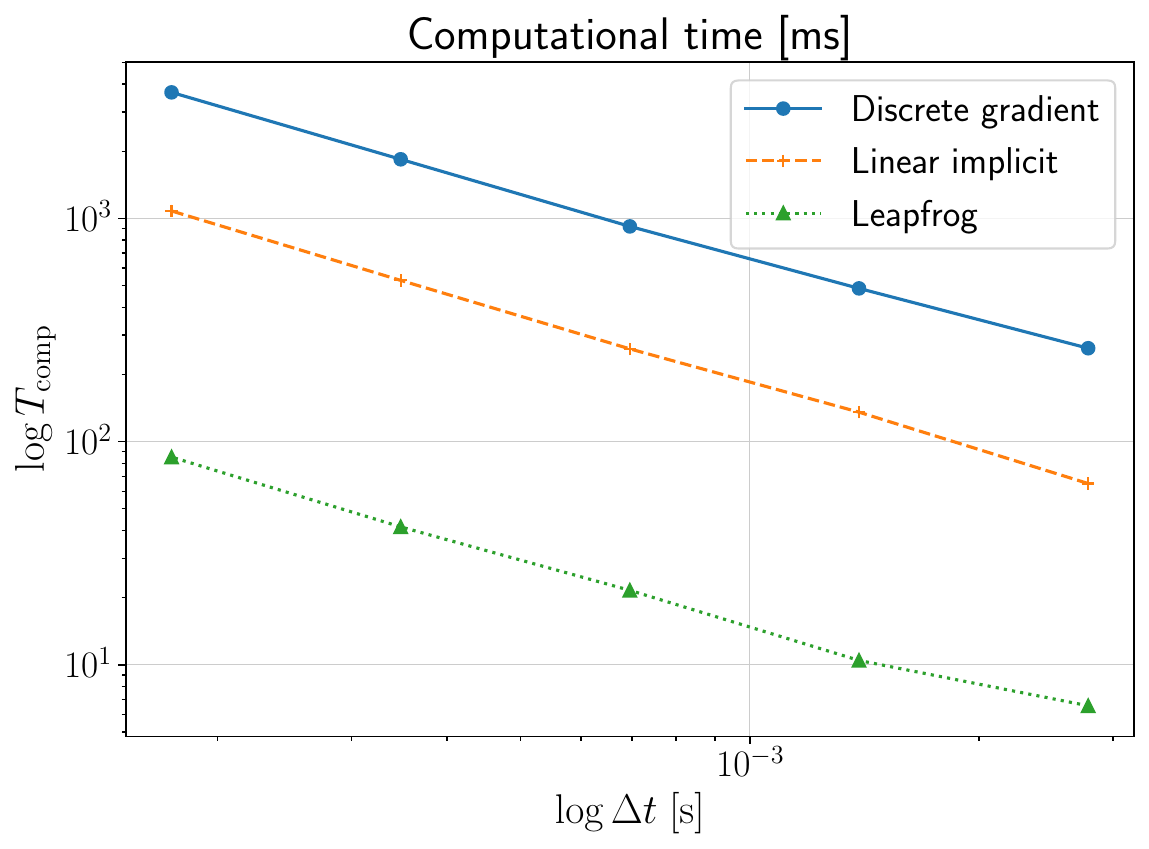}
    \caption{Computational time for Duffing}
    \label{fig:comp_time_Duffing}
    \end{minipage}
\end{figure*}

\rewone{Figures \ref{fig:trend_variables_Duffing} shows the time signals of position of velocity for the chosen parameters. The nonlinear oscillations deviating from a sinusoidal trend are clearly visible.}
In Fig. \ref{fig:err_energy_Duffing} the error with respect to the exact energy is shown. One can notice that the leapfrog method only approximately conserves the energy where the discrete gradient and linear implicit method preserve it up to machine precision. Concerning the total computational time, reported in Fig. \ref{fig:comp_time_Duffing},  the linear implicit method stays in between the leapfrog and discrete gradient methods. In this example, however, there is no mass matrix coming from a finite element discretization, so it is natural that a fully explicit method is by one order of magnitude faster.

\subsection{Vibrations of a von-K\'arm\'an beam}

\begin{table}[htbp]
\centering
\begin{tabular}{|c|c|}
\hline
Parameter & Value \\
\hline
$\rho$ & $2700 \; \mathrm{[kg/m^3]}$ \\
$E$ & $70 \; \mathrm{[GPa]}$ \\
$L$ & $1 \; \mathrm{[m]}$ \\
$d$ & $2 \; \mathrm{[mm]}$ \\
$A_z$ & $d$ \\
$N_{\rm el}$ & 50 \\
$T_{1, \mathrm{bending}}$ & $\displaystyle 2\frac{L^2}{\pi} \sqrt{\frac{\rho d^2}{EI}}$ \\
$\Delta t_{\rm base}$ & $17 \; [\mu \mathrm{s}]$ \\
$T_{\rm end}$ & $5\, T_{1, \mathrm{bending}}$ \\
\hline
\end{tabular}
\caption{Parameters for the von-K\'arm\'an beam}
\label{tab:par_vonkarman}
\end{table}

In this section the free vibration of a geometrically nonlinear beam are analyzed. The actual values for the literal that will appear hereafter are reported in Table \ref{tab:par_vonkarman}. \rewtwo{A prismatic beam of length $L$ and square cross section whose side has length $d$ is considered}. As boundary conditions, \rewone{the beam is taken to be simply supported}
\begin{equation*}
\begin{aligned}
    q_x(x=0, \; t) = q_x(x=L, \; t) = 0, \\
    q_z(x=0, \; t) = q_z(x=L, \; t) = 0.
\end{aligned}
\end{equation*}
For the initial conditions, we consider initial displacement given by the first mode of vibration for both the vertical displacement, and zero velocities
\begin{equation*}
    \begin{aligned}
        q_x(x, \; 0) &= 0, \\
        v_x(x, \; 0) &= 0.
    \end{aligned} \qquad 
    \begin{aligned}
        q_z(x, \; 0) &= A_z \sin\left(\pi x/L\right), \\
        v_z(x, \; 0) &= 0.
    \end{aligned}
\end{equation*}
The beam is discretized using $N_{\rm el} = 50$ finite elements. The total simulation time is taking to be five times the period of the first bending mode, i.e. \rewtwo{$T_{\rm end} = 5\, T_{1, \rm bending}$}. The time step for the coarsest simulation is taken to be $\Delta t_{\rm base} = \frac{1}{4} \Delta t_{\rm CFL, bend} = 17 \; [\mu \mathrm{s}]$ where 
$$\Delta t_{\rm CFL, bend} = \frac{1}{2}\left(\frac{L}{N_{\rm el}}\right)^2 \sqrt{\frac{EI}{\rho h^2}}$$ 
is the minimum time step for stability in finite difference simulations (in space and time) for bending of beams, cf. \cite[Chapter 5]{bilbao2001phdthesis}. \rewone{Because of mass matrix arising in the discretization, in finite elements the value is reduced by a factor to obtain a stable simulation.} \\

\begin{figure*}[htbp]
\centering
\subfloat[][Horizontal displacement $q_x$]{%
	\label{fig:ref_q_x}%
\includegraphics[width=0.45\textwidth]{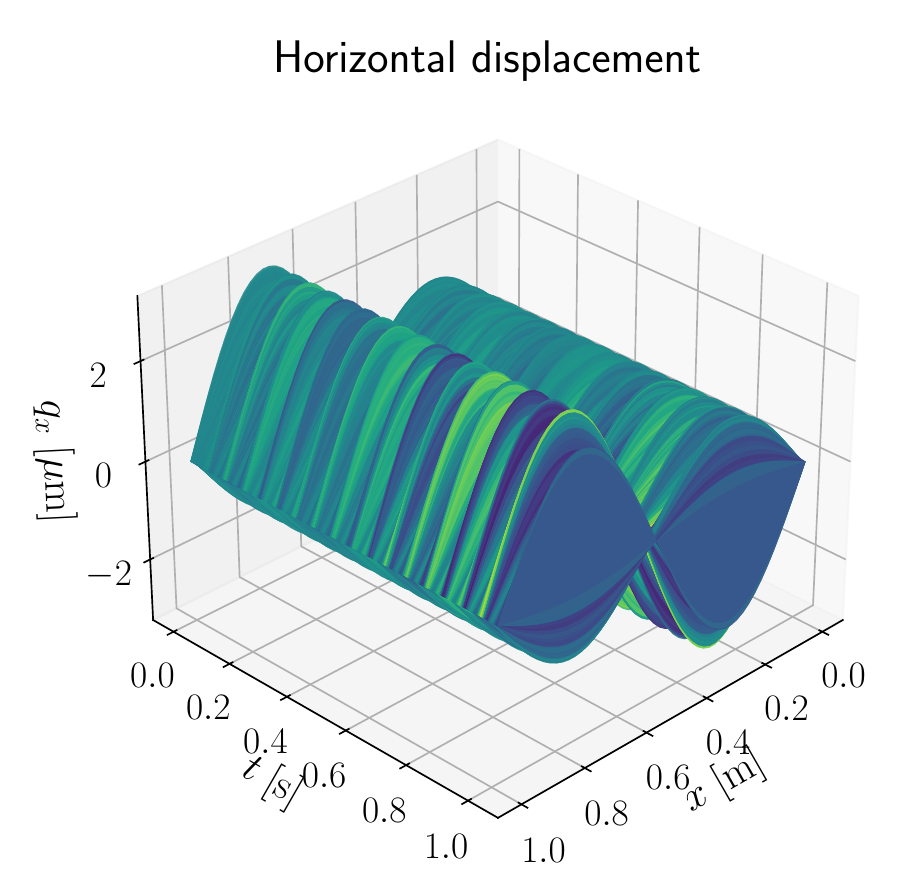}}%
\hspace{8pt}
\subfloat[][Vertical displacement $q_z$]{%
	\label{fig:ref_q_z}%
\includegraphics[width=0.45\textwidth]{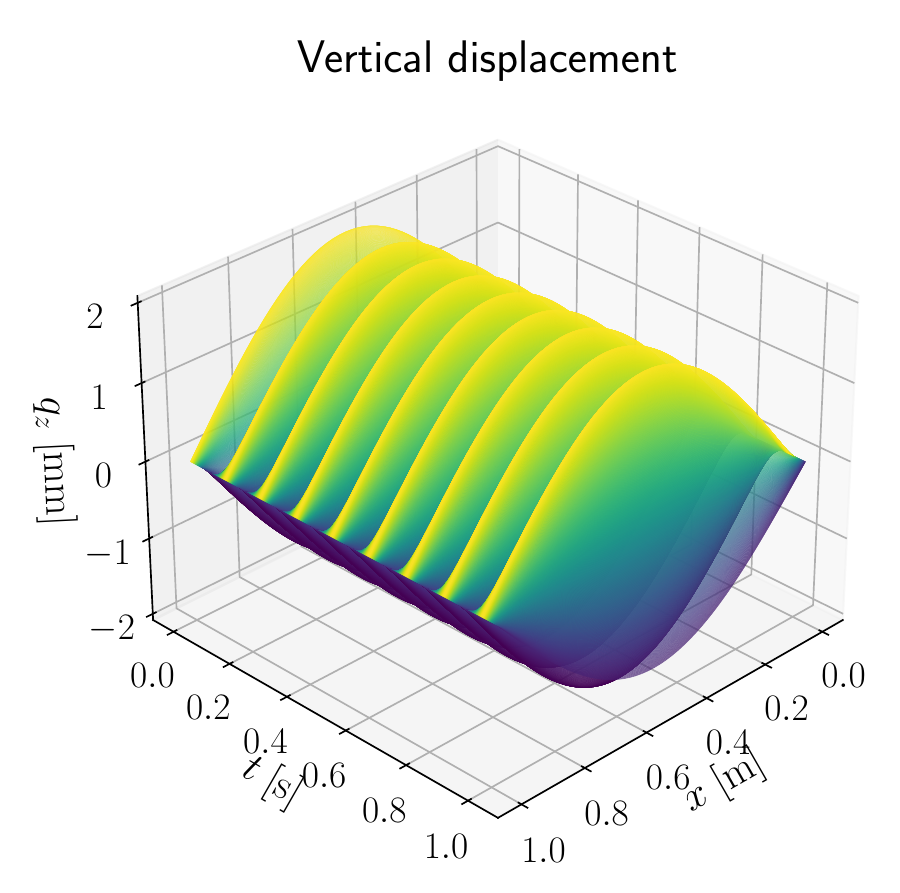}}
\caption{Reference solution for the von-K\'arm\'an beam obtained using Leapfrog with $\Delta t_{\rm base}/2^6$}%
\label{fig:ref_vonkarman}%
\end{figure*}

\begin{figure*}[htbp]
\centering
\subfloat[][Linear Horizontal displacement $q_x$]{%
	\label{fig:linear_q_x}%
\includegraphics[width=0.45\textwidth]{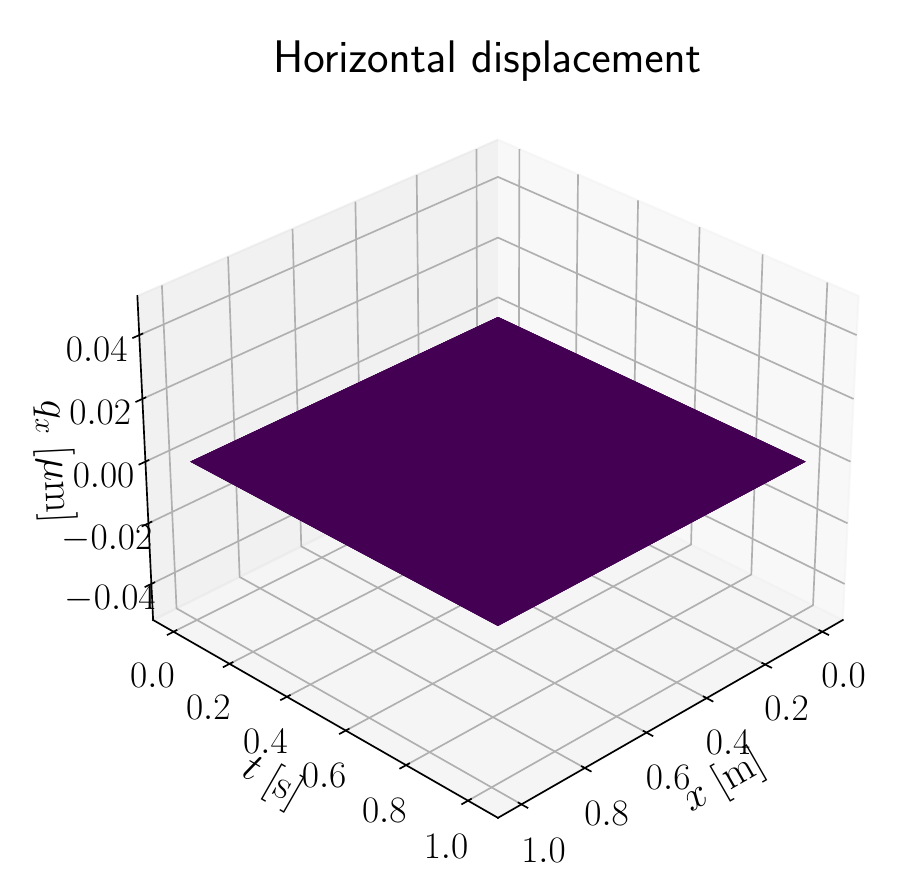}}%
\hspace{8pt}
\subfloat[][Linear Vertical displacement $q_z$]{%
	\label{fig:linear_q_z}%
\includegraphics[width=0.45\textwidth]{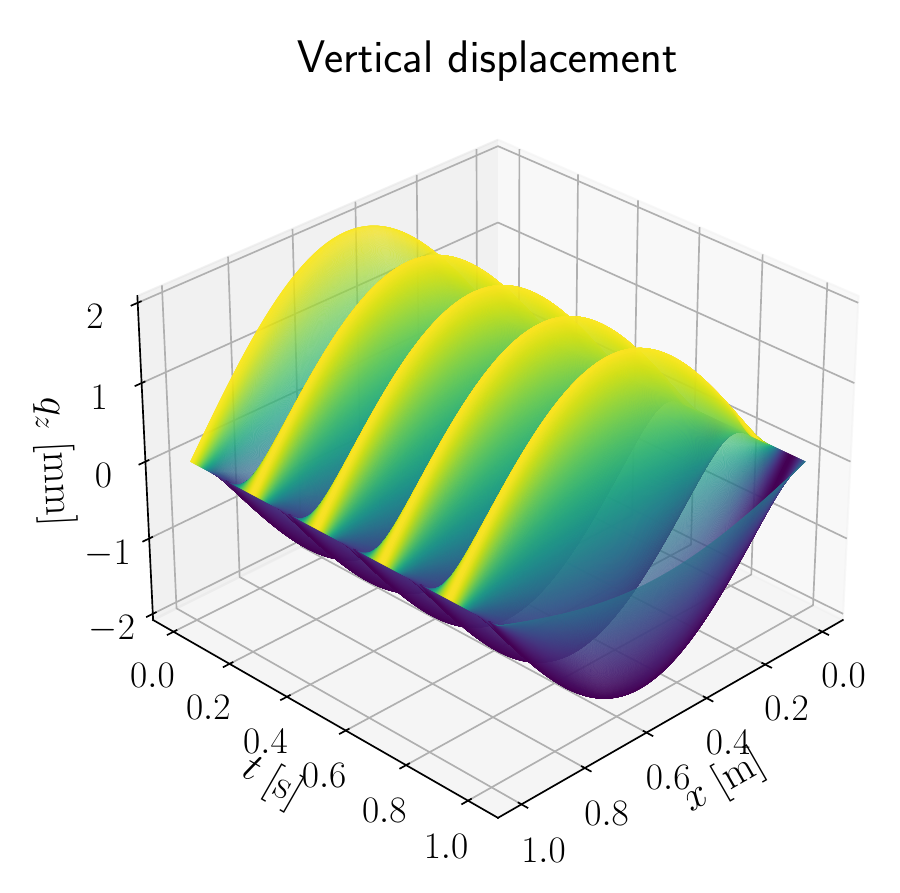}} 
\caption{Linear solution for the von-K\'arm\'an beam}%
\label{fig:linear_vonkarman}%
\end{figure*}

\begin{figure*}[htbp]
\centering
\subfloat[][Horizontal displacement $q_x(L/4, t)$]{%
	\label{fig:ref_q_x_at_point}%
\includegraphics[width=0.45\textwidth]{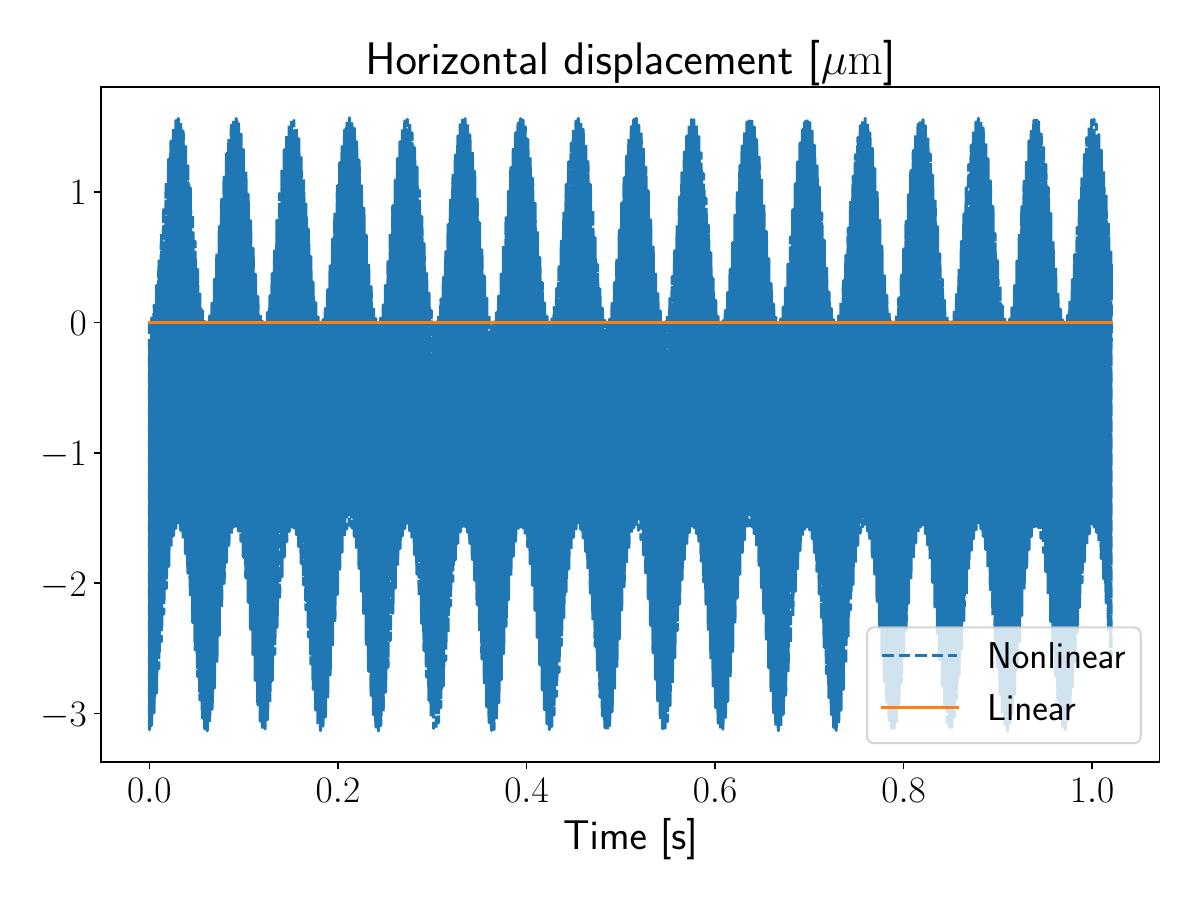}}%
\hspace{8pt}
\subfloat[][Linear Vertical displacement $q_z(L/4, t)$]{%
	\label{fig:ref_q_z_at_point}%
\includegraphics[width=0.45\textwidth]{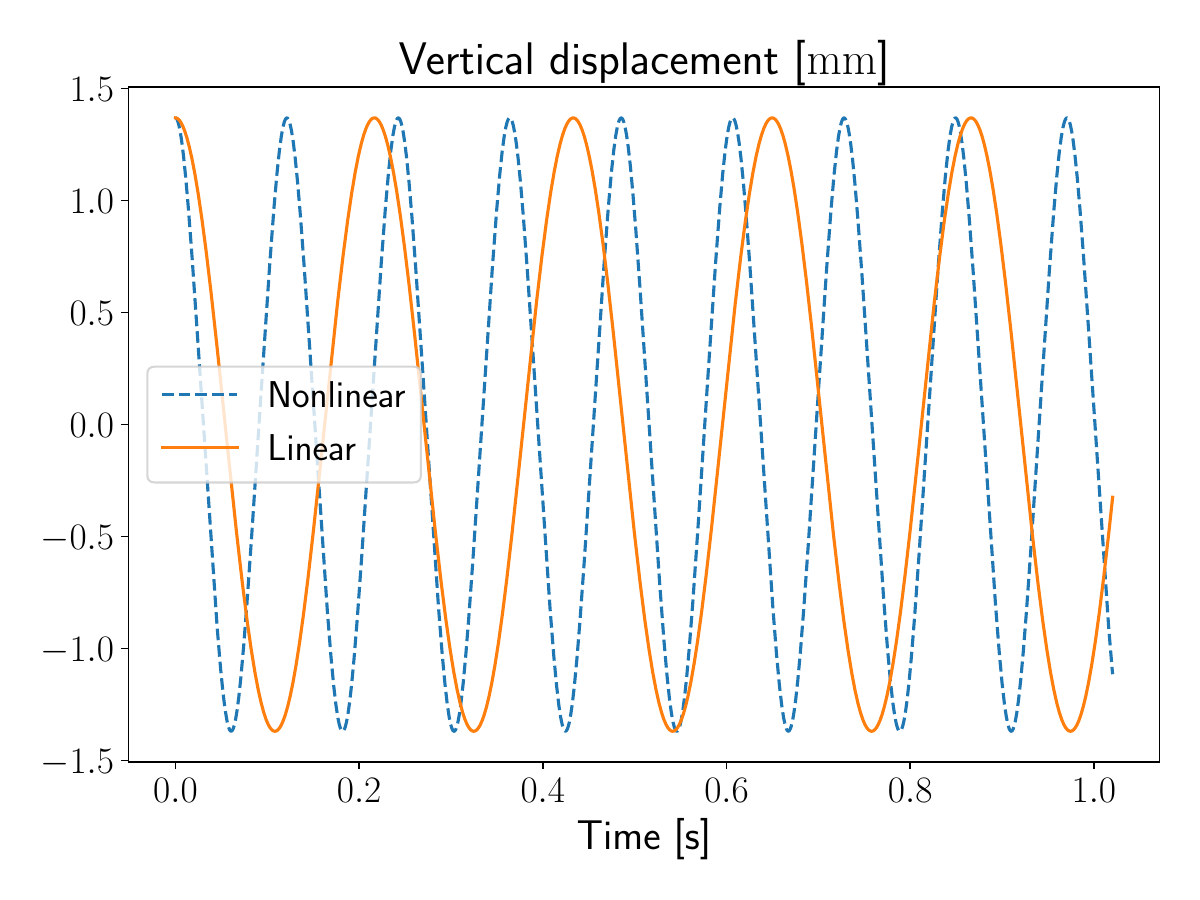}}
\caption{Comparison between linear and nonlinear case at $x=L/4$ for the von-K\'arm\'an beam}%
\label{fig:vonkarman_lin_nonlin}%
\end{figure*}

\begin{figure*}[htbp]
\centering
\subfloat[][Horizontal displacement $q_x$]{%
	\label{fig:err_q_x_vk}%
\includegraphics[width=0.45\textwidth]{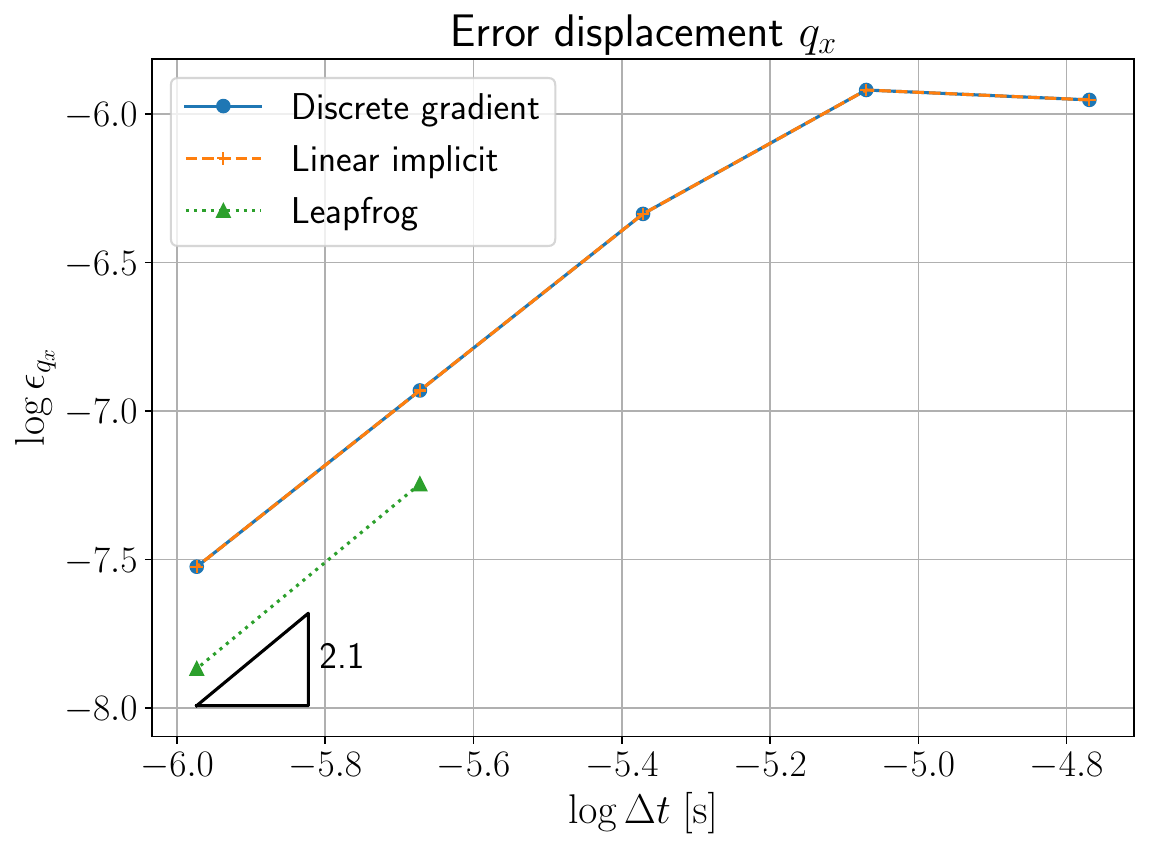}}%
\hspace{8pt}
\subfloat[][Vertical displacement $q_z$]{%
	\label{fig:err_q_z_vk}%
\includegraphics[width=0.45\textwidth]{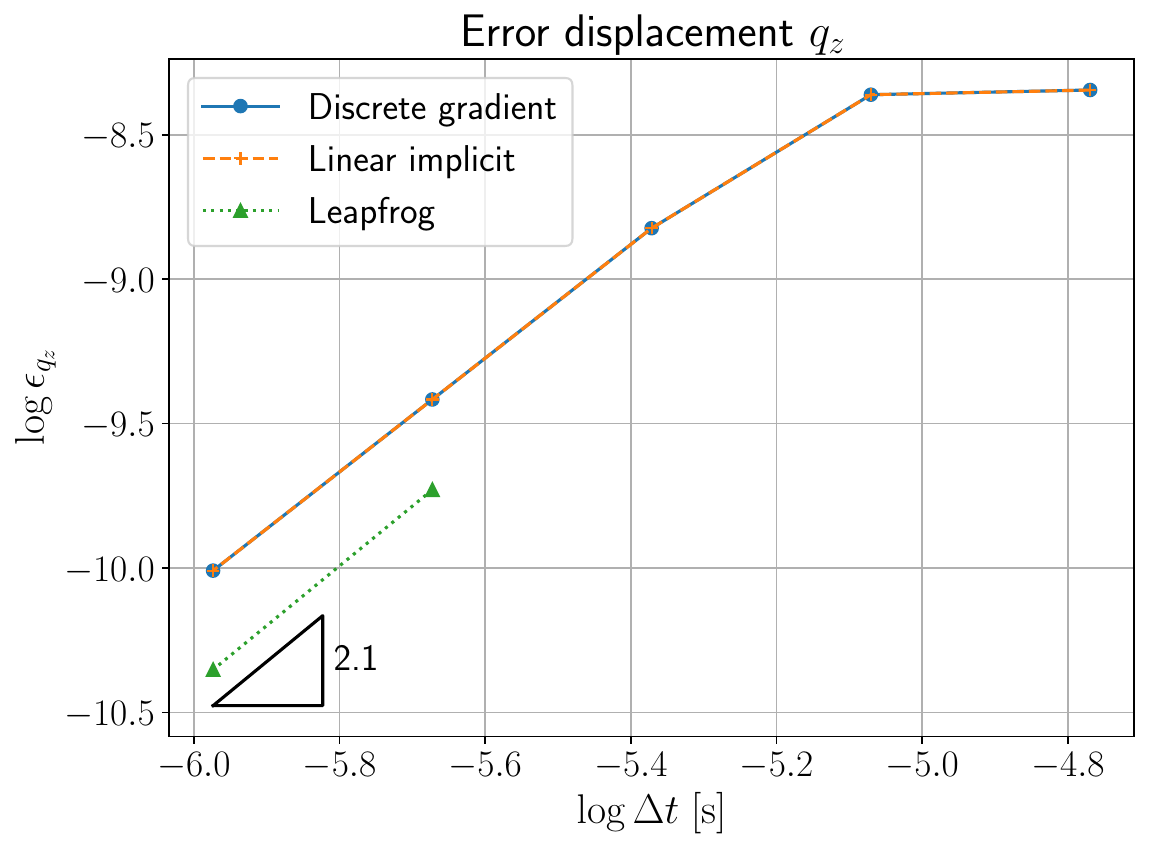}} \\
\subfloat[][Horizontal velocity $v_x$]{%
	\label{fig:err_v_x_vk}%
\includegraphics[width=0.48\textwidth]{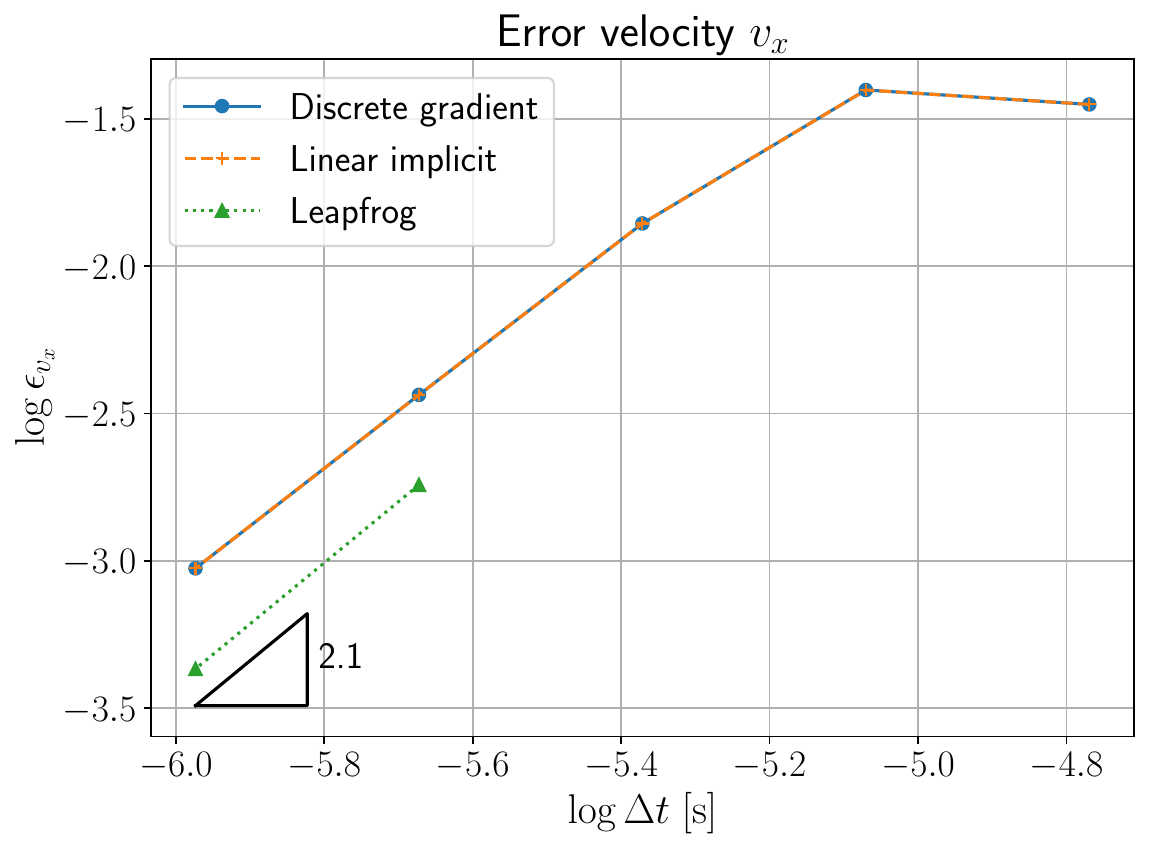}}%
\hspace{8pt}
\subfloat[][Vertical velocity $v_z$]{%
	\label{fig:err_v_z_vk}%
\includegraphics[width=0.45\textwidth]{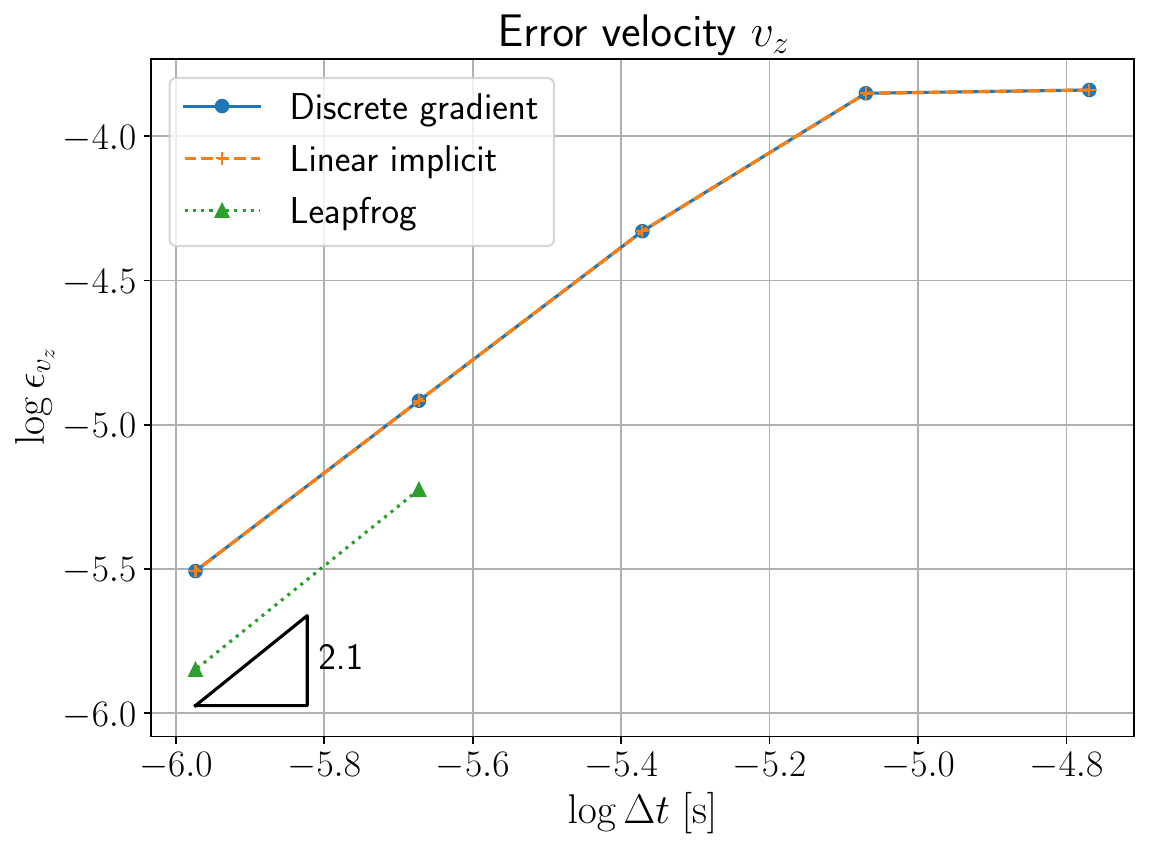}}
\caption{Convergence rate for $q_x, q_z, v_x, v_z$ in the von-K\'arm\'an beam considering the time span $[0, T_{1, \rm bend}/10]$}%
\label{fig:conv_vonkarman}%
\end{figure*}

\begin{figure*}[htbp]
    \centering
    \begin{minipage}{0.45\textwidth}
    \includegraphics[width=0.98\textwidth]{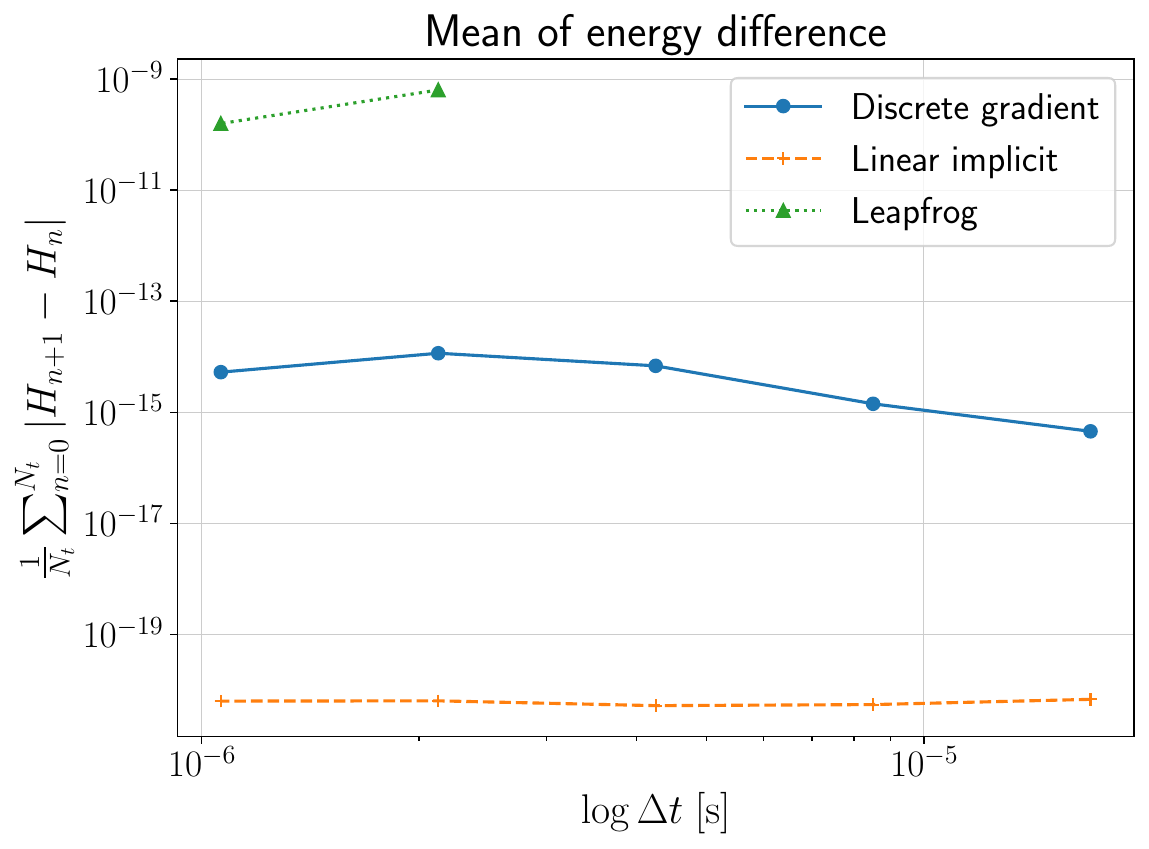}
    \caption{Mean of $|H_{n+1}-H_n|$ for the von-K\'arm\'an beam}
    \label{fig:err_energy_vonkarman}
    \end{minipage}\hspace{8pt}
    \begin{minipage}{0.45\textwidth}
    \includegraphics[width=0.98\textwidth]{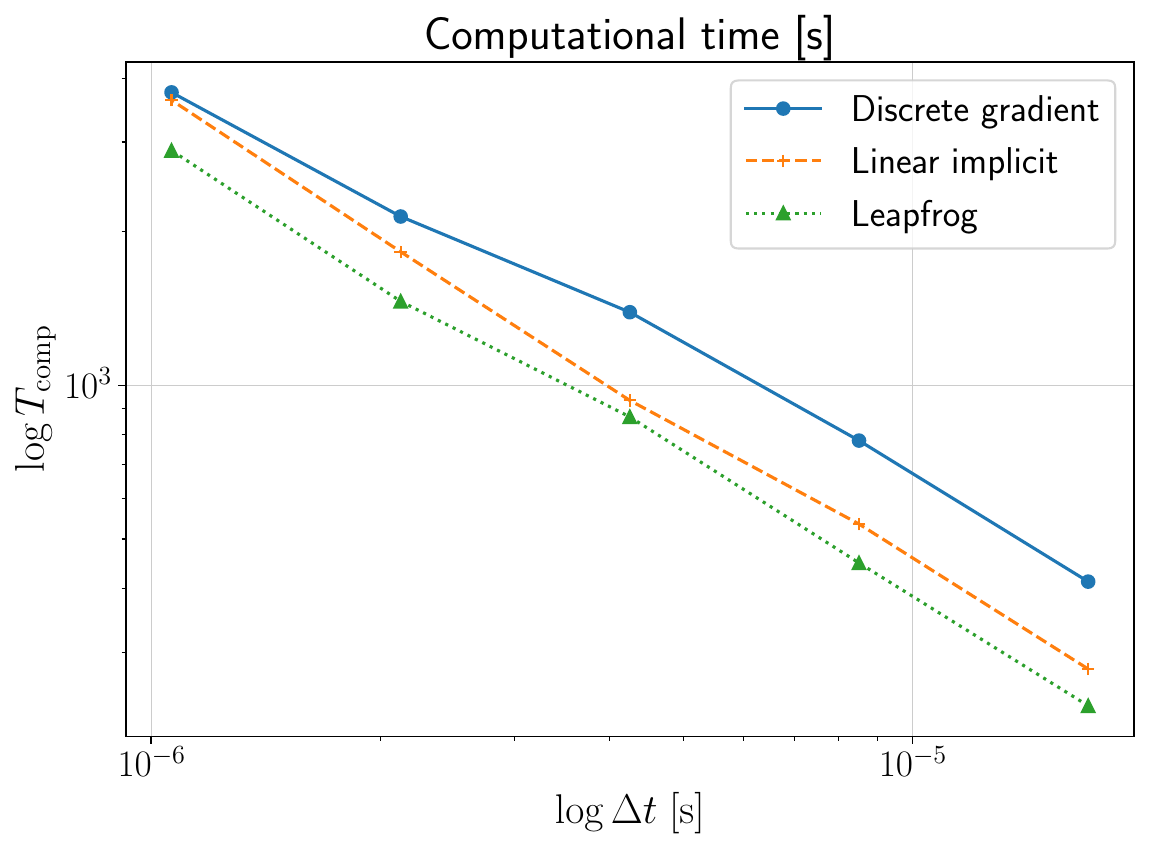}
    \caption{Computational time for the von-K\'arm\'an beam}
    \label{fig:comp_time_vonkarman}
    \end{minipage}
\end{figure*}

Given the nonlinear coupling, the axial and bending dynamics interact and the reference solution, obtained using the leapfrog method with a time step $\Delta t_{\rm ref} = \Delta t_{\rm base}/2^6$, is reported in Fig. \ref{fig:ref_vonkarman}. \rewtwo{The linear solution is shown for comparison in Fig. \ref{fig:linear_vonkarman}. The different trajectories for the nonlinear and linear solutions at point $x=L/4$ can be observed in Fig. \ref{fig:vonkarman_lin_nonlin}. It can be noticed that the vertical displacement follows the mode shape with a period different then the one expected from the linear analysis.} The convergence  of the different variables with respect to the reference solution is reported in Fig.~\ref{fig:conv_vonkarman}. It can be noticed that all methods converge with a second order convergence and the discrete gradient and linear implicit integrators have identical precision. The leapfrog method is unstable for the time step $\Delta t_{\rm base}$ as the axial dynamics is way faster that the bending one and thus requires a time step of the order of $\Delta_{\rm CFL, ax} = \frac{L}{N_{\rm el}} \sqrt{\frac{E}{\rho}}$. Since the reference solution is computed with the leapfrog method, it is not surprising that this method achieves a better performance. For what concerns the energy conservation in Fig. \ref{fig:err_energy_vonkarman} the mean of the energy different between two time steps is reported. The proposed linearly implicit method respects energy conservation to machine precision. The discrete gradient method respects energy conservation but not as accurately, probably because of tolerances settings in the Newton method to compute the nonlinear solution. The computational time taken by the three methods is reported in Fig. \ref{fig:comp_time_vonkarman}. No static condensation is applied as the system to be solved is rather small and no noticeable improvements are obtained. Since now mass matrices are included in the formulation the different methods take comparable time for the same time step. Of course the considering mass lumping strategies would be of substantial help for the leapfrog method, but also for the linearly implicit method, since a conjugate gradient method would be much more effective when preconditioning with the (lumped) mass matrix \cite{wathen2015prec}. The incorporation of mass lumping strategies is not considered in the present contribution. It can be noticed 

\subsection{Geometrically nonlinear elasticity}

\begin{figure*}[htbp]
\centering
\begin{minipage}[b]{0.45\textwidth}
\centering
\begin{tabular}{|c|c|}
\hline
Parameter & Value \\
\hline
$\rho$ & $1100 \; \mathrm{[kg/m^3]}$ \\
$E$ & $17 \; \mathrm{[MPa]}$ \\
$\nu$ & 0.3 \\
$\mu$ & $\displaystyle \frac{E}{2(1+\nu)}$ \\
$\lambda$ &  $\displaystyle \frac{E \nu}{(1-2\nu)(1+\nu)}$ \\
$\kappa$ &  $\displaystyle  \lambda + \frac{2}{3}\mu$ \\
$c_l$ &  $\displaystyle  \sqrt{\frac{\kappa + 4/3\mu}{\rho}}$ \\
$\Delta t_{\rm base}$ & $1.16 \; \mathrm{[ms]}$ \\
$T_{\rm end}$ & $0.5 \; \mathrm{[s]}$ \\
\hline
\end{tabular}
\caption{Parameters for the geometrically nonlinear elasticity problem}
\label{tab:par_elasticity}
\end{minipage}
\hspace{8pt}
\begin{minipage}[b]{0.45\textwidth}
\centering
\includegraphics[width=.75\textwidth]{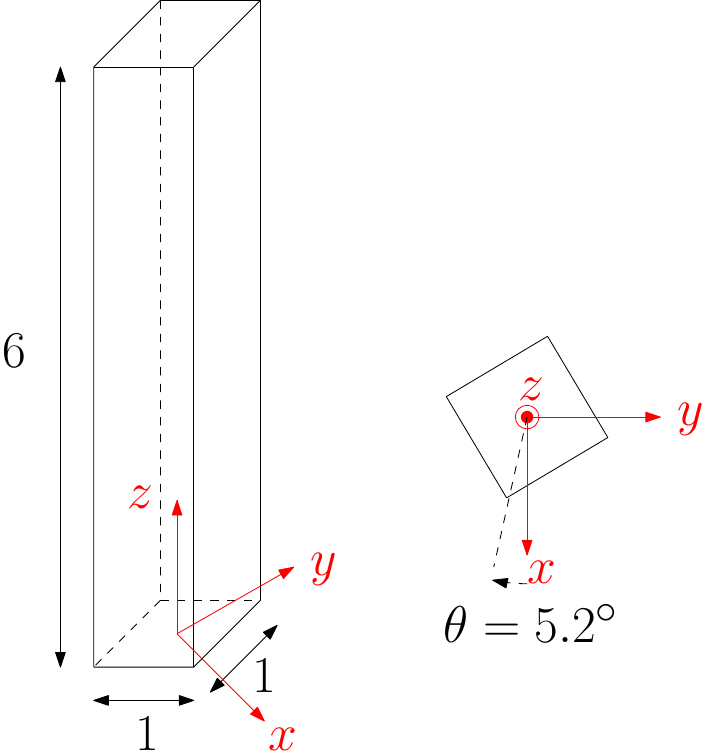}
\caption{Geometry for the geometrically nonlinear elasticity problem}
\label{fig:geo_column}
\end{minipage}
\end{figure*}

Non symmetrical oscillations of a column are considered in this test. The example is taken from \cite{scovazzi2016} and the only difference with respect to the simulation therein is the fact that a Saint-Venant constitutive model is considered here instead of a neo-Hookean one. The geometry and reference frame are reported in Fig. \ref{fig:geo_column} and distances are reported in meters. The column is considered clamped at its base and the initial conditions are taken to be
\begin{equation*}
    \bm{q}(\bm{x}, 0) = \bm{0}, \qquad \bm{v}(\bm{x}, 0) = \begin{pmatrix}
        \frac{5}{3} z & 0 & 0
    \end{pmatrix}^\top \; \mathrm{[m/s]}.
\end{equation*}
The mesh is divided into $6 \times 6 \times 36$ cubes each composed of 6 tetrahedral, leading to a mesh size of $h = 1/6 \; \mathrm{[m]}$. The total simulation time is taken to be $T_{\rm end} = 0.5 \; \mathrm{[s]}$ and the time step for the coarsest simulation is $\Delta t_{\rm base} = h/c_l = 1.16 \; \mathrm{[ms]}$, where $c_l$ is the speed of propagation of longitudinal waves in the solid (ignoring nonlinear effects), reported in Table \ref{tab:par_elasticity}, together with the actual value of the physical parameters. \\

\begin{figure*}[htbp]
\centering
\subfloat[][$t=\frac{1}{4} T_{\rm end}$]{%
	\label{fig:disp_elasticity_1_4}%
\includegraphics[width=0.4\textwidth]{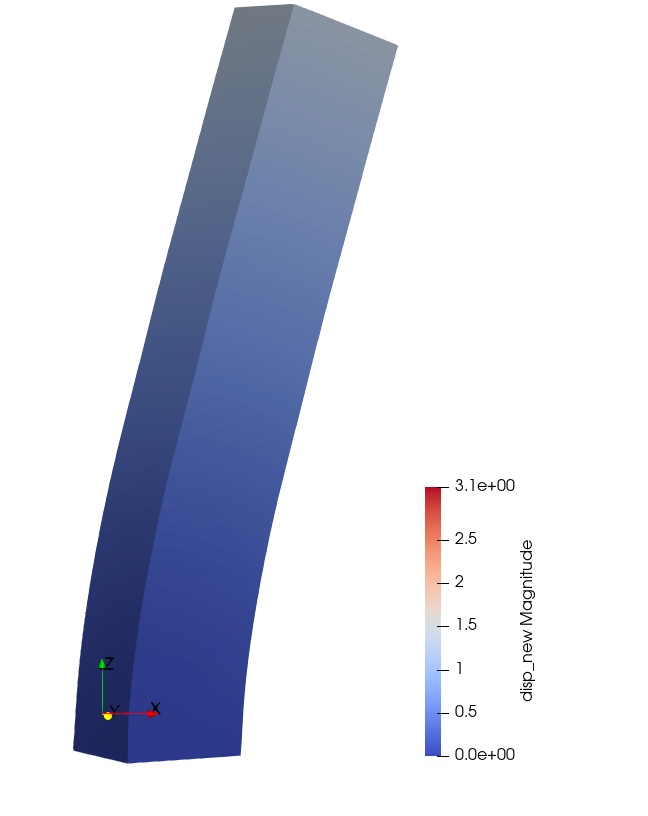}}%
 \hspace{8pt}%
\subfloat[][$t=\frac{1}{2} T_{\rm end}$]{%
	\label{fig:disp_elasticity_2_4}%
\includegraphics[width=0.4\textwidth]{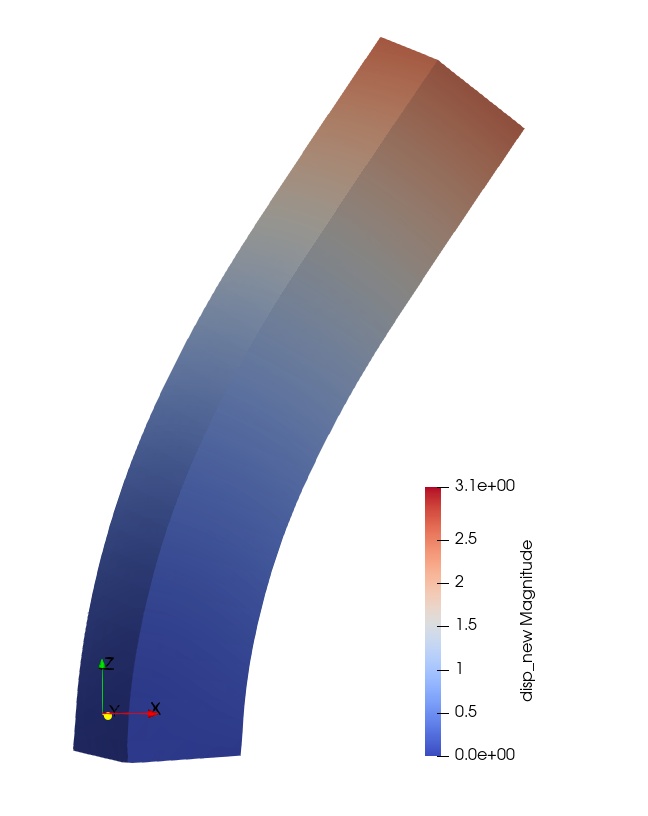}} \\
 \subfloat[][$t=\frac{3}{4} T_{\rm end}$]{%
	\label{fig:disp_elasticity_3_4}%
\includegraphics[width=0.4\textwidth]{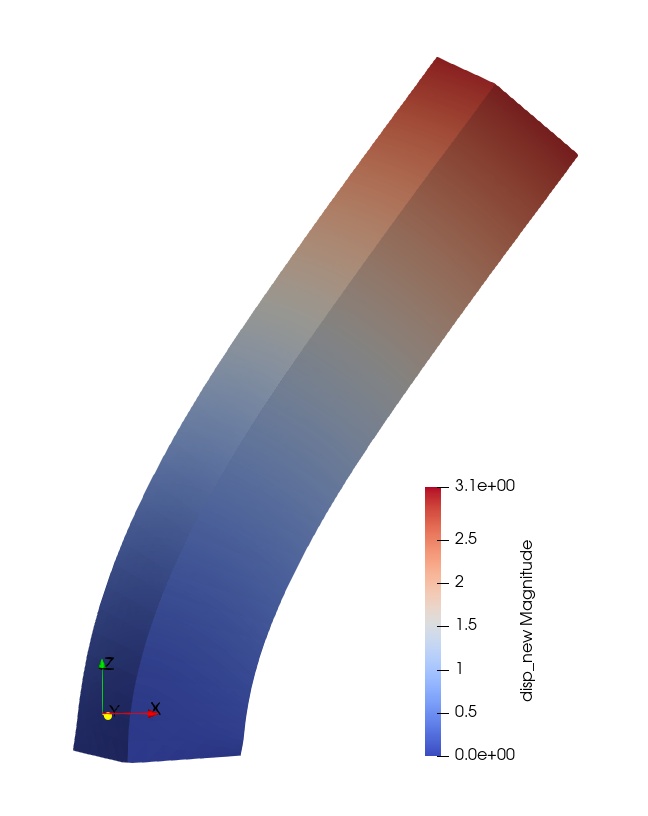}} 
\hspace{8pt}%
\subfloat[][$t=T_{\rm end}$]{%
	\label{fig:disp_elasticity_4_4}%
\includegraphics[width=0.4\textwidth]{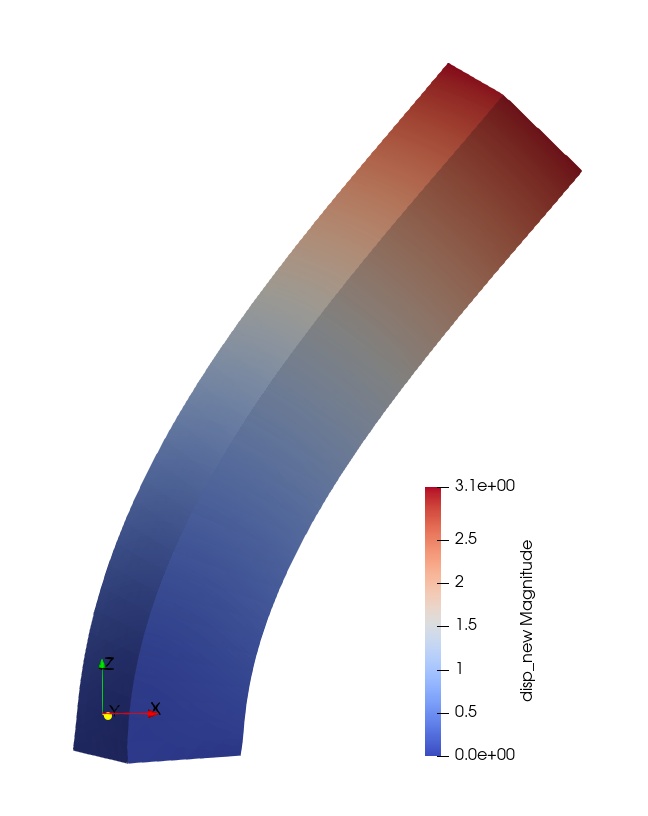}}
\caption{Screenshots of the displacement solution for geometrically nonlinear elasticity}%
\label{fig:screenshots_elasticity}%
\end{figure*}

\begin{figure*}[htbp]
\centering
\subfloat[][$t=\frac{1}{4} T_{\rm end}$]{%
	\label{fig:stress_elasticity_1_4}%
\includegraphics[width=0.4\textwidth]{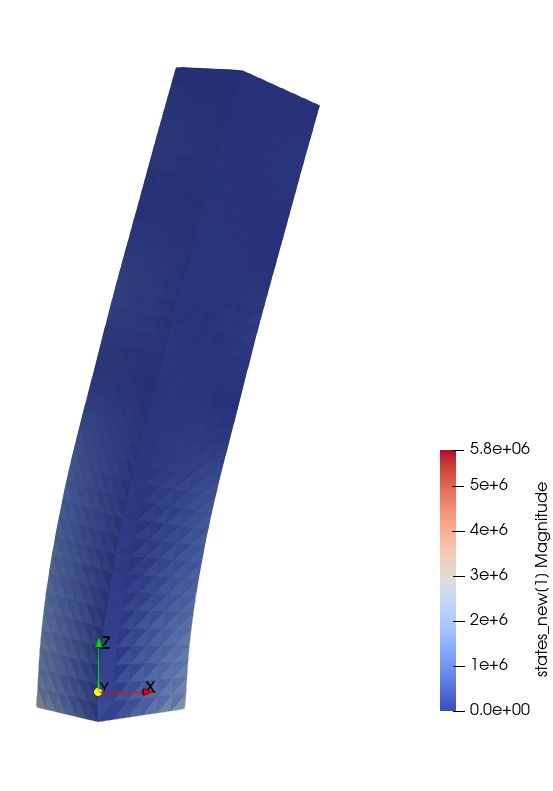}}%
 \hspace{8pt}%
\subfloat[][$t=\frac{1}{2} T_{\rm end}$]{%
	\label{fig:stress_elasticity_2_4}%
\includegraphics[width=0.4\textwidth]{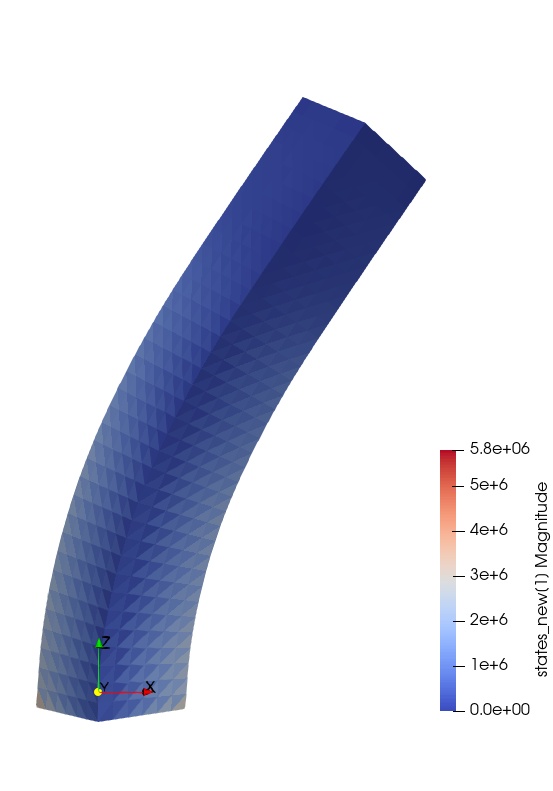}} \\
 \subfloat[][$t=\frac{3}{4} T_{\rm end}$]{%
	\label{fig:stress_elasticity_3_4}%
\includegraphics[width=0.4\textwidth]{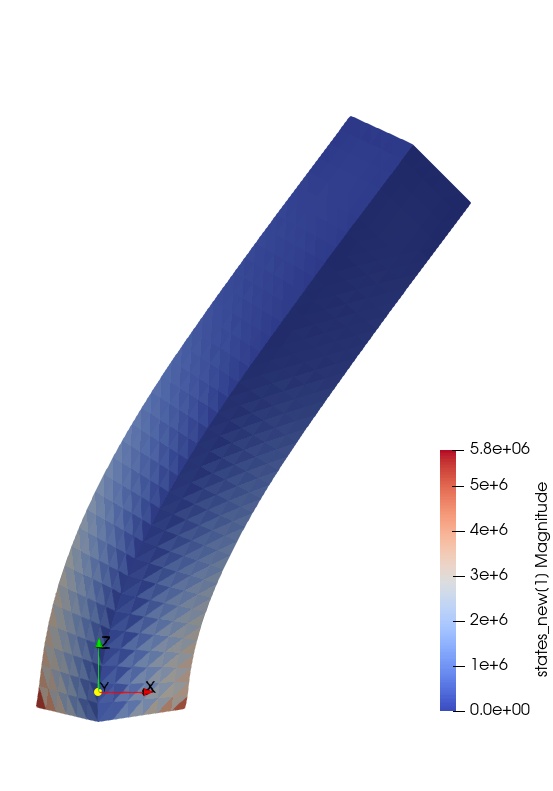}} 
\hspace{8pt}%
\subfloat[][$t=T_{\rm end}$]{%
	\label{fig:stress_elasticity_4_4}%
\includegraphics[width=0.4\textwidth]{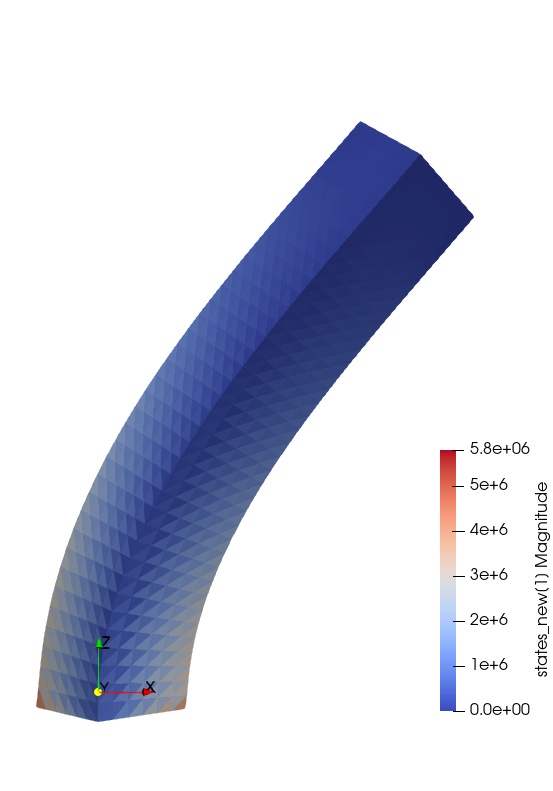}}
\caption{Screenshots of the Frobenius norm of the stress, i.e. $||\bm{S}||_F$ for geometrically nonlinear elasticity using the linear implicit scheme with $\Delta t = \Delta t_{\rm base}/2^5$ where
$\Delta t_{\rm base}= 1.16 \; \mathrm{[ms]}$. }%
\label{fig:screenshots_stress_elasticity}%
\end{figure*}

\begin{figure*}[htbp]
\centering
\subfloat[][Displacement $\bm{q}$]{%
	\label{fig:err_q_elasticity}%
\includegraphics[width=0.45\textwidth]{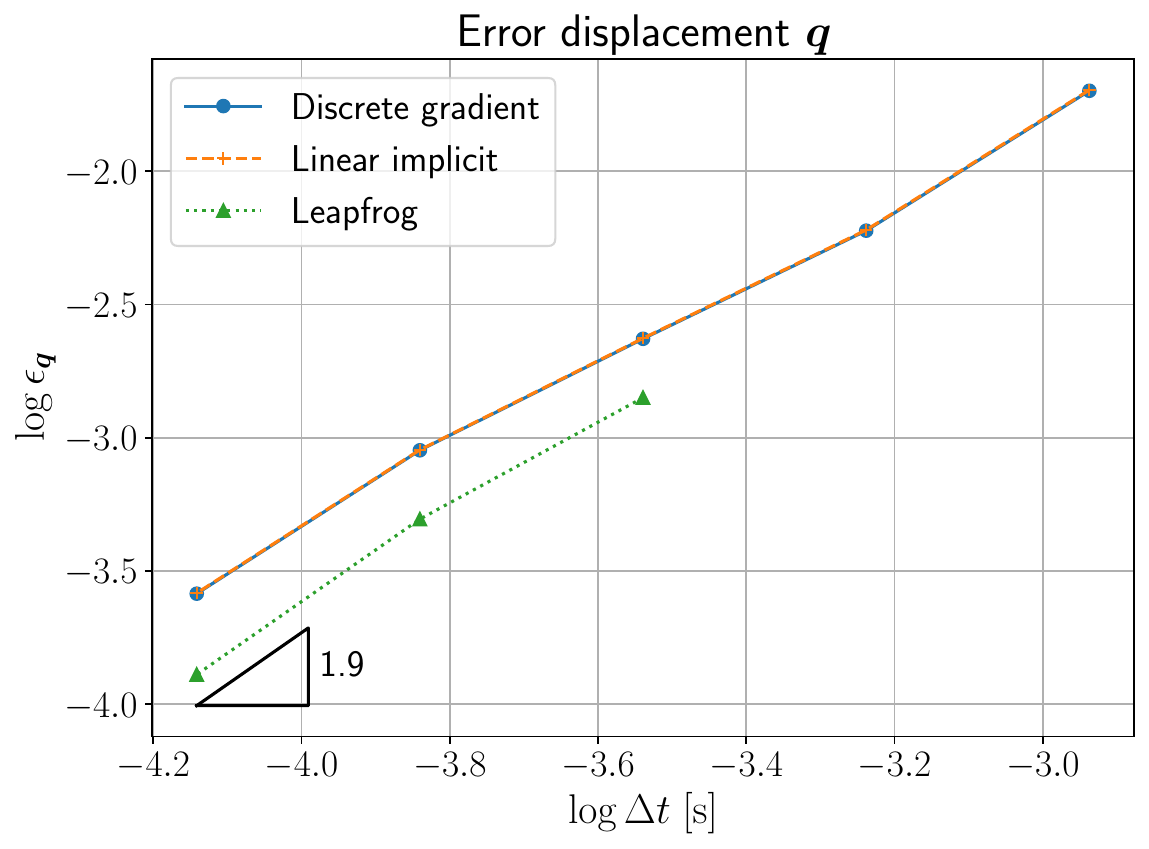}}%
 \hspace{8pt}%
\subfloat[][Velocity $\bm{v}$]{%
	\label{fig:err_v_elasticity}%
\includegraphics[width=0.45\textwidth]{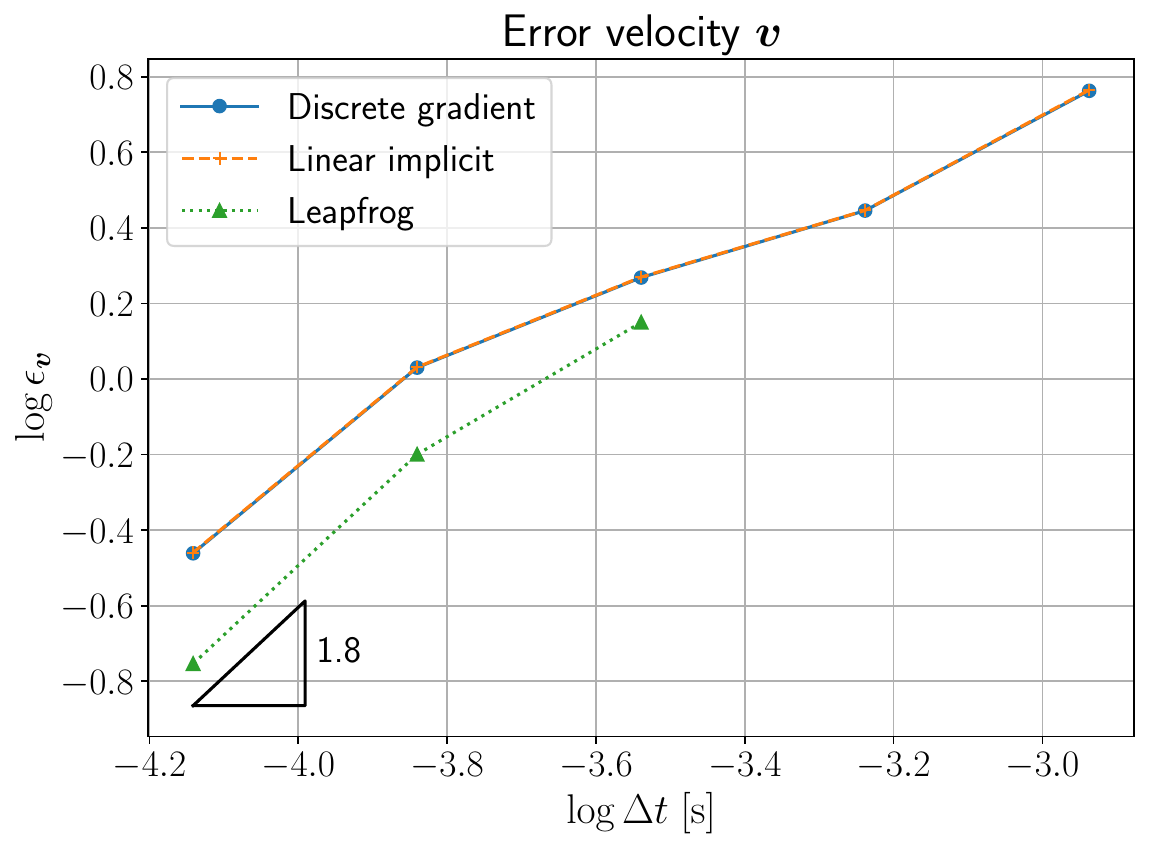}}
\caption{Convergence rate for $\bm{q}, \; \bm{v}$ in geometrically nonlinear elasticity}%
\label{fig:conv_elasticity}%
\end{figure*}

\begin{figure*}[htbp]
    \centering
    \begin{minipage}{0.45\textwidth}
    \includegraphics[width=0.98\textwidth]{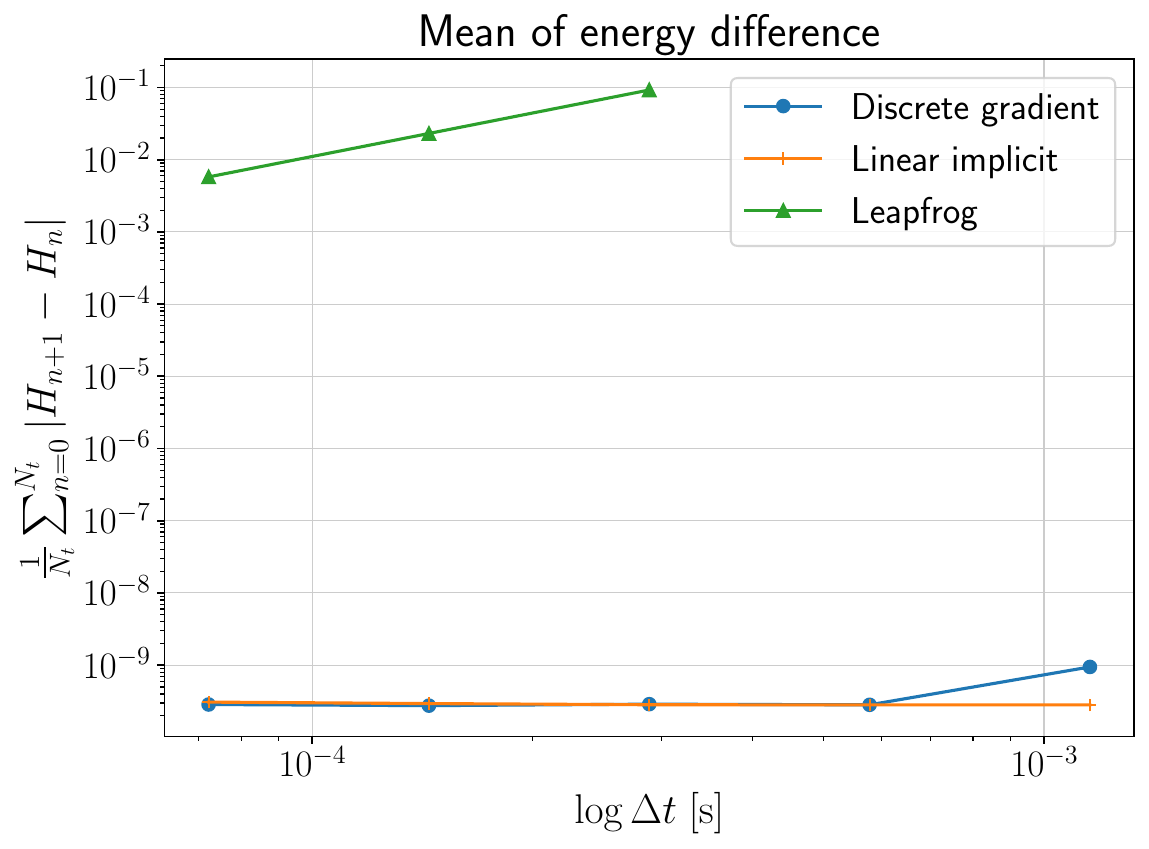}
    \caption{Mean of $|H_{n+1}-H_n|$ for the geometrically nonlinear elasticity problem}
    \label{fig:err_energy_elasticity}
    \end{minipage}\hspace{8pt}
    \begin{minipage}{0.45\textwidth}
    \includegraphics[width=0.98\textwidth]{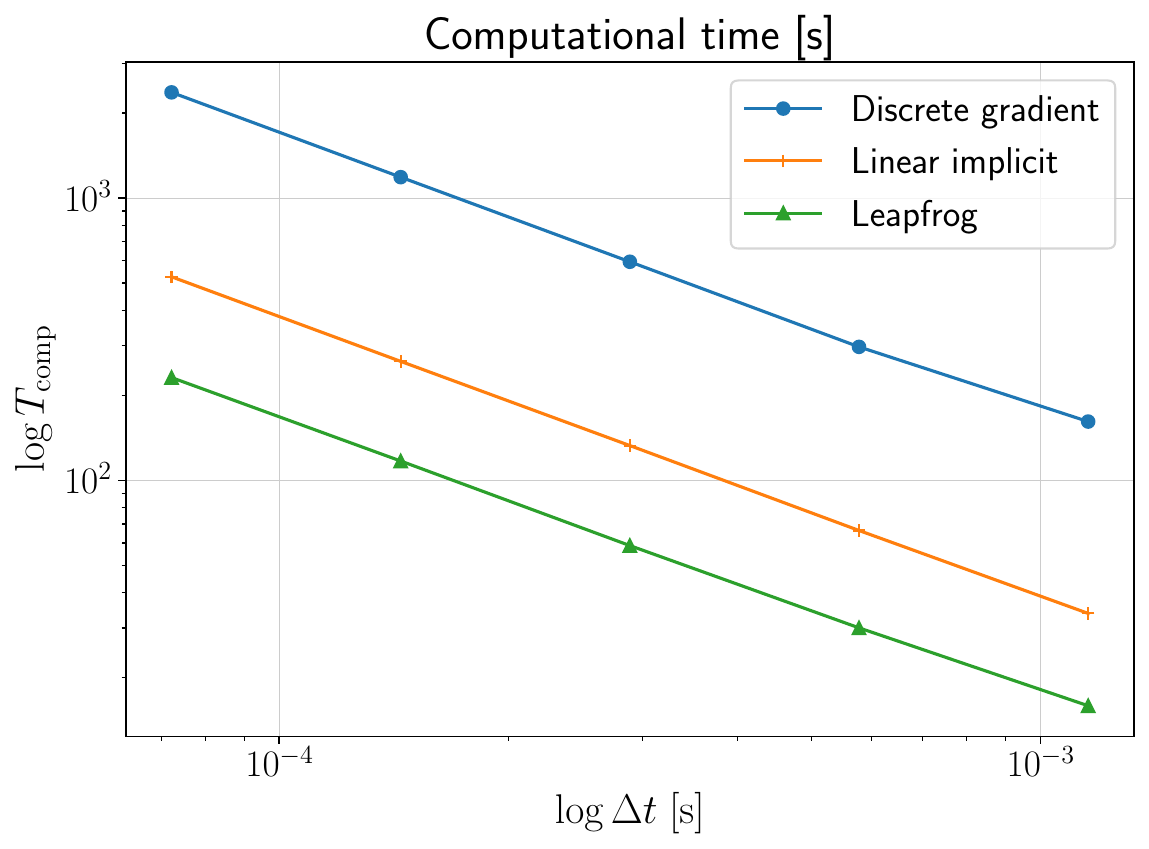}
    \caption{Computational time for geometrically nonlinear elasticity problem}
    \label{fig:comp_time_elasticity}
    \end{minipage}
\end{figure*}

Snapshots of the displacement solution at different time instants are reported in Fig. \ref{fig:screenshots_elasticity}. The snapshots are computed via the leapfrog method with a time step $\Delta t_{\rm ref} = \Delta t_{\rm base}/2^7$ (this is the reference solution considered for the convergence analysis). \rewtwo{The Frobenius norm of the stress, i.e. $||\bm{S}||_F$ using linear implicit scheme with $\Delta t = \Delta t_{\rm base}/2^5$ is reported in Fig. \ref{fig:screenshots_stress_elasticity}.} The convergence analysis reported in Fig. \ref{fig:conv_elasticity}. The accuracy for the discrete gradient and the linearly implicit method is again comparable and the leapfrog method is again unstable for the coarsest time step and requires $\Delta t <0.25 \Delta t_{\rm base}$ to be stable. The mean of the energy difference between adjacent time steps is reported in Fig. \ref{fig:err_energy_elasticity}. The energy conservation is verified for the linearly implicit method and discrete gradient, where the leapfrog is much less accurate in this respect. The computational time is reported in Fig. \ref{fig:comp_time_elasticity}. For this example the system to be solved is in the order of 5000 degrees of freedom and the leapfrog method is only 2 to 3 faster than the linear implicit method. The discrete gradient is one order of magnitude slower that the leapfrog method. 

\begin{figure*}[htbp]
    \centering
        \includegraphics[width=0.3\textwidth]{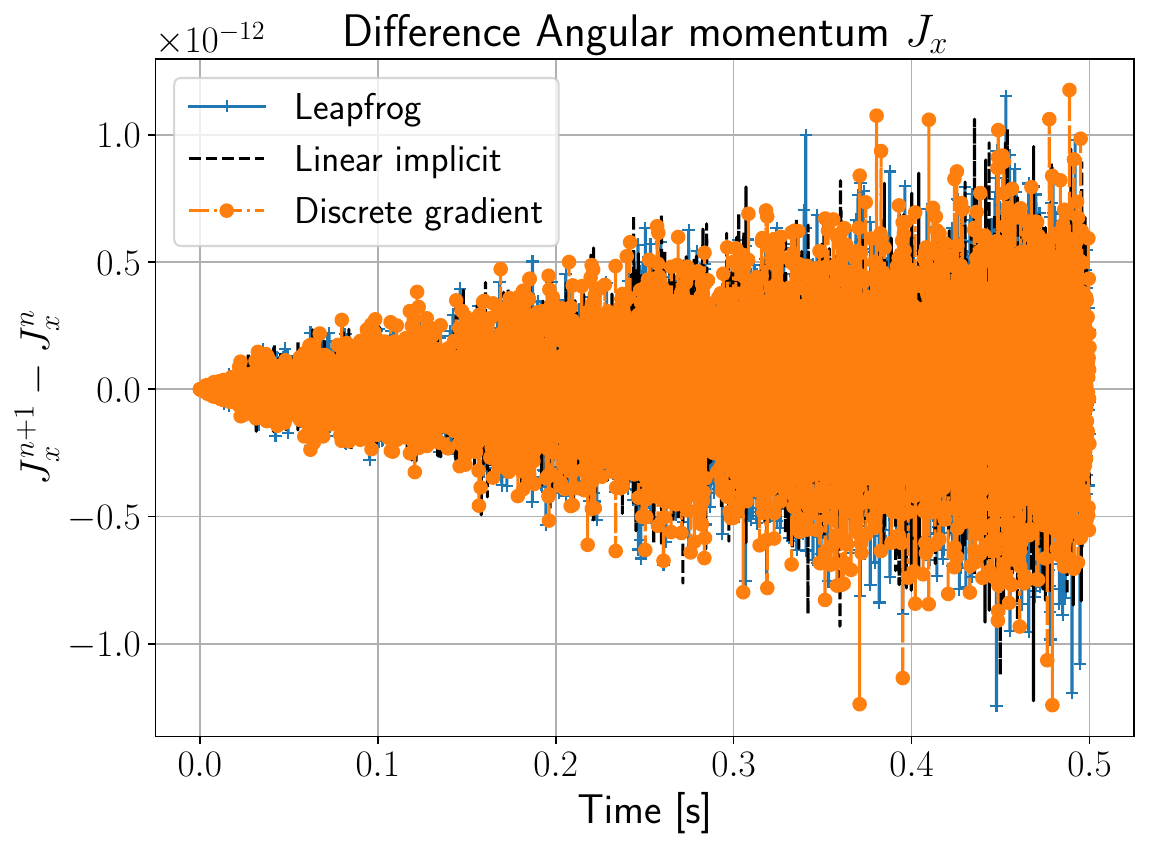}
    \hfill
        \includegraphics[width=0.3\textwidth]{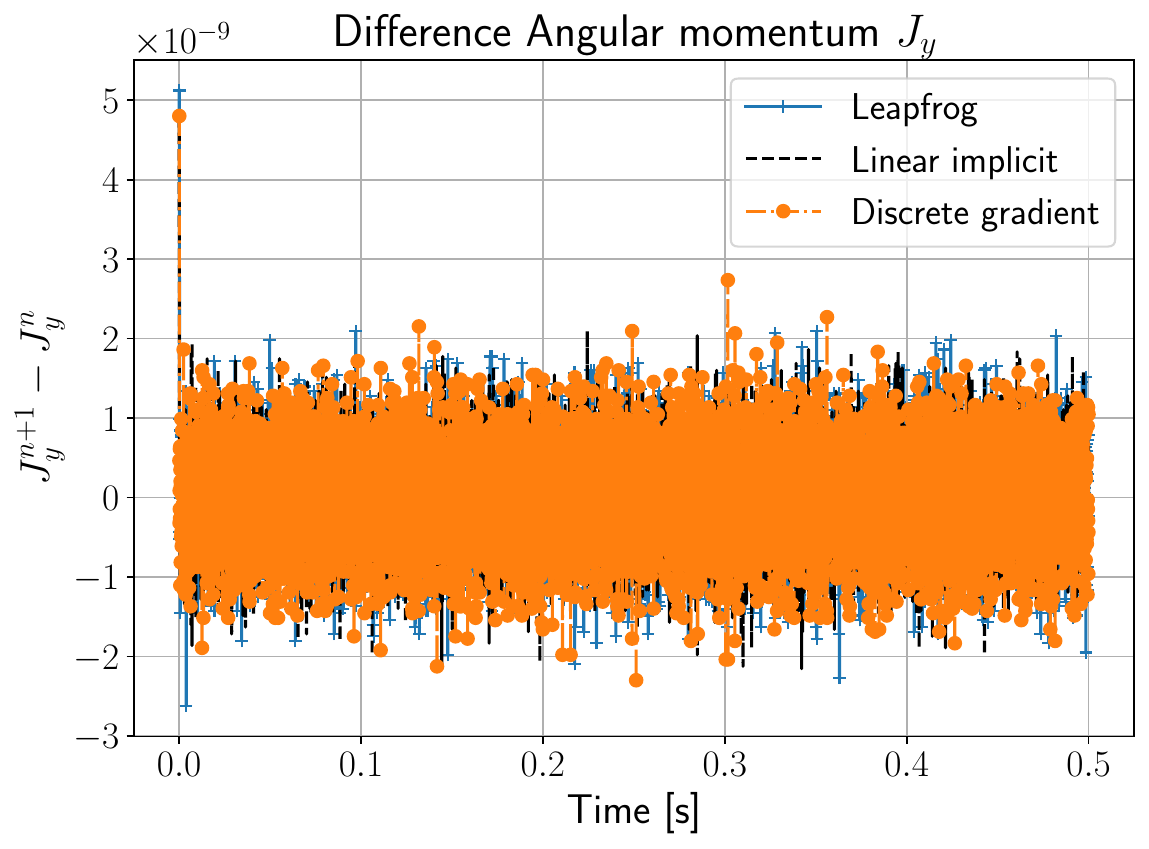}
    \hfill
        \includegraphics[width=0.3\textwidth]{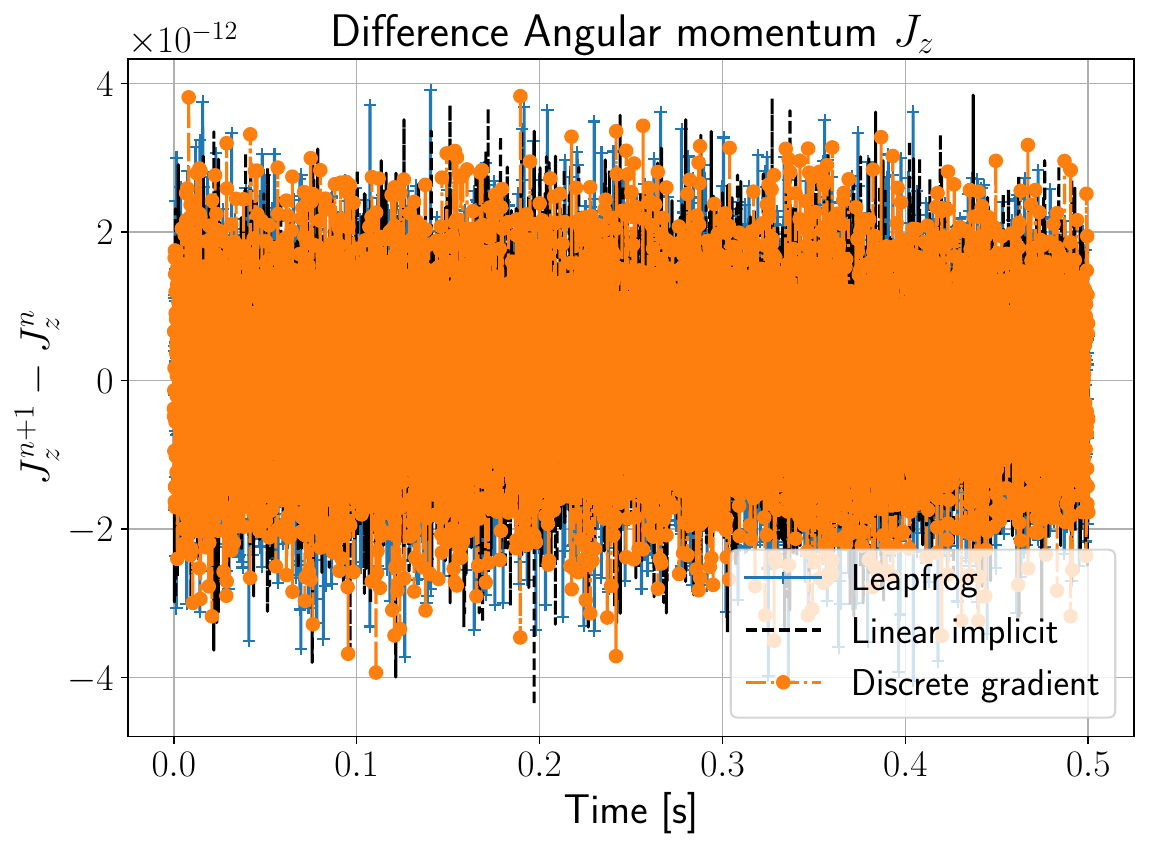}
    \caption{Angular momentum conservation for the bending column without Dirichlet boundary conditions.}
    \label{fig:angular_mom}
\end{figure*}

\rewtwo{To show the angular momentum conservation, we consider the same simulation without Dirichlet boundary conditions. The time trend of the angular momentum is shown in Fig. \ref{fig:angular_mom}. The leapfrog and energy-momentum preserving scheme are known to preserve the angular momentum \cite{hairer2003verlet,simo1992conserving} and so does the proposed integrator.}

\section{Conclusion}
\label{sec:concl}
In this contribution we have presented a general framework to recast geometrically nonlinear problems into a non-canonical Hamiltonian formulation. The procedure is detailed for moderate rotations of beams and plates described by von-K\'arm\'an nonlinearities and for geometrically nonlinear elasticity. The methodology is readily applicable to general geometrically nonlinear problems arising in continuum mechanics, like rods, beams, plates, shells and solid mechanics. The non-canonical Hamiltonian structure can be readily discretized by mixed finite elements where the stress field is an additional unknown of the problem. The finite element space used for its approximation needs to be rich enough in order to accurately describe the geometrically nonlinear effects. Since this variable is discretized using a discontinuous space it can be statically condensed at the discrete time level. The time discretization combines the St\"ormer Verlet and implicit midpoint methods to achieve an exact energy conservation without requiring the solution of a nonlinear system. It is shown that the time integration strategy is essentially analogous to the one proposed in \cite{bilbao2023explicit} in the context of the scalar auxiliary variable method. The method exhibits higher accuracy than both the leapfrog and discrete gradient methods when the error is measured against an analytical solution. However, additional tests are required to verify this assertion. The scheme is more efficient than the discrete gradient scheme as it does not require the solution of a nonlinear system and the linear system arising from discrete time is always solvable. As highlighted in the different examples, it is more stable than the leapfrog method, even though generally slower. The numerical results indicate that the methodology is a valid alternative to established time integration strategy in the context of nonlinear elasticity. Future developments may include linear algebra strategies to make the method computationally more efficient. Furthermore, a finite elements and finite differences coupling may be exploited to leverage the intrinsic advantages of the two methods.

\section*{Code Availability}
The code used for the present work is hosted at: \\
\url{https://github.com/a-brugnoli/hamiltonian-geometrically-nonlinear-elasticity/tree/main}

\section*{Acknowledgments}
The first author would like to thank Jeremy Bleyer for pointing out the possibility of statically condensed the stress variable. The authors would like to thank Michele Ducceschi for many insightful discussions on the topic. 

\bibliography{biblio}

\appendix

\section{Two-dimensional differential  operators}\label{app:diff_2D}
Let $\bbV= \bbR^2, \; \bbM = \bbR^{2\times 2}, \; \bbS = \mathrm{sym}(\bbM)$. The rotation of a vector $\bm{v}=(v_1 \; v_2)^\top$ is denoted by 

$$
\begin{pmatrix}
v_1 \\ v_2
\end{pmatrix}^\perp :=
\begin{pmatrix}
v_2 \\ -v_1
\end{pmatrix}, \quad 
\bm{v}^\perp = \bm{J} \bm{v}, \quad \bm{J}:=\begin{pmatrix}
    0 & 1 \\
    -1 & 0
\end{pmatrix}.
$$
The $\curl$ operator of a scalar function provides a vector given by
$$ 
\curl v = \begin{pmatrix}
 \partial_2 v \\ 
-\partial_1 v
\end{pmatrix}
= (\grad{v})^\perp = \bm{J}\grad v.
$$

For a vector $\bm{v}=(v_1 \; v_2)^\top$ the $\rot: C^\infty(\bbV) \rightarrow C^\infty(\bbR)$ operator reads
$$
\rot \bm{v} = \pdv{v_2}{x_1}-\pdv{v_1}{x_2} = \div \bm{v}^\perp =\div(\bm{J}\bm{v}).
$$
The $\curl$ and $\rot$ operators are adjoint operators i.e. for vanishing boundary conditions of either variable it holds
\begin{equation}\label{eq:curl_rot_vectors}
\begin{aligned}
    \inpr[\Omega]{\bm{u}}{\curl v} &= \int_{\Omega} \{u_1 \partial_2 v - u_2 \partial_1 v\} \d\Omega, \\
    &= \int_{\Omega} (-\partial_2 u_1 + \partial_1 u_2)  v\, \d\Omega, \\
    &= \inpr[\Omega]{\rot \bm{u}}{v}.
\end{aligned}
\end{equation}
When applied to a vector $\bm{v}=(v_1 \; v_2)^\top$ the curl operator gives a matrix
$$
\curl \bm{v}  = 
\begin{bmatrix}
    \partial_2 v_1 & -\partial_1 v_1 \\
    \partial_2 v_2 & -\partial_1 v_2 \\
\end{bmatrix}.
$$

The divergence of a matrix is defined row-wise as 
$$
\div{\bm{M}} = \begin{pmatrix}
\partial_1 M_{11} + \partial_2 M_{12}\\ 
\partial_1 M_{21} + \partial_2 M_{22} \\ 
\end{pmatrix}.
$$
The $\rot$ operator applied to a matrix  is defined by
$$
\rot \bm{M} = \begin{pmatrix}
    \partial_1 M_{12} -\partial_2 M_{11} \\ 
    \partial_1 M_{22} -\partial_2 M_{21} \\ 
\end{pmatrix}.
$$
Thus the $\rot$ and $\div$ operators for matrices are related by $\rot\bm{M} = \div(\bm{M}\bm{J}^\top)$. 
Not surprisingly, the $\curl$ operator for vectors and the $\rot$ operator for matrices are adjoint operators, i.e. for vanishing boundary conditions
\begin{equation}\label{eq:curl_rot_matrices}
\inpr[\Omega]{\bm{M}}{\curl \bm{v}} = \inpr[\Omega]{\rot \bm{M}}{\bm{v}}.
\end{equation}
The $\Air: C^\infty(\bbR) \rightarrow C^\infty(\bbS)$ operator (applied to a scalar function) is then defined as
$$
\Air v = \curl \curl v = \begin{bmatrix}
    \partial_{22}v & -\partial_{12} v \\
    -\partial_{12}v & \partial_{11} v
\end{bmatrix}.
$$
Given $\bm{S} \in C^{\infty}(\bbS)$, and from Eqs. \eqref{eq:curl_rot_vectors} and \eqref{eq:curl_rot_matrices} the adjoint of the $\Air= \curl\curl$ operator can then be obtained as $\mathrm{Air}^* = \rot\rot$
$$
\begin{aligned}
\inpr[\Omega]{\bm{S}}{\Air v} &= \inpr[\Omega]{\bm{S}}{\curl\curl v} \\
&=\inpr[\Omega]{\rot \bm{S}}{\curl v} \\
&=\inpr[\Omega]{\rot\rot\bm{S}}{v}.    
\end{aligned}
$$

The $\rot\rot$ operator can be seen as a rotated double divergence
\begin{equation}\label{eq:rotrot_divdiv}
\rot\rot\bm{S} = \div\div(\bm{JSJ}^\top).
\end{equation}

\section{Elasticity complex in \texorpdfstring{$\mathbb{R}^2$}{R}}\label{app:complex}
A complex is a sequence of vector spaces connected by differential operators such that the composition of two consecutive operators vanishes. In two dimension, the elasticity complex simplifies into two complex. The first one is the Airy complex
\begin{equation}\label{eq:air_complex}
\begin{tikzcd}[column sep=normal]
\bbP_1 \arrow[r, "\subset"] & C^\infty\arrow[r, "\Air"] & C^\infty \otimes \mathbb{S}  \arrow[r, "\div"] &  C^\infty \otimes \mathbb{V} 
\end{tikzcd}    
\end{equation}
where $\bbP_1:=\bbR + \bm{x}^\perp \cdot \bbR^2$ is the space of first order polynomial.
The second complex, i.e. the adjoint of the Airy complex, is given by the $\rot\rot$ complex
\begin{equation}\label{eq:rotrot_complex}
    \begin{tikzcd}[column sep=normal]
\mathbf{RM} \arrow[r, "\subset"] & C^\infty \otimes \mathbb{V} \arrow[r, "\mathrm{def}"] & C^\infty \otimes \mathbb{S}  \arrow[r, "\rot\rot"] &  C^\infty
\end{tikzcd}
\end{equation}
where $\mathbf{RM}=\bbR^2 + \bm{x}^\perp \bbR$ is the space of rigid body motion. The $\rot\rot$ is the adjoint of the $\Air$ operator. The key property of a complex in this case reads
\begin{equation*}
    \div\Air \equiv 0, \qquad \rot\rot\sgrad \equiv 0.
\end{equation*}

\end{document}